\documentclass[a4paper,11pt]{article}
\usepackage{ifstyle}   
\usepackage{graphicx}
\usepackage{epsfig} 
\numberwithin{equation}{section}  

\usepackage[francais]{babel}               
\usepackage[T1]{fontenc}
\usepackage[latin1]{inputenc}  
\usepackage{amssymb,diagrams}             
\usepackage{array}                 
\usepackage{fullpage}
\newtheorem{theoreme}{Théorème} 
\newtheorem{corol}{Corollaire}
\newtheorem{lem}[subsubsection]{Lemme}

\newcommand{\N}{\mathbb N}
\newcommand{\Z}{\mathbb Z}  

\newcommand{\R}{\mathbb R}
\newcommand{\C}{\mathbb C}
\newcommand{\K}{\mathbb K}
\newcommand{\I}{\mathbb I}

\begin{document} 
{\def\thefootnote{\relax}
\footnote{\hskip-0.6cm
{\bf{Mots-cl\'es}} : faisceaux analytiques cohérents, algèbre homologique non-commutative, EDP quasi-linéaires, methode Nash-Moser. \\
{\bf{Classification AMS}} : 32C35, 35N10.}}       
\begin{center} 
\Huge{\bf{Faisceaux $\bar{\partial}$-cohérents sur les variétés complexes}}
\\
\vspace{0.6cm}
\huge{Nefton Pali}
\\
\Large{Institut Fourier, UMR 5582, Université Joseph Fourier
\\
BP 74, 38402 St-Martin-d'Hères cedex, France
\\
E-mail: \textit{nefton.pali@ujf-grenoble.fr}}
\end{center} 
\vspace{\fill} 
{\bf Résumé}.-Cet article est une version, destinée à une audience plus large, de notre récent travail (\cite{pal}). Dans ces travaux nous généralisons au contexte des faisceaux analytiques cohérents un résultat classique de Grothendieck-Koszul-Malgrange (\cite{ko-mal}) concernant l'intégrabilité des connexions de type $(0,1)$ sur un fibré vectoriel ${\cal C}^{\infty}$ au dessus d'une variété complexe. En introduisant la notion de faisceau $\bar{\partial}$-cohérent, qui est une notion qui vit dans le contexte ${\cal C}^{\infty}$, nous montrons l'existence d'une équivalence (exacte) entre la catégorie des faisceaux analytiques cohérents et la catégorie des faisceaux $\bar{\partial}$-cohérents. La difficulté essentielle de la preuve de ce résultat consiste à résoudre un problème différentiel quasi-linéaire dont le terme principal est un opérateur $\bar{\partial}$ usuel. La solution de ce problème est obtenue en utilisant un schéma de convergence rapide de type Nash-Moser. Un autre ingrédient important pour la conclusion de la preuve est un résultat profond de Malgrange qui affirme la fidélité plate de l'anneau des germes des fonctions ${\cal C}^{\infty}$ à valeurs complexes en un point, sur l'anneau des germes des fonctions holomorphes en ce point.
\\
\\
\\
{\bf{Abstract}}.-This paper is a large audience version of our recent work (\cite{pal}) in which we give a generalization, in the context of coherent analytic sheaves, of a classical result of Grothendieck-Koszul-Malgrange (\cite{ko-mal}) concerning the integrability of connections of type $(0,1)$ over a ${\cal C}^{\infty}$ vector bundle over a complex manifold. We introduce the notion of $\bar{\partial}$-coherent sheaf, which is a ${\cal C}^{\infty}$ notion, and we prove the existence of an (exact) equivalence between the category of coherent analytic sheaves and the category of $\bar{\partial}$-coherent sheaves. The principal difficulty of the proof is the solution of a quasi-linear differential equation with standard $\bar{\partial}$ as its principal term. We are able to find a solution of this differential equation, using a rapidly convergent iteration scheme of Nash-Moser type. We also use a deep result of Malgrange asserting that the ring of germs of complex differentiable functions at a point is faithfully flat over the ring of germs of holomorphic functions at the same point.
\newpage
\tableofcontents
\addtocounter{section}{-1}
\newpage
\section{Introduction }
Le résultat de Grothendieck-Koszul-Malgrange mentionné dans le résumé (voir \cite{ko-mal}) affirme qu'un fibré vectoriel complexe différentiable au dessus d'une variété complexe, qui admet une connexion $\bar{\partial}$ de type $(0,1)$ telle que $\bar{\partial}^2=0$ possède une structure de fibré vectoriel holomorphe. En d'autres termes, en utilisant l'équivalence entre les notions de fibré vectoriel et de faisceau localement libre sur une variété, on observe que le noyau de la connexion  $\bar{\partial}$ est un faisceau de ${\cal O}$-modules localement libre. Le point essentiel de ce résultat consiste à prouver l'existence d'une solution de l'équation différentielle quasi-linéaire  $g^{-1}\bar{\partial}_{_J }g=A$ (où $\bar{\partial}_{_J }$ est la (0,1)-connexion canonique sur le faisceau des fonctions ${\cal C}^{\infty}$ et $A$ représente la $(0,1)$-forme de connexion locale de $\bar{\partial}$) avec la condition d'intégrabilité  $\bar{\partial}_{_J }A +A\wedge A=0$. Nous généraliserons cette équation pour prouver notre caractérisation différentielle laquelle introduit comme objet nouveau la notion de faisceau $\bar{\partial}$-cohérent. Précisément un faisceau  $\bar{\partial}$-cohérent est un faisceau ${\cal G}$ de modules de fonctions ${\cal C}^{\infty}$ à valeurs complexes, muni d'une connexion $\bar{\partial}$ de type $(0,1)$ telle que $\bar{\partial}^2=0$, et qui admet localement  des résolutions de longueur finie, par des  modules de fonctions ${\cal C}^{\infty}$ à valeurs complexes. Notre caractérisation affirme essentiellement que le noyau de la connexion $\bar{\partial}$ est un faisceau analytique cohérent. On aura alors une équivalence exacte (l'exactitude est due à la fidélité plate de l'anneau des germes des fonctions ${\cal C}^{\infty}$ à valeurs complexes sur l'anneau des germes des fonctions holomorphes, voir \cite{mal}) entre la catégorie des faisceaux analytiques cohérents et la catégorie des faisceaux  $\bar{\partial}$-cohérents. La difficulté essentielle de la preuve consiste à montrer que quel que soit le choix de la résolution locale du faisceau ${\cal G}$, on peut trouver au voisinage de chaque point de l'ouvert sur lequel on considère la résolution locale, une autre résolution locale constituée de matrices holomorphes. En d'autres termes on cherche une résolution locale  admettant des formes de connexion nulles. Pour atteindre cet objectif on introduit la notion de recalibration, laquelle généralise la notion classique de changement de jauge pour les formes locales d'une connexion de type $(0,1)$ intégrable sur un faisceau localement libre. La notion de recalibration ne représente rien d'autre qu'une action d'un semi-groupe sur l'ensemble des formes qui représentent localement la condition d'intégrabilité de la connexion $\bar{\partial}$. La notion en question permet de traduire notre problème en termes d'un système différentiel quasi-linéaire dont le terme principal est un opérateur $\bar{\partial}$ usuel. Les conditions d'intégrabilité de ce système ne sont rien d'autre que les expressions locales de la condition d'intégrabilité de la connexion $\bar{\partial}$. Les solutions du système différentiel seront obtenues à l'aide d'un procédé itératif de type Nash-Moser, dont chaque étape  est déterminée par une recalibration obtenue en fonction de l'étape précédente. La preuve de l'existence de solutions du système différentiel en question est exposée dans la sous-section $\ref{s3.4}$, et constitue la partie essentielle de la preuve de notre caractérisation des faisceaux analytiques cohérents. La technique qui consiste à utiliser des schémas itératifs pour montrer l'existence des solutions de problèmes non linéaires est désormais bien connue en analyse complexe. On peut citer par exemple les travaux de Webster \cite{We-1} et \cite{We-2}, qui utilisent des techniques de type Nash-Moser, (voir \cite{Mos-1} et \cite{Mos-2}) pour prouver l'existence des solutions de deux problèmes différentiels fondamentaux en géométrie complexe. L'ingrédient final qui permet de conclure notre preuve est le  résultat profond de Malgrange mentionné précédemment, lequel permet aussi une généralisation du théorème de Dolbeault au cas des faisceaux analytiques cohérents ($\bar{\partial}$-cohérents).  Enfin on remarque aussi un résultat d'intégrabilité pour les connexions sur les faisceaux admettant des résolutions locales de longueur finie sur les variétés différentiables.
\newpage
\section{Faisceaux $\bar{\partial}$-cohérents sur les variétés complexes}\label{S1} 
On va commencer par rappeler un résultat classique bien connu du à Grothendieck-Koszul-Malgrange. On a besoin d'abord de quelques notations et remarques. Soit $(X,J)$ une variété complexe, où $J$ est le tenseur de la structure presque-complexe supposé  intégrable. On désigne par ${\cal E}_X$ le faisceau des fonctions ${\cal C}^{\infty} $ à valeurs complexes, par  ${\cal E}^{0,q}_X$ le faisceau des $(0,q)$-formes et par
$\bar{\partial}_{_J }\in {\cal H} om_{_{{\cal O}_X } }({\cal E}^{0,q}_X, {\cal E}^{0,q+1}_X )(X)$ la composante de type $(0,1)$ de la différentielle. Soit $F\longrightarrow X$ un fibré vectoriel complexe ${\cal C}^{\infty}$. On rappelle que sur une variété complexe il y a une équivalence entre les notions suivantes.
\begin{eqnarray*} 
\{\mbox{Fibrés vectoriels complexes ${\cal C}^{\infty}$ } \} &\longleftrightarrow &\{\mbox{Faisceaux de ${\cal E}$-modules localement libres} \}
\\
\\
F &\longmapsto &{\cal E}(F):=\mbox{Faisceau des sections ${\cal C}^{\infty}$ de $F$ }  
\end{eqnarray*} 
L'application inverse envoie ${\cal G}\mapsto (F,(\bar{\psi}_{\alpha})_{\alpha })$
où
$$
F:=\coprod_{x\in X}  {\cal G}_x/m({\cal E}_x)\cdot {\cal G}_x
$$
et les trivialisations locales 
$$
\bar{\psi}_{\alpha}:U_{\alpha} \times \C^r\stackrel{\simeq}{\longrightarrow} F_{|_{U_{\alpha} } }
$$ 
sont obtenues de façon naturelle à partir des trivialisations locales 
$\psi_{\alpha }: {\cal E}_{_{U_{\alpha } }  }^{\oplus r }\stackrel{\simeq}{\longrightarrow}
{\cal G}_{{|_{U_{\alpha } } } } $ du faisceau localement libre ${\cal G} $. L'équivalence précédente est aussi valable pour les fibrés vectoriels holomorphes et les faisceaux de ${\cal O}$-modules localement libres définis sur une variété complexe. On considère la définition suivante.
\begin{defi}
Une connexion de type $(0,1)$ sur un faisceau $\cal G$ de
${\cal E}_X$-modules est un morphisme de faisceaux de groupes additifs $\bar{\partial}:{\cal G}\longrightarrow{\cal G} \otimes_{_{{\cal E}_X }}{\cal E}^{0,1}_X$
 tel 
que $\bar{\partial}(g\cdot f)=( \bar{\partial} g)\cdot f +g\otimes\bar{\partial}_{_J } f$, pour tout $g\in {\cal G}_x$ et $f\in {\cal E}_{X,x}$.
\end{defi}
La donnée d'une connexion de type $(0,1)$ sur le faisceau de ${\cal E}_X$-modules ${\cal G}$ détermine de façon univoque une dérivation $\bar{\partial}$ de type $(0,1)$ sur le complexe $({\cal G} \otimes_{_{{\cal E}_X }}{\cal E}^{0,q}_X)_{ q \geq 0}$. En effet on peut définir l'extension 
$
\bar{\partial}:{\cal G}\otimes_{_{{\cal E}_X }}{\cal E}^{0,q}_X\longrightarrow  {\cal G} \otimes_{_{{\cal E}_X }}{\cal E}^{0,q+1}_X 
$
 par la formule classique
\begin{eqnarray*} 
&(\bar{\partial}\omega )(\xi _0,...,\xi _q):=\displaystyle{\sum_{0\leq j \leq q}(-1)^j }\bar{\partial}(\omega (\xi _0,...,\widehat{\xi _j},..., \xi _q))(\xi _j)+&
\\
\\
&+\displaystyle{\sum_{0\leq j<k \leq q}(-1)^{j+k}}\omega ([\xi _j,\xi _k],\xi _0,...,\widehat{\xi _j},...,\widehat{\xi _k},..., \xi _q)& 
\end{eqnarray*} 
avec $\omega\in ({\cal G} \otimes_{_{{\cal E}_X }}{\cal E}^{0,q}_X)(U)$ et $\xi _j\in{\cal E}(T^{0,1}_X)(U)$, sur un ouvert $U$ quelconque, ou de façon équivalente, par la règle de Leibnitz 
$$
\bar{\partial}(g\otimes\alpha ):=\bar{\partial}g \wedge \alpha +g\otimes\bar{\partial}_{_J }\alpha  
$$
avec $g\in{\cal G}_x$ et $\alpha \in {\cal E}^{0,q}_{X,x} $ pour tout $x\in X$. Souvent on pense une connexion en termes de sa extension au complexe $({\cal G} \otimes_{_{{\cal E}_X }}{\cal E}^{0,q}_X)_{ q \geq 0}$. Sur tout faisceau $\cal F$ de ${\cal O}_X$-modules on peut considérer
la connexion canonique 
$$\bar{\partial}_{_{\cal F} }:=\I_{_{\cal F}} \otimes _{_{{\cal O}_X }}\bar{\partial}_{_J }
:{\cal F}\otimes_{_{{\cal O}_X }}{\cal E}^{0,q}_X 
\longrightarrow  {\cal F} \otimes_{_{{\cal O}_X }}{\cal E}^{0,q+1}_X $$
qui est vue comme une connexion de type $(0,1)$ sur le faisceau de  ${\cal E}_X$-modules 
$$
{\cal F}^{\infty}:={\cal F} \otimes_{_{{\cal O}_X }}{\cal E}_X .
$$ 
Bien évidement la définition précédente est une généralisation immédiate de la notion classique de connexion de type $(0,1)$ canonique associée à un fibré vectoriel holomorphe, (voir par exemple les ouvrages \cite{Dem}, (\cite{Gri-Ha} et(\cite{Wel}).
La donnée d'une connexion de type $(0,1)$ sur ${\cal G}$ détermine aussi le tenseur de 
courbure de la connexion qu'on 
notera par
$$
\Theta _{\bar{\partial}}\in \Big({\cal E}nd_{_{{\cal E}_{_X}}}({\cal G}) \otimes_{_{{\cal E}_{_X} }}{\cal E}^{0,2}_X\Big)(X)
$$ 
et qu'on définit par la formule $\Theta _{\bar{\partial}}(\xi ,\eta)\cdot g:=(\bar{\partial}^2g)(\xi ,\eta)$ 
avec $g\in{\cal G}(U)$ et $\xi ,\eta\in{\cal E}(T^{0,1}_X)(U)$. On note de plus $\xi _{\bar{\partial} }\cdot g:=\bar{\partial}g(\xi )$ la dérivée covariante de la section $g$ calculée le long du champ $\xi $ et on remarque que la première définition de l'extension de la connexion $\bar{\partial}$ implique de façon triviale la formule
$$
\xi _{\bar{\partial} }\,.\,(\eta _{\bar{\partial} }\,.\, g)- \eta _{\bar{\partial} }\,.\,(\xi _{\bar{\partial} }\,.\,g )=[\xi ,\eta]_{\bar{\partial} }\,.\,g +\Theta _{\bar{\partial}}(\xi ,\eta)\,.\,g ,
$$
où $[\xi ,\eta]\in{\cal E}(T^{0,1}_X)(U)$, grâce à l'hypothèse d'intégrabilité du tenseur de la structure presque-complexe $J\in {\cal C}^{\infty} (T_X^*\otimes_{_{_\R}} T_X)(X)$. Le tenseur de courbure $\Theta _{\bar{\partial}}$ exprime donc le défaut de commutation des dérivées covariantes secondes des sections de ${\cal G}$ le long des champs de type $(0,1)$. Nous porterons un intérêt particulier aux connexions de type $(0,1)$ intégrables, c'est à dire aux connexions telles que $\bar{\partial}^2=0$. La formule précédente caractérise alors ce type de connexions $(0,1)$ comme étant celles pour lesquelles les dérivées covariantes secondes, calculées le long de deux champs qui commutent, commutent également. En termes explicites on a l'égalité $\xi _{\bar{\partial} }\cdot(\eta _{\bar{\partial} }\cdot g)= \eta _{\bar{\partial} }\cdot(\xi _{\bar{\partial} }\cdot g )$ si $[\xi ,\eta]=0$. Un exemple de $(0,1)$-connexion intégrable est évidemment la connexion $\bar{\partial}_{_{\cal F}}$ introduite précédemment.
Avec les notations introduites précédemment on peut  énoncer le résultat de Grothendieck-Koszul-Malgrange (\cite{ko-mal}) sous la forme suivante.
\begin{theoreme}\label{Gro0} 
Soit $F\longrightarrow X$ un fibré vectoriel complexe ${\cal C}^{\infty}$ sur une variété complexe $X$. Alors l'existence d'une structure holomorphe sur le fibré $F$ est équivalente à l'existence d'une connexion $\bar{\partial}:{\cal E}(F)\longrightarrow  {\cal E}(F) \otimes_{_{{\cal E}_X }}{\cal E}^{0,1}_X$
de type $(0,1)$ intégrable $(i.e.\; \bar{\partial}^2=0)$ sur le faisceau de ${\cal E} $-modules ${\cal E}(F)$.
\end{theoreme} 
En utilisant l'équivalence entre les notions de fibrés vectoriels holomorphes et faisceaux de ${\cal O}$-modules localement libres sur une variété complexe, on peut reformuler en termes équivalents le théorème  $\ref{Gro0}$ sous la forme suivante.
\begin{theoreme}\label{Gro1}
Soit $F\longrightarrow X$ un fibré vectoriel complexe ${\cal C}^{\infty}$ sur une variété complexe $X$ muni d'une connexion $\bar{\partial}:{\cal E}(F)\longrightarrow  {\cal E}(F) \otimes_{_{{\cal E}_X }}{\cal E}^{0,1}_X$ telle que $\bar{\partial}^2=0$. Alors le noyau $Ker \bar{\partial}\subset  {\cal E}(F)$ de la connexion est un faisceau de ${\cal O}$-modules localement libre tel que 
$(Ker\bar{\partial})\cdot {\cal E}_X={\cal E}(F)$ $($ceci signifie que les générateurs locaux du  noyau $Ker \bar{\partial}$ sur le faisceau ${\cal O}_X $ sont aussi des générateurs locaux de ${\cal E}(F)$ sur le faisceau ${\cal E}_X)$.
\end{theoreme}
On a en conclusion que le noyau de la connexion $\bar{\partial}$ est le faisceau des sections holomorphes ${\cal O}(F)$ du fibré $F$. Dans le cas des faisceaux de ${\cal E}_X$-modules inversibles  qui admettent une connexion $\bar{\partial}_0$ de type $(0,1)$ telle que $\bar{\partial}^{2}_0=0$ on sait que toutes les connexions de ce type, et seulement celles ci, sont de la forme  $\bar{\partial}_0+A\otimes$ où $A\in{\cal E}^{0,1}_X(X)$ est une $(0,1)$-forme $\bar{\partial}_{_J }$-fermée. On a alors que si $L$ est un fibré en droites holomorphe les structures holomorphes sur $L$ sont en bijection avec les $(0,1)$-formes globales $\bar{\partial}_{_J }$-fermées. Avant d'énoncer le résultat qu'on se propose de démontrer, on remarque que si $\cal F$ est un faisceau analytique cohérent, le théorème des syzygies (voir \cite{kob}, chapitre 5)  implique l'existence d'une ${\cal O}$-résolution de longueur finie 
$$
0\rightarrow {\cal O}_{_U}^{\oplus p_m }\stackrel{\varphi_m }{\longrightarrow}  
{\cal O}_{_U}^{\oplus p_{m-1} }\stackrel{ \varphi _{m-1} }{\longrightarrow} \cdot \cdot\cdot
\stackrel{\varphi_2}{ \longrightarrow } {\cal O}_{_U}^{\oplus p_1 }\stackrel{\varphi_1}{ \longrightarrow }{\cal O}_{_U}^{\oplus p_0 }\stackrel{\psi}{\longrightarrow}
{\cal F}_{_{|_U} }\rightarrow 0
$$
dans la catégorie des faisceaux $\C$-analytiques cohérents. En rappelant que ${\cal F}^{\infty}:={\cal F} \otimes_{_{{\cal O}_X }}{\cal E}_X$ on obtient un diagramme commutatif suivant dont toutes les directions horizontales et verticales sont exactes.    
\begin{diagram}[height=1cm,width=1cm]
&         &0&                 &0&                                         &0&                                             &0&                                    \\
&         &\uTo&              &\uTo&                                      &\uTo&                                          &\uTo&                                 \\
0&  \rTo  & {\cal F}
          _{|_U}&  \rTo   &{\cal F}^{\infty}_{|_U}&    \rTo^{\bar{\partial} _{_{\cal F}}} 
 & {\cal F}^{\infty}_{|_U}\otimes_{_{{\cal E}_{_U}  }}{\cal E}^{0,1}_{_U} &   \rTo^{\bar{\partial}_{_{\cal F}}} & {\cal F}^{\infty}_{|_U}\otimes_{_{{\cal E}_{_U}   }}{\cal E}^{0,2}_{_U} &\rTo& \cdot \cdot \cdot
\\
&        &\uTo_{\psi} &              &\uTo_{\psi} &                                      &\uTo_{\psi \otimes \I_{(0,1)} } &                                          &\uTo_{\psi \otimes \I_{(0,2)} } &                                 \\
0&\rTo  & {\cal O}
          _{_U}^{\oplus p_0} &  \rTo   &{\cal E}^{\oplus p_0}_{_U}&    \rTo^{\bar{\partial}_{_J }}
 & ({\cal E}^{0,1}_{_U} )^{\oplus p_0}&   \rTo^{\bar{\partial}_{_J }} & ({\cal E}^{0,2}_{_U} )^{\oplus p_0}&\rTo&\cdot \cdot \cdot
\\
&        &\uTo_{\varphi _1}  &              &\uTo_{\varphi _1} &                                      &\uTo_{\varphi _1\otimes \I_{(0,1)} } &                                          &\uTo_{\varphi _1\otimes \I_{(0,2)} } &                                 \\
0&\rTo  & {\cal O}
          _{_U}^{\oplus p_1} &  \rTo   &{\cal E}^{\oplus p_1}_{_U}&    \rTo^{\bar{\partial}_{_J }}
 & ({\cal E}^{0,1}_{_U} )^{\oplus p_1}&   \rTo^{\bar{\partial}_{_J }} & ({\cal E}^{0,2}_{_U} )^{\oplus p_1}&\rTo&\cdot \cdot \cdot
\\
&        &\uTo_{\varphi _2}  &              &\uTo_{\varphi _2} &                                      &\uTo_{\varphi _2\otimes \I_{(0,1)} } &                                          &\uTo_{\varphi _2\otimes \I_{(0,2)} } &                                                  
\\
&         &\vdots&                 & \vdots&                                         &\vdots&                                             &\vdots&             
\\
&        &\uTo_{\varphi_{m-1} }  &              &\uTo_{\varphi_{m-1} } &                                      &\uTo_{\varphi_{m-1} \otimes \I_{(0,1)} } &                                          &\uTo_{\varphi_{m-1} \otimes \I_{(0,2)} } &                                                  
\\
0&\rTo  & {\cal O}
          _{_U}^{\oplus p_{m-1} } &  \rTo   &{\cal E}^{\oplus p_ {m-1}}_{_U}&    \rTo^{\bar{\partial}_{_J }}
 & ({\cal E}^{0,1}_{_U} )^{\oplus p_{m-1} }&   \rTo^{\bar{\partial}_{_J }} & ({\cal E}^{0,2}_{_U} )^{\oplus p_{m-1} }&\rTo&\cdot \cdot \cdot
\\
&        &\uTo_{\varphi _m}  &              &\uTo_{\varphi _m} &                                      &\uTo_{\varphi _m\otimes \I_{(0,1)} } &                                          &\uTo_{\varphi _m\otimes \I_{(0,2)} } &                                 \\
0&\rTo  & {\cal O}
          _{_U}^{\oplus p_m} &  \rTo   &{\cal E}^{\oplus p_m}_{_U}&    \rTo^{\bar{\partial}_{_J }}
 & ({\cal E}^{0,1}_{_U} )^{\oplus p_m}&   \rTo^{\bar{\partial}_{_J }} & ({\cal E}^{0,2}_{_U} )^{\oplus p_m}&\rTo&\cdot \cdot \cdot
\\
&        &\uTo&              &\uTo&                                      &\uTo&                                          &\uTo&                                 \\
&         &0&                 &0&                                         &0&                                             &0&                                    \\       
\end{diagram}
\\ 
La raison de l'exactitude est la suivante. La platitude de l'anneau ${\cal E}_{_{X,x} }$ sur l'anneau ${\cal O}_{_{X,x} }$ (voir l'ouvrage de Malgrange \cite{mal}) implique l'exactitude des autres flèches verticales. L'exactitude du dernier complexe $(({\cal E}^{0,q}_{_U} )^{\oplus p_m};\bar{\partial}_{_J } )_{ q \geq 0}$ implique l'exactitude du complexe  
$$
({\cal R}^{\cal E}(\varphi _{m-1})   \otimes_{_{{\cal E}_U }}{\cal E}^{0,q}_U;\bar{\partial}_{_J } )_{ q \geq 0} ,
$$ 
où ${\cal R}^{\cal E}(\varphi _{m-1})$ désigne le faisceau des relations de $\varphi _{m-1}$ sur le faisceau ${\cal E}_{_{X} }$. En procédant par récurrence décroissante et en utilisant l'exactitude des complexes en $\bar{\partial}_{_J }$ et l'exactitude des flèches verticales on obtient finalement l'exactitude du complexe 
$$
({\cal F}^{\infty}_{|_U}\otimes_{_{{\cal E}_{_U}  }}{\cal E}^{0,q}_{_U};\bar{\partial} _{_{\cal F}})_{ q \geq 0} . 
$$
Faisons maintenant le point de la situation obtenue jusqu'ici. On est parti d'un faisceau analytique cohérent ${\cal F}$ pour obtenir un faisceau de  ${\cal E}$-modules ${\cal F}^{\infty}$ admettant des ${\cal E}$-résolutions locales de longueur finie lequel est muni d'une connexion $\bar{\partial}_{_{\cal F} }$ de type $(0,1)$ intégrable telle que le noyau de celle-ci soit le faisceau analytique cohérent de départ ${\cal F}$. De manière générale on a la caractérisation différentielle suivante.
\begin{theoreme}\label{cdif} 
Soit X une variété complexe et soit ${\cal G}$ un faisceau de  ${\cal E}_X$-modules qu'on suppose muni d'une connexion 
$\bar{\partial}:{\cal G}\longrightarrow  {\cal G} \otimes_{_{{\cal E}_X }}{\cal E}^{0,1}_X $ de type $(0,1)$ telle que $\bar{\partial}^{2}=0$. Si de plus le faisceau ${\cal G}$ admet des ${\cal E}$-résolutions locales de longueur finie, alors le faisceau de  ${\cal O}_X$-modules $Ker\bar{\partial}\subset\cal G$ est analytique cohérent, on a les égalités 
$
{\cal G}=(Ker\bar{\partial})\cdot {\cal E}_X\cong
(Ker\bar{\partial})\otimes_{_{{\cal O}_{_X}}}{\cal E}_X
$ et  la connexion $\bar{\partial}$ coïncide, à isomorphisme canonique près, avec l'extension naturelle $\bar{\partial}_{Ker\bar{\partial}} $ associée au faisceau analytique cohérent $Ker\bar{\partial} $.
\end{theoreme}
Le théorème précédent montre donc qu'on est dans la même situation que celle décrite précédemment. Bien évidemment le théorème $\ref{cdif}$ constitue une généralisation du théorème $\ref{Gro1}$. Considérons maintenant la définition suivante.
\begin{defi}
Un couple $({\cal G},\bar{\partial})\equiv{\cal G}_{\bar{\partial}}$ où ${\cal G}$ et $\bar{\partial}$ vérifient les hypothèses du théorème 1 est appelé faisceau $\bar{\partial}$-cohérent. Un morphisme $\varphi :{\cal A}_{\bar{\partial}_1}\longrightarrow {\cal B}_{\bar{\partial}_2}$ de faisceaux $\bar{\partial}$-cohérents est un 
morphisme de faisceaux de ${\cal E}_X$-modules tels que le diagramme suivant soit commutatif
\begin{diagram}[height=1cm,width=1cm]
{\cal A}\otimes_{_{{\cal E}_{_X}}}{\cal E}^{0,1}_X&\rTo^{\varphi\otimes \I_{(0,1)} }&{\cal B}\otimes_{_{{\cal E}_{_X}}}{\cal E}^{0,1}_X
\\
\uTo^{\bar{\partial}_1}&              &\uTo^{\bar{\partial}_2}  
\\
{\cal A}&\rTo^{\varphi}&\cal B
\end{diagram}  
\end{defi}
Le théorème $\ref{cdif} $ et la fidélité plate du faisceau ${\cal E}_X$ sur ${\cal O}_X$ montrent que sur une variété complexe on a une équivalence exacte entre la catégorie ${\cal O}{\bf Coh} $ des faisceaux analytiques cohérents et la catégorie $\bar{\partial}{\bf Coh} $ des faisceaux $\bar{\partial}$-cohérents. Plus explicitement on a le foncteur $\infty$ qui agit de la façon suivante:
\begin{eqnarray*}
{\cal F}\in{\cal O}{\bf Coh} &\stackrel{\infty}{\longmapsto}&{\cal F}^{\infty}_{\bar{\partial}_{_{\cal F}}}\in\bar{\partial}{\bf Coh} 
\\
 \varphi \in Hom_{_{{\cal O}_X}}({\cal A},{\cal B})&\stackrel{\infty}{ \longmapsto}&\varphi\otimes \I\in 
Hom({\cal A}^{\infty}_{\bar{\partial}_{_{\cal A}}},{\cal B}^{\infty}_{\bar{\partial}_{_{\cal B}}})
\end{eqnarray*}
et son inverse:
\begin{eqnarray*}
{\cal G}_{\bar{\partial}} \in\bar{\partial}{\bf Coh}  & \stackrel{\infty^{-1}}{ \longmapsto}&Ker\bar{\partial} \in {\cal O}{\bf Coh} 
\\
 \varphi \in Hom({\cal A}_{\bar{\partial}_1},{\cal B}_{\bar{\partial}_2}) &\stackrel{\infty^{-1}}{ \longmapsto}   
 & \varphi_{|..}\in Hom_{_{{\cal O}_X}}(Ker\bar{\partial}_1,Ker\bar{\partial}_2)
\end{eqnarray*}
Considérons maintenant le cas des faisceaux d'ideaux. Soit ${\cal I}\subseteq {\cal O}_X$ un faisceau d'ideaux de fonctions holomorphes (non nécessairement cohérent). On considère le faisceau d'ideaux de fonctions ${\cal C}^{\infty}$ à valeurs complexes
${\cal I}^{\infty}:= {\cal I}\cdot{\cal E}_X \subseteq {\cal E}_X$ et on remarque que la règle de Leibnitz implique que pour tout germes de $(0,1)$-champs $\xi \in{\cal E}(T^{0,1}_X)_x$ on a l'inclusion $\xi .\,{\cal I}^{\infty} _x\subseteq {\cal I}^{\infty}_x $, pour tout $x\in X$. De manière générale on a la définition suivante: 
\begin{defi}
Un faisceaux d'ideaux ${\cal J}\subseteq {\cal E}_X$ de fonctions ${\cal C}^{\infty}$ à valeurs complexes est dit $\bar{\partial}_{_J }$-stable si pour tout germe de $(0,1)$-champs $\xi \in{\cal E}(T^{0,1}_X)_x$ on a l'inclusion $\xi .\,{\cal J}_x\subseteq {\cal J}_x $, pour tout $x\in X$.
\end{defi}
La $\bar{\partial}_{_J }$-stabilité d'un faisceaux d'ideaux ${\cal J}\subseteq {\cal E}_X$ de fonctions ${\cal C}^{\infty}$ à valeurs complexes implique évidemment qu'on peut considérer l'opérateur $\bar{\partial}_{_J }$ comme une connexion 
$\bar{\partial}_{_J } :{\cal J}\longrightarrow  {\cal J} \otimes_{_{{\cal E}_X }}{\cal E}^{0,1}_X$
de type $(0,1)$ intégrable sur le faisceaux ${\cal J}$. Dans le cas où ${\cal J}={\cal I}^{\infty}$ on a par conséquence de la platitude de l'anneau ${\cal E}_{X,x}$ sur l'anneau ${\cal O}_{X,x}$ que la connexion en question coïncide, à isomorphisme canonique prés, avec la connexion canonique $\bar{\partial}_{_{\cal I}}$ associée au faisceau ${\cal I}$.
Une conséquence immédiate du théorème précédent est le corollaire suivant:
\begin{corol}\label{corol1}   Soit ${\cal J}\subseteq {\cal E}_X$ un faisceaux d'ideaux de fonctions ${\cal C}^{\infty}$ à valeurs complexes  $\bar{\partial}_{_J }$-stable admettant des ${\cal E}$-résolutions locales de longueur finie. Alors le faisceaux d'ideaux ${\cal J}\cap {\cal O}_X$ est analytique cohérent et $({\cal J}\cap {\cal O}_X)\cdot{\cal E}_X={\cal J} $, $($autrement dit le faisceau d'ideaux
${\cal J}\cap {\cal O}_X$ est un ${\cal O}$-module localement de type fini et ces générateurs locaux sur ${\cal O}_X$ sont aussi des générateurs locaux du faisceau ${\cal J}$ sur ${\cal E}_X)$.
\end{corol} 
Concrètement pour vérifier la $\bar{\partial}_{_J }$-stabilité  du faisceau ${\cal J}$ il suffit de faire un choix arbitraire de repères locales 
$(\xi_1,...,\xi_n)\in{\cal E}(T^{0,1}_X)^{\oplus n}(U)$, de générateurs $(\psi_1,...,\psi_p)\in{\cal J}^{\oplus p}(U)$ et de montrer, pour tout $x\in U$, l'existence de germes de fonctions $f_{k,l,j}\in {\cal E}_{_{X,x} }$  qui vérifient les égalités 
$$
\xi _{k,x} \,.\psi_{l,x}=\sum_{j=1}^p\,f_{k,l,j}\cdot\psi_{j,x} .
$$
On remarque qu'en général le fait qu'un faisceau d'ideaux ${\cal J}\subseteq {\cal E}_X$ soit un $\cal E$-module localement de type fini n'implique pas nécessairement que le faisceaux d'ideaux ${\cal J}\cap {\cal O}_X$ soit un $\cal O$-module localement de type fini. On a le contre-exemple suivant.
\\
{\bf{Contre-exemple.} } Soit $X=\C$ et ${\cal J}:={\cal E}(\psi)\subseteq {\cal E}_{_\C}$ le faisceaux d'ideaux  de fonctions ${\cal C}^{\infty}$ à valeurs complexes engendré sur ${\cal E}_{_\C}$ par la fonction ${\cal C}^{\infty}$ sur $\C$, 
$
\psi(z):=\exp(-1/x^2)\cdot \sin(1/x)+i\,y
$, $(z=x+i\,y)$.
On a que $({\cal J}\cap {\cal O}_{_\C})_z=0$ pour 
$z=0,\;({\cal J}\cap {\cal O}_{_\C})_z=m({\cal O}_{_{\C,z} })$ pour tout $z=1/(k\pi),\;k\in\Z$ et 
$({\cal J}\cap {\cal O}_{_\C})_z={\cal O}_{_{\C,z} }$ pour $ z\not=0,\,1/(k\pi),\;k\in\Z$. Le faisceau d'ideaux ${\cal J}\cap {\cal O}_{_\C}$ n'est pas un $\cal O$-module localement de type fini. En effet pour tout voisinage ouvert $U\subset \C$ tel que $0\in \overline{U}$ on a $({\cal J}\cap {\cal O}_{_\C})(U)=0$. Ceci signifie bien evidemment que tous les morphismes 
$\varphi :{\cal O}^{\oplus r} _{_U}\longrightarrow ({\cal J}\cap {\cal O}_{_\C})_{_{|_U}}$ sont nuls. 
\\
\\
Le corollaire $\ref{corol1}$ montre donc l'existence d'une équivalence exacte entre la catégorie des faisceaux d'ideaux de fonctions holomorphes cohérents et la  catégorie des faisceaux d'ideaux ${\cal J}\subseteq {\cal E}_X$ de fonctions ${\cal C}^{\infty}$ à valeurs complexes admettant des ${\cal E}$-résolutions locales de longueur finie, qui sont stables par rapport aux dérivations le long des champs de vecteurs de type (0,1). On remarque enfin que la cohomologie des faisceaux cohérents ($\bar{\partial}$-cohérents) sur une variété complexe peut se calculer, grâce à l'isomorphisme fonctoriel de De Rham-Weil (voir par exemple les ouvrages \cite{Dem}, (\cite{Gri-Ha} et(\cite{Wel}),  par la formule suivante:
$$
H^q(X,{\cal G}_{\bar{\partial}}):=H^q (X,Ker\bar{\partial})\cong H^q(\Gamma(X,{\cal G} \otimes_{_{{\cal E}_X }}{\cal E}^{0,*}_X);\bar{\partial}) , 
$$ 
qui constitue une généralisation du théorème de Dolbeault. Un cas particulier (ou une généralisation si on veut) de la formule précédente est la suivante: 
\begin{eqnarray*} 
H^q(X,{\cal G}_{\bar{\partial}}  \otimes_{_{{\cal E}_X }}{\cal E}^{p,0}_{\bar{\partial}_{_{J,p} }  }) :=H^q(X,(Ker\bar{\partial})\otimes_{_{{\cal O}_X }}{\cal O}(\Omega ^p_X)) \cong
 H^q(\Gamma(X,{\cal G} \otimes_{_{{\cal E}_X }}{\cal E}^{p,*}_X);\bar{\partial}_{\pi} )=:H^{p,q}(X,{\cal G}_{\bar{\partial}})
\end{eqnarray*} 
où $\bar{\partial}_{_{J,p} }:=(-1)^p  \bar{\partial}_{_J }$, ${\cal G}_{\bar{\partial}}  \otimes_{_{{\cal E}_X }}{\cal E}^{p,0}_{\bar{\partial}_{_{J,p} } }:=({\cal G} \otimes_{_{{\cal E}_X }}
{\cal E}^{p,0}_X \,,\bar{\partial}_{\pi })$ et 
$
\bar{\partial}_{\pi }:{\cal G}\otimes_{_{{\cal E}_X }}{\cal E}^{p,0}_X\longrightarrow  {\cal G} \otimes_{_{{\cal E}_X }}{\cal E}^{p,1}_X 
$ 
désigne la connexion sur le produit $ {\cal G}\otimes_{_{{\cal E}_X }}{\cal E}^{p,0}_X$, laquelle est définie par la règle de Leibnitz.
\newpage
\section{Idée de la preuve du théorème $\ref{cdif}$  dans le cas des faisceaux de ${\cal E}$-modules localement libres avec la technique de type Nash-Moser}
\subsection{Expression locale de la condition d'intégrabilité $\bar{\partial}^2=0$ dans le cas des faisceaux de ${\cal E}$-modules localement libres}\label{explocm0}  
A partir de maintenant on va noter par $M_{k,l}({\cal E}^{0,q}_X(U))$ l'espace des matrices $k\times l$ à coefficients dans l'espace des $(0,q)$ formes ${\cal E}^{0,q}_X(U)$. Soit $\psi: {\cal E}_{_{U }  }^{\oplus r }\stackrel{\simeq}{\longrightarrow}
{\cal G}_{{|_{U } } } $ une trivialisation locale du faisceau ${\cal G}$. Le fait que $\psi$ est surjective implique l'existence d'une matrice $\omega ^{0,0}\in M_{p_0,p_0}({\cal E}^{0,1}_X(U))$ telle que 
$$
\bar{\partial}\psi=\psi\cdot \omega^{0,0}
$$
On obtient alors le diagramme commutatif suivant:
\begin{diagram}[height=1cm,width=1cm]
{\cal G}_{_{|_U} }&\rTo^{\bar{\partial}} &{\cal G}_{_{|_U}}\otimes_{_{{\cal E}_{_U}}}{\cal E}^{0,1}_{_U}
                                                                               &\rTo^{\bar{\partial}} &    {\cal G}_{_{|_U}}
                                                                                                        \otimes_{_{{\cal E}_{_U}}}{\cal E}^{0,2}_{_U} 
\\
\uTo^{\psi}_{\wr} &                             &\uTo^{\wr}_{\psi \otimes \I_{(0,1)} }&                               &\uTo^{\wr}_{\psi \otimes \I_{(0,2)} }
\\
{\cal E}^{\oplus p_0} _{_U}&\rTo^{\bar{\partial}_{_J }+\omega ^{0,0} }&({\cal E}^{0,1}_{_U})^{\oplus p_0} &\rTo^{\bar{\partial}_{_J }+\omega ^{0,0} }&({\cal E}^{0,2}_{_U})^{\oplus p_0} 
\end{diagram} 
L'hypothèse d'intégrabilité $\bar{\partial}^2=0$ est équivalente localement à l'égalité $0=\bar{\partial}^2\psi$. En explicitant celle-ci on a:
\begin{eqnarray*} 
0=\bar{\partial}^2\psi=\bar{\partial}(\psi\cdot \omega^{0,0})=\bar{\partial}\psi\wedge \omega^{0,0}+\psi \cdot \bar{\partial}_{_J }\omega ^{0,0} 
=\psi(\bar{\partial}_{_J }\omega ^{0,0} +\omega ^{0,0}  \wedge \omega ^{0,0})
\end{eqnarray*}
Le fait que $\psi$ soit injective implique alors la relation
$$
\bar{\partial}_{_J }\omega ^{0,0} +\omega ^{0,0} \wedge \omega ^{0,0} =0
$$
Donc l'hypothèse d'intégrabilité $\bar{\partial}^{2}=0$ s'exprime localement par cette relation. Dans la suite de cette section on désignera par $\Omega (U)\subset M_{p_0,p_0}({\cal E}^{0,1}_X(U))$ l'ensemble constitué par des matrices $\omega ^{0,0}$ qui vérifient la condition en question. Les éléments de cet ensemble seront appelés $calibrations$.
\subsection{Formulation du problème différentiel dans le cas des faisceaux de ${\cal E}$-modules localement libres} \label{forpbm0} 
On veut trouver pour chaque $x\in U$ un voisinage ouvert $V$ de $x$ et un élément $\eta^{0,0}\in M_{p_0,p_0}({\cal E}_X(V))$ tel que 
$g_0\equiv g_0(\eta):=\I_{p_0}+\eta^{0,0}\in GL(p_0, {\cal E}(V)),\;\eta\equiv \eta^{0,0}$, qui soit solution de l'équation différentielle
$$
\bar{\partial}(\psi\cdot g_0)=0.
$$
Si on atteint ce but on obtiendra le diagramme commutatif suivant:
\begin{diagram}[height=1cm,width=1cm]
0&\rTo  &(Ker\,\bar{\partial})_{_{|_V}}&\rTo& {\cal G}_{_{|_V} }&\rTo^{\bar{\partial}} &{\cal G}_{_{|_V}}\otimes_{_{{\cal E}_{_V}}}{\cal E}^{0,1}_{_V}&
\\
&      &\uTo^{\psi_{\eta} }_{\wr}&            &\uTo^{\psi_{\eta}}_{\wr}&          &\uTo^{\wr}_{\psi_{\eta} \otimes \I_{(0,1)} }& 
\\
0&\rTo &{\cal O}^{\oplus p_0}_{_V}& \rTo &{\cal E}^{\oplus p_0}_{_V}& \rTo^{\bar{\partial}_{_J }} &({\cal E}^{0,1}_{_V})^{\oplus p_0}& 
\end{diagram} 
où $\psi_{\eta}=\psi\cdot g_0(\eta)$, lequel permet de conclure dans le cas des faisceaux localement libres. L'équation précédente est équivalente à l'équation
$
\psi(\omega ^{0,0}\cdot g_0+\bar{\partial}_{_J }\eta^{0,0})=0 .
$
L'injectivité de $\psi$ implique alors que l'équation précédente est équivalente à l'équation
$$
(S_{\omega})\;:\; \bar{\partial}_{_J }\eta^{0,0}+\omega ^{0,0}\wedge  \eta^{0,0}+\omega^{0,0}=0 .
$$ 
On a la proposition suivante.
\begin{prop}
La condition $\bar{\partial}_{_J }\omega ^{0,0} +\omega ^{0,0} \wedge \omega ^{0,0}=0$ est la condition d'intégrabilité du problème différentiel quasi-linéaire $(S_{\omega})$, $($le problème est quasi-linéaire car on cherche $\eta^{0,0}$ dans un espace non linéaire$)$.
\end{prop}  
Il est élémentaire de vérifier la nécessité de la condition $\bar{\partial}_{_J }\omega ^{0,0} +\omega ^{0,0} \wedge \omega ^{0,0}=0$. La matrice
$$
\omega ^{0,0}_{\eta}:=g_0^{-1}(\bar{\partial}_{_J }\eta^{0,0}+\omega ^{0,0}\wedge  \eta^{0,0}+\omega^{0,0})   
$$
vérifie les relations
\begin{eqnarray*}
&\bar{\partial}\psi_{\eta} =\psi_{\eta} \cdot \omega^{0,0}_{\eta} , &
\\
\\
&\bar{\partial}_{_J }\omega ^{0,0}_{\eta}  +\omega ^{0,0}_{\eta}  \wedge \omega ^{0,0} _{\eta} =0 .&
\end{eqnarray*}  
On voudrait alors appliquer un procédé itératif pour faire de sorte que $\omega^{0,0}_{\eta}=0 $, ce qui équivaut à résoudre le système $(S_{\omega} )$. Dans la suite on appellera élément de recalibration un élément $\eta^{0,0}\in M_{p_0,p_0}({\cal E}_X(U))$ tel que $g_0=\I_{p_0}+\eta^{0,0}\in GL(p_0, {\cal E}(U))$ et on désignera par ${\cal P}(U)$ l'ensemble constitué par tels éléments. Venons-en maintenant  à un  préliminaire technique avant d'exposer l'idée de la preuve de l'existence des solutions pour le système différentiel $(S_{\omega} )$.

\subsection{Choix des normes et opérateur de Leray-Koppelman}\label{intLK} 
A partir de maintenant on va supposer dans cette section que $U=B_1$ est la boule  unité. Soit $u=\sum'_{|I|=q}u_I\, d\bar{z}_I$ est une $(0,q)$-forme à coefficients des $(k,l)$-matrices à coefficients dérivables jusqu'à l'ordre $h\geq 0$. On définit une norme de Hölder invariante par changement d'échelle
$$
\|u\|_{r,\,h,\,\mu ,\,q}:=\sum _{\displaystyle
\scriptstyle   |I|=q 
\atop
\scriptstyle |\alpha |\leq h }
S_{|\alpha |} \,r^{|\alpha |+q}\,\|\partial^\alpha u_I\|_{r,\mu } 
$$    
où $\|\cdot\|_{r,\mu }$ est la norme de Hölder invariante usuelle d'une fonction, $\mu \in (0,1)$ une constante fixée une fois pour toutes dans notre problème et $(S_k)_{k\geq 0}\subset (0,\infty)$ une suite de poids (cette suite sera choisie à décroissance assez rapide de façon à rendre en particulier les séries convergentes).
\\
On remarque que si le degré $q$ est $\geq 1$ on a que la norme $\|u\|_{r,h,\mu ,q}$ tend vers zéro lorsque le rayon $r$ tend vers zéro. 
\\
On considère maintenant l'opérateur de Leray-Koppelman classique de la boule de rayon $r$ (le lecteur peut consulter avec profit les ouvrages classiques de Henkin-Leiterer \cite{He-Le}, de Range \cite{Ra} et l'article de Harvey-Polkin \cite{Ha-Po})
$$
T_{r,q} : {\cal C}^{h,\mu} _{0,q+1}(\bar{B} _r, M_{k,l}(\C))\longrightarrow  {\cal C}^{h,\mu } _{0,q}(B_r, M_{k,l}(\C))
$$
Il existe une suite de poids $S=(S_k)_{k\geq 0}$ de la norme de Hölder introduite précédemment telle que pour toute forme différentielle $u\in {\cal C}^{h,\mu } _{0,q+1}(\bar{B} _r , M_{k,l}(\C))$ on a l'estimation intérieure:
\begin{eqnarray}\label{Webineq} 
\|T_{r,q} \,u\|_{r(1-\sigma ),\,h+1,\,\mu ,\,q} \leq C\cdot \sigma ^{-s(h)}\cdot\|u\|_{r,\,h,\,\mu ,\,q+1 }   
\end{eqnarray}
avec $C>0$ une constante $ind\acute{e} pendante$ de la régularité $h$ et $\sigma \in (0,1)$, $s(h)\in \N$ une fonction affine strictement croissante. Pour simplifier les notations on identifiera dans la suite de cette section 
$\|\cdot\|_{r,h,\mu ,q}\equiv \|\cdot\|_{r,h},\;T_{r,q}\equiv T_r $  et $\|\partial^h f\|_{\bullet}\equiv\sum_{|\alpha |=h}\|\partial^{\alpha }  f\|_{\bullet}$. 
\subsection{Esquisse du schéma de convergence rapide de type Nash-Moser dans le cas des faisceaux de ${\cal E}$-modules localement libres}\label{NMm0} 

\subsubsection{Estimation fondamentale du schéma de convergence rapide dans le cas des faisceaux de ${\cal E}$-modules localement libres}\label{estm0}  
On désigne par $\varepsilon \in (0,1/2)$ une constante fixée telle que pour toutes les matrices $A\in M_{p_0,p_0}(\C)$ telle que $\|A\|<\varepsilon$ on a l'inversibilité de la matrice $\I_{p_0}+A $. 
\begin{prop}
Supposons donnés $\omega ^{0,0}\in \Omega (B_1),\; r,\,\sigma \in (0,1),\;h\in \N$ et les poids $S_j>0,\;j=0,...,h+1$ de la norme de Hölder tels que l'estimation  $(\ref{Webineq})$ soit satisfaite. Supposons que le rayon $r$ soit suffisamment petit pour assurer l'estimation 
$$
C\cdot \sigma^{-s(h)}\cdot \|\omega ^{0,0}\|_{r,h} <\varepsilon .  
$$
Supposons de plus que  le poids  $S_{h+1}$ soit suffisamment petit pour pouvoir assurer l'estimation
\begin{eqnarray*}
S_{h+1} \|\partial^{h+1} \omega^{0,0}_I \|_{r ,\mu } \leq \| \omega^{0,0}_I \|_{r ,\mu }
\end{eqnarray*}
où $\omega ^{0,0}_I$ sont les coefficients de la forme $\omega ^{0,0}$ par rapport aux coordonnées choisies. Alors $\eta^{0,0}:=-T_r\,\omega ^{0,0}$ est un paramètre de recalibration 
tel que
\begin{eqnarray*}
\|\omega^{0,0}_{\eta} \|_{r(1-\sigma ),h+1} \leq 6C \cdot\sigma^{-s(h)}\cdot \|\omega ^{0,0}\|^2_{r,h}
\end{eqnarray*}
\end{prop} 
Sans l'hypothèse sur le poids $S_{h+1}$ l'estimation précédente est valable pour $h$ à la place de $h+1$, pour toutes les matrices $\omega ^{0, 0}$ qui vérifient la relation $\bar{\partial}_{_J }\omega ^{0,0} +\omega ^{0,0} \wedge \omega ^{0,0} =0$. L'hypothèse sur les poids $S_{h+1}$, comme on verra mieux ensuite, joue un rôle fondamental pour la convergence vers une solution $\ci$ du problème $(S_{\omega})$.
\\
$Preuve$. En utilisant la définition du paramètre $\eta^{0,0}$ et la formule d'homotopie pour l'opérateur $ \bar{\partial}_{_J }$ on a:
\begin{eqnarray*}
&\omega ^{0,0}_{\eta}=g_0^{-1}(T_r\,  \bar{\partial}_{_J }\omega ^{0,0}-\omega ^{0,0}\wedge T_r\,\omega ^{0,0})=&
\\
\\
&=-  g_0^{-1}(T_r\,(\omega ^{0,0}\wedge \omega ^{0,0}) +\omega ^{0,0}\wedge T_r\,\omega ^{0,0})&
\end{eqnarray*}   
Le fait que $\varepsilon\in (0,1/2) $ implique que $\|g_0^{-1}\|_{r(1-\sigma ),h+1}<2$ (voir les détails dans la preuve complète, prop $\ref{3.4}$ dans la sous-section $\ref{ss3.4.1}$). L'hypothèse sur le poids $S_{h+1}$ implique l'inégalité suivante:
$$
\|\omega ^{0,0} \wedge T_r\,\omega ^{0,0}\|_{r(1-\sigma ),h+1}\leq 2\|\omega ^{0,0}\|_{r,h}\cdot \|T_r\,\omega ^{0,0}\|_{r(1-\sigma ),h+1} .
$$
On obtient alors l'inégalité:
\begin{eqnarray*}
\|\omega^{0,0}_{\eta} \|_{r(1-\sigma ),h+1} \leq 2\Big(C\cdot \sigma^{-s(h)}\cdot \|\omega ^{0,0}\|^2_{r,h} 
+ 2C\cdot \sigma^{-s(h)}\cdot \|\omega ^{0,0}\|^2_{r,h}\Big)
\end{eqnarray*}
laquelle permet de conclure. \hfill $\Box$
\subsubsection{Idée du procédé itératif dans le cas des faisceaux de ${\cal E}$-modules localement libres}\label{procitm0} 
Les calibrations $\omega^{0,0}_k\in \Omega (\bar{B}_{r_k})$ obtenues au $k$-ième pas du procédé itératif 
sont définies par la formule récursive $\omega^{0,0} _{k+1}=\omega^{0,0}  _{k,\,\eta_{k+1} }$
où $r_{k+1}:=r_k(1-\sigma _k)$ et où $\sigma _k\in (0,1)$ est un paramètre qui contrôle la décroissance des rayons des boules, (le rayon initial $r_0$ étant choisi suffisamment petit). On choisit les quantités $\sigma _k\in (0,1)$ de telle sorte que la série $\sum \sigma _k$ soit convergente. Le rayon limite $r_{\infty}:=\lim_{k\rightarrow +\infty}r_k$ est alors non nul. Le paramètre $\eta^{0,0} _{k+1} \in {\cal P}(\bar{B}_{r_{k+1}})$
qui contrôle la recalibration des éléments $\omega^{0,0}  _k$, $k\geq 0$ au $k$-ième pas du procédé itératif est défini par la formule  $$\eta^{0,0}_{k+1}:=-T_{r_k}\omega^{0,0}_k .$$
 Les poids sont choisis de telle sorte que les estimations suivantes soient satisfaites pour tout  entier $k\geq 0$.
\begin{eqnarray*}
&S_{k+1} \|\partial^{k+1} \omega^{0,0}_{k,\,I}  \|_{r_k ,\,\mu }\leq\| \omega^{0,0}_{k,\,I}  \|_{r_k ,\,\mu }& 
\\\nonumber
\\
&S_{k+1} \|\partial^{k+1} g_0(k)^{\pm 1}  \|_{r_k ,\,\mu }\leq 2^{-k-1}\|g_0(k)^{\pm 1}   \|_{r_k ,\,\mu }  &
\end{eqnarray*}
où 
$$
g_0(k) :={\displaystyle \prod_{0\leq j\leq k }^{\longrightarrow}}g_{0,j} 
$$
(on pose par définition $g_0:=\I_{p_0} $. Le symbole de produit avec une flèche vers la droite désigne le produit non commutatif de termes qui sont écrits en ordre croissant de l'indice vers la droite). 
La dernière inégalité sert à assurer un bon fonctionnement du procédé itératif, plus précisément elle permet d'appliquer la proposition précédente à toutes les étapes du procédé. On pose par définition
$$
a_{k} :=\|\omega ^{0,0}_k\|_{r_k,\,k}
\qquad
\mbox{et} 
\qquad
b_k:=H\cdot\sigma_k^{- s(k) }\cdot a_k
$$
Avec les notations introduites précédemment on a la proposition suivante.
\begin{prop} 
Pour tout entier $k\geq 0$ on a les estimations suivantes;
\begin{eqnarray*}
&a_{k+1} \leq H\cdot\sigma_k^{-s(k)}\cdot a^2_{k}\leq 1 ,&
\\\nonumber
\\
&\|\eta^{0,0}_{k+1} \|_{r_{k+1}  ,\,k+1 } \leq  b_k   <\varepsilon <1/2&
\end{eqnarray*}
et les quantités $a_k,\;b_k$ tendent $($avec la bonne vitesse$)$ vers zéro lorsque $k$ tend vers plus l'infini.
\end{prop}
En conclusion la limite 
$$
\eta^{0,0}=-\I_{p_0}+\lim _{k\rightarrow \infty}g_0(k)
$$
est une solution du système différentiel $(S_{\omega })$.
\section{Introduction au cas de ${\cal E}$-résolution locale de profondeur homologique égale à un $(m=1)$}
\subsection{Expression locale de la condition d'intégrabilité $\bar{\partial}^2=0$ dans le cas $m=1$} \label{explocm1} 
Nous commençons par prouver le lemme élémentaire suivant.
\begin{lem} 
Soit X une variété complexe et soit ${\cal G}$ un faisceau de  ${\cal E}_X$-modules. Si le faisceau ${\cal G}$ admet des ${\cal E}$-présentations locales, soit par exemple
$$
 {\cal E}_{_U}^{\oplus p_1 }\stackrel{\varphi}{ \longrightarrow }{\cal E}_{_U}^{\oplus p_0 }\stackrel{\psi}{\longrightarrow}
{\cal G}_{_{|_U} }\rightarrow 0
$$
une ${\cal E}$-présentation au dessus d'un ouvert $U$, alors l'existence d'une connexion 
$\bar{\partial}$ de type $(0,1)$ sur le faisceau ${\cal G}_{_{|_U} }$ telle que 
$\bar{\partial}^{2}=0$, implique l'existence de matrices $\omega ^{s,0}\in M_{p_s,p_s}({\cal E}^{0,1}_X(U)),\;s=0,1$ et 
$\omega ^{0,1}\in M_{p_1,p_0}({\cal E}^{0,2}_X(U))$ telles que 
$
\bar{\partial}\psi=\psi\cdot \omega^{0,0}
$ et les relations
\begin{eqnarray}
&\bar{\partial}_{_J }\varphi+\omega ^{0,0}\cdot \varphi =  \varphi\cdot \omega ^{1,0}\label{1*} , &
\\\nonumber
\\
&\bar{\partial}_{_J }\omega ^{0,0} +\omega ^{0,0}\wedge \omega ^{0,0}=\varphi \cdot\omega ^{0,1}
\label{2*}  &
\end{eqnarray}
soient satisfaites. Réciproquement l'existence de matrices $\omega ^{s,0},\;s=0,1 $ et $ \omega ^{0,1} $ qui vérifient les relations $(\ref{1*})$ et $(\ref{2*} )$, implique l'existence d'une connexion $\bar{\partial}$ de type $(0,1)$ sur le faisceau ${\cal G}_{_{|_U} }$ telle que 
$
\bar{\partial}\psi=\psi\cdot \omega^{0,0}
$ et
$\bar{\partial}^{2}=0
$.
\end{lem} 
$Preuve$. La ${\cal E}$-présentation de ${\cal G}_{_{|_U}}$ considérée dans l'hypothèse implique l'existence des ${\cal E}$-présentations
$$
({\cal E}^{0,q}_{_U})^{\oplus p_1}\longrightarrow ({\cal E}^{0,q}_{_U})^{\oplus p_0}\longrightarrow  {\cal G}_{_{|_U}}\otimes_{_{{\cal E}_{_U}}}{\cal E}^{0,q}_{_U}\longrightarrow 0
$$ 
pour $q\geq 0$. On aura alors l'existence d'une matrice $\omega ^{0,0}\in M_{p_0,p_0}({\cal E}^{0,1}_X(U))$ telle que 
$
\bar{\partial}\psi=\psi\cdot \omega^{0,0}
$.
 En appliquant la connexion $\bar{\partial}$ à l'identité $\psi \circ \varphi=0$ on obtient la relation
$$
\psi(\bar{\partial}_{_J }\varphi +\omega ^{0,0}\cdot  \varphi)=0 .
$$
L'exactitude des 
${\cal E}$-présentations précédentes, (pour $q=1$), implique alors l'existence d'une matrice 
$
\omega ^{1,0}\in M_{p_1,p_1}({\cal E}^{0,1}_X(U))
$ 
telle que la relation $(\ref{1*} )$ soit satisfaite. On obtient alors le diagramme commutatif suivant, ayant des flèches verticales exactes:
\begin{diagram}[height=1cm,width=1cm]
0&                                       &0&                                                                &0
\\
\uTo&                                    &\uTo&                                                            &\uTo
\\
{\cal G}_{_{|_U} }&\rTo^{\bar{\partial}} &{\cal G}_{_{|_U}}\otimes_{_{{\cal E}_{_U}}}{\cal E}^{0,1}_{_U}
                                                                               &\rTo^{\bar{\partial}} &    {\cal G}_{_{|_U}}
                                                                                                        \otimes_{_{{\cal E}_{_U}}}{\cal E}^{0,2}_{_U} 
\\
\uTo^{\psi}&                             &\uTo_{\psi \otimes \I_{(0,1)} }&                               &\uTo_{\psi \otimes \I_{(0,2)} }
\\
{\cal E}^{\oplus p_0} _{_U}&\rTo^{\bar{\partial}_{_J }+\omega ^{0,0} }&({\cal E}^{0,1}_{_U})^{\oplus p_0} &\rTo^{\bar{\partial}_{_J }+\omega ^{0,0} }&({\cal E}^{0,2}_{_U})^{\oplus p_0} 
\\
\uTo^{\varphi_1}&                             &\uTo_{\varphi_1\otimes \I_{(0,1)} }&                               &\uTo_{\varphi_1 \otimes \I_{(0,2)} }
\\
{\cal E}^{\oplus p_1} _{_U}&\rTo^{\bar{\partial}_{_J }+\omega ^{1,0} }&({\cal E}^{0,1}_{_U})^{\oplus p_1}&\rTo^{\bar{\partial}_{_J }+\omega ^{1,0} }&({\cal E}^{0,2}_{_U})^{\oplus p_1} 
\end{diagram}
\\ 
L'hypothèse d'intégrabilité $\bar{\partial}^2=0$ implique 
$$
\psi(\bar{\partial}_{_J }\omega ^{0,0} +\omega ^{0,0}  \wedge \omega ^{0,0} )=0
$$
d'où l'existence d'une matrice $\omega^{0,1}  \in M_{p_1,p_0}({\cal E}^{0,2}_X(U))$ telle que la relation $(\ref{2*} )$ soit satisfaite.
\\
Pour prouver la réciproque du lemme il suffit de considérer la connexion quotient $\bar{\partial} $ obtenue par la connexion $\bar{\partial}_{_J }+\omega ^{0,0}\wedge \bullet $.\hfill $\Box$
\\
\\
En appliquant l' opérateur $\bar{\partial}_{_J }$ à la relation $(\ref{1*})$  on obtient l'égalité:
$$
\bar{\partial}_{_J }\omega ^{0,0}\cdot \varphi -\omega ^{0,0}\wedge  \bar{\partial}_{_J }\varphi =\bar{\partial}_{_J }\varphi\wedge \omega ^{1,0}+\varphi \cdot \bar{\partial}_{_J }\omega ^{1,0} .
$$
En utilisant la relation $(\ref{1*})$ dans l'égalité précédente on obtient 
\begin{eqnarray*} 
&\bar{\partial}_{_J }\omega ^{0,0}\cdot \varphi +\omega ^{0,0}\wedge \omega ^{0,0}  \cdot\varphi -\omega ^{0,0}\wedge \varphi \cdot \omega ^{1,0} =&
\\
\\
& =-\omega ^{0,0}\cdot \varphi \wedge \omega ^{1,0} + \varphi ( \omega ^{1,0}\wedge \omega ^{1,0}+  \bar{\partial}_{_J }\omega ^{1,0}) .&
\end{eqnarray*} 
En utilisant la relation $(\ref{2*})$ on a:
$$
\varphi (\bar{\partial}_{_J }\omega ^{1,0}+\omega ^{1,0}\wedge \omega ^{1,0}-\omega ^{0,1}\cdot \varphi)=0 .
$$
Dans notre cas $\varphi$ est injective. On déduit alors la relation
$$
\bar{\partial}_{_J }\omega ^{1,0}+\omega ^{1,0}\wedge \omega ^{1,0}=\omega ^{0,1}\cdot \varphi .
$$
En appliquant l'opérateur $\bar{\partial}_{_J }$ à la relation $(\ref{2*})$ on a:
$$
\bar{\partial}_{_J }\omega ^{0,0}\wedge \omega ^{0,0}  -\omega ^{0,0}\wedge  \bar{\partial}_{_J }\omega ^{0,0}  =\bar{\partial}_{_J }\varphi\wedge \omega ^{0,1}+\varphi \cdot \bar{\partial}_{_J }\omega ^{0,1} .
$$
En utilisant les relations $(\ref{1*})$ et $(\ref{2*})$ dans l'égalité précédente on obtient 
\begin{eqnarray*} 
&-\omega ^{0,0}\wedge \omega ^{0,0} \wedge \omega ^{0,0} +\varphi \cdot\omega ^{0,1}\wedge \omega ^{0,0} +\omega ^{0,0}\wedge \omega ^{0,0} \wedge \omega ^{0,0} -\omega ^{0,0}\wedge \varphi \cdot \omega ^{0,1} =&
\\
\\
& =-\omega ^{0,0}\cdot \varphi \wedge \omega ^{0,1} + \varphi ( \omega ^{1,0}\wedge \omega ^{0,1}+  \bar{\partial}_{_J }\omega ^{0,1}) .&
\end{eqnarray*} 
On a donc
$$
\varphi (   \bar{\partial}_{_J }\omega ^{0,1}- \omega ^{0,1}\wedge \omega ^{0,0}    +\omega ^{1,0}\wedge \omega ^{0,1} )=0 .
$$
L'injectivité de $\varphi $ implique alors la relation
$$
\bar{\partial}_{_J }\omega ^{0,1}- \omega ^{0,1}\wedge \omega ^{0,0}    +\omega ^{1,0}\wedge \omega ^{0,1} =0 .
$$
On a obtenu en conclusion les relations
$$
\bar{\partial}_{_J }\varphi+\omega ^{0,0}\cdot \varphi =  \varphi\cdot \omega ^{1,0} 
$$
et
$$
(*)
\left  \{
\begin{array}{lr}
\bar{\partial}_{_J }\omega ^{0,0} +\omega ^{0,0}\wedge \omega ^{0,0}=\varphi \cdot\omega ^{0,1} 
\\
\\
\bar{\partial}_{_J }\omega ^{1,0} +\omega ^{1,0}\wedge \omega ^{1,0}=\omega ^{0,1} \cdot \varphi 
\\
\\
\bar{\partial}_{_J }\omega ^{0,1} -\omega ^{0,1}\wedge \omega ^{0,0}+\omega ^{1,0}\wedge \omega ^{0,1}=0  .
\end{array}
\right.
$$
Il faut remarquer que dans la dernière expression on n'a pas de termes du type $\varphi \cdot \omega ^{\bullet,\bullet}$ ou $\omega ^{\bullet,\bullet}\cdot \varphi$. Les expressions $(*)$ constituent les expressions locales de la condition d'intégrabilité $\bar{\partial}^2=0$ dans le cas de longueur $m=1$ de la ${\cal E}$-résolution locale. La relation $(\ref{1*})$ est simplement une identité de commutation. 
\subsection{Introduction à la formulation du problème différentiel dans le cas $m=1$}\label{forpbm1}  
La partie principale de la preuve consiste à prouver l'existence, pour tout $x\in U$, d'un voisinage ouvert $V \subset U$ de $x$ et $g_0=\I_{p_0}+\eta^{0,0}\in GL(p_0,{\cal E}(V) ),\;g_1=\I_{p_1}+\eta^{1,0}\in GL(p_1,{\cal E}(V))$ solution du système différentiel
$$
(\Sigma )\;
\left  \{
\begin{array}{lr}
\bar{\partial}(\psi\cdot g_0)=0 
\\
\\
\bar{\partial}_{_J}(g_0^{-1} \cdot \varphi \cdot g_1)=0  .
\end{array}
\right.
$$
Si on atteint ce but on obtiendra le diagramme commutatif suivant:
\begin{diagram}[height=1cm,width=1cm]
&         &                              &                 &0&                                         &0                                                       \\
&         &                              &              &\uTo&                                      &\uTo                                                       \\
0&  \rTo  &(Ker\bar{\partial})_{_{|V} } &  \rTo   &{\cal G}_{|_V}& \rTo^{\bar{\partial}} 
 & {\cal G}_{|_V}\otimes_{_{{\cal E}_{_V}}}{\cal E}^{0,1}_{_V} 
\\
&          &\uTo^{\tilde{\psi}_{|..}}&                     &\uTo_{\tilde{\psi}}&                      &\uTo_{\tilde{\psi}\otimes\I_{(0,1)}}                   \\
0&\rTo  & {\cal O} _{_V}^{\oplus p_0} &  \rTo   &{\cal E}^{\oplus p_0}_{_V}&    \rTo^{\bar{\partial}_{_J }}
 & ({\cal E}^{0,1}_{_V} )^{\oplus p_0}&  
\\
&             &\uTo^{\tilde{\varphi}_{|..}} &                             &\uTo_{\tilde{\varphi}}&               &\uTo_{\tilde{\varphi}\otimes \I_{(0,1)} }
\\
0& \rTo  & {\cal O}
          _{_V}^{\oplus p_1} &  \rTo   &{\cal E}^{\oplus p_1}_{_V}&    \rTo^{\bar{\partial}_{_J }}
 & ({\cal E}^{0,1}_{_V} )^{\oplus p_1} 
\end{diagram} 
\\
avec $\tilde{\psi}:=\psi_{\eta}:=\psi\cdot g_0$ et $\tilde{\varphi} :=\varphi_{\eta}:=g_0^{-1} \cdot \varphi \cdot g_1$. En utilisant la fidélité plate de l'anneau ${\cal E}_{X,x}$ sur l'anneau ${\cal O}_{X,x}$ (voir la preuve complète, section $\ref{s3.5}$ pour plus de détails) on peut conclure. En reprenant un calcul fait dans le cas $m=0$, mais valable en tous les cas, on a:
$$
\bar{\partial}(\psi\cdot g_0)=\psi(\bar{\partial}_{_J }\eta^{0,0}+\omega ^{0,0}\wedge  \eta^{0,0}+\omega^{0,0}) .
$$
Si $\bar{\partial}(\psi\cdot g_0)=0$ alors l'hypothèse d'exactitude implique l'existence d'une matrice $\eta^{0,1}\in M_{p_1,p_0}({\cal E}^{0,1}_X(U))$ solution de l'équation
$$
\bar{\partial}_{_J }\eta^{0,0}+\omega ^{0,0}\wedge  \eta^{0,0}+\varphi \cdot \eta^{0,1} +\omega^{0,0}=0 .
$$
On a donc que l'équation précédente est équivalente avec l'équation $\bar{\partial}(\psi\cdot g_0)=0$. Si on pose par définition
$$
\omega ^{0,0}_{\eta}:=g_0^{-1}(\bar{\partial}_{_J }\eta^{0,0}+\omega ^{0,0}\wedge  \eta^{0,0}+\varphi \cdot \eta^{0,1} +\omega^{0,0} )   
$$
on aura la validité de la relation $\bar{\partial}\psi_{\eta} =\psi_{\eta} \cdot \omega^{0,0}_{\eta}$. De façon analogue, si on pose par définition
\begin{eqnarray*}
&\omega ^{1,0}_{\eta}:=g_1^{-1}(\bar{\partial}_{_J }\eta^{1,0}+\omega ^{1,0}\wedge  \eta^{1,0}+\eta^{0,1}\cdot\varphi _{\eta}  +\omega^{1,0} ) ,  &
\\
\\
&\omega ^{0,1}_{\eta}:=g_1^{-1}(\bar{\partial}_{_J }\eta^{0,1}+\omega ^{0,1}\wedge  \eta^{0,0}+\omega ^{1,0}\wedge  \eta^{0,1}+\eta^{0,1}\wedge \omega ^{0,0}_{\eta} +\omega^{0,1} )   &
\end{eqnarray*}  
on aura la validité des  relations (voir article pour les détails des calculs en général prop $\ref{3.2}$ de la sous-section $\ref{ss3.2}$)
$$
\bar{\partial}_{_J }\varphi_{\eta}+\omega ^{0,0}_{\eta}\cdot \varphi_{\eta} =  \varphi_{\eta}\cdot \omega ^{1,0} _{\eta}
$$
et
$$
\left  \{
\begin{array}{lr}

\bar{\partial}_{_J }\omega ^{0,0}_{\eta} +\omega ^{0,0}_{\eta}\wedge \omega ^{0,0}_{\eta}=\varphi _{\eta}\cdot\omega ^{0,1} _{\eta}
\\
\\
\bar{\partial}_{_J }\omega ^{1,0}_{\eta} +\omega ^{1,0}_{\eta}\wedge \omega ^{1,0}_{\eta}=\omega ^{0,1}_{\eta} \cdot \varphi _{\eta}
\\
\\
\bar{\partial}_{_J }\omega ^{0,1}_{\eta} -\omega ^{0,1}_{\eta}\wedge \omega ^{0,0}_{\eta}+\omega ^{1,0}_{\eta}\wedge \omega ^{0,1}_{\eta}=0 
\end{array}
\right.
$$
lesquelles sont analogues à la relation $(\ref{1*})$ et aux relations $(*)$ considérées précédemment. En utilisant l'hypothèse d'exactitude on obtient le lemme suivant.
\begin{lem}
Pour tout choix de matrices $\omega ^{s,0},\;\omega ^{0,1}, s=0,1$ qui vérifient les relations $\eqref{1*}$ et $(*)$ on a que l'existence d'une solution $g:=(g_0,g_1)$ du système différentiel $(\Sigma)$ est équivalente à l'existence d'une solution $\eta:=(\eta^{0,0},\eta^{1,0},\eta^{0,1}),\; g=g(\eta)$, du système différentiel quasi-linéaire
$$
(S_{\omega} )\;
\left  \{
\begin{array}{lr}
\bar{\partial}_{_J }\eta^{0,0}+\omega ^{0,0}\wedge  \eta^{0,0}+\varphi \cdot \eta^{0,1} +\omega^{0,0}=0
\\
\\
\bar{\partial}_{_J }\eta^{1,0}+\omega ^{1,0}\wedge  \eta^{1,0}+\eta^{0,1}\cdot\varphi _{\eta}  +\omega^{1,0} =0
\\
\\
\bar{\partial}_{_J }\eta^{0,1}+\omega ^{0,1}\wedge  \eta^{0,0}+\omega ^{1,0}\wedge  \eta^{0,1} +\omega^{0,1}=0
\end{array}
\right.
$$ 
qui n'est rien d'autre que le système différentiel  
$$
\left  \{
\begin{array}{lr}
\omega^{s,0}_{\eta}=0
\\
\\
\omega^{0,1}_{\eta}=0
\\
\\
s=0,1 .
\end{array}                                                                           
\right.
$$
\end{lem} 
On a besoin de l'hypothèse d'exactitude de la ${\cal E} $-résolution locale seulement pour prouver que les systèmes $(\Sigma)$ et $(S_{\omega } )$ sont équivalents. A partir du moment où on s'intéresse seulement au système $(S_{\omega} )$ l'hypothèse d'exactitude n'a plus aucun intérêt. En termes précis on a la proposition suivante.
\begin{prop}
Supposons données des matrices $\omega ^{s,0},\;\omega ^{0,1}, s=0,1$ et $\varphi$ telles que 
$$
\bar{\partial}_{_J}\varphi +\omega ^{0,0}\cdot \varphi = \varphi \cdot \omega^{1,0} .
$$
Alors les relations 
$$
(*)
\left  \{
\begin{array}{lr}
\bar{\partial}_{_J }\omega ^{0,0} +\omega ^{0,0}\wedge \omega ^{0,0}=\varphi \cdot\omega ^{0,1} 
\\
\\
\bar{\partial}_{_J }\omega ^{1,0} +\omega ^{1,0}\wedge \omega ^{1,0}=\omega ^{0,1} \cdot \varphi 
\\
\\
\bar{\partial}_{_J }\omega ^{0,1} -\omega ^{0,1}\wedge \omega ^{0,0}+\omega ^{1,0}\wedge \omega ^{0,1}=0 
\end{array}
\right.
$$
constituent les conditions d'intégrabilité du système différentiel $(S_{\omega})$.
\end{prop} 

\section{Idée de la preuve du théorème $\ref{cdif}$ dans le cas général d'une ${\cal E}$-résolution locale de longueur arbitraire}
\subsection{Première étape: présentation de l'expression locale de la condition d'intégrabilité $\bar{\partial}^2=0$}\label{exploc}  
On pose par définition $I_m:=\{(s,k)\;|\;s=0,...,m,\; k=-1,...,m-s,\;(s,k)\not= (0,-1) \} $. On utilisera dans la suite la convention qui consiste à négliger les termes d'une somme ou d'un produit si l'ensemble des indices sur lesquels on effectue ces opérations est vide. On a le lemme suivant.
\begin{lem} 
Soit donnée
$
0\rightarrow {\cal E}_{_U}^{\oplus p_m }\stackrel{\varphi_m }{\longrightarrow}  
{\cal E}_{_U}^{\oplus p_{m-1} }\stackrel{ \varphi _{m-1} }{\longrightarrow} \cdot \cdot\cdot
\stackrel{\varphi_2}{ \longrightarrow } {\cal E}_{_U}^{\oplus p_1 }\stackrel{\varphi_1}{ \longrightarrow }{\cal E}_{_U}^{\oplus p_0 }\stackrel{\psi}{\longrightarrow}
{\cal G}_{_{|_U} }\rightarrow 0
$
une  ${\cal E}$-résolution locale de longueur finie. Alors l'existence d'une connexion 
$\bar{\partial}$ de type $(0,1)$ sur le faisceau ${\cal G}$ telle que 
$\bar{\partial}^{2}=0$, implique l'existence des matrices $\omega^{s,k}  \in M_{p_{s+k} ,p_s}({\cal E}^{0,k+1}_X(U))$ pour $(s,k)\in I_m$ telles que si on utilise l'identification 
$
\varphi _s\equiv\omega ^{s,-1}
$, 
on aura la validité des relations 
$
\bar{\partial}\psi=\psi\cdot \omega^{0,0}
$ et
\begin{eqnarray}\label{3*} 
\bar{\partial}_{_J }\omega^{s,k} +\sum_{j=-1}^{k+1} (-1)^{k-j}  \omega^{s+j,k-j}\wedge \omega^{s,j}=0 . 
\end{eqnarray}
\end{lem}
On obtient alors le diagramme commutatif suivant, ayant des flèches verticales exactes:
\begin{diagram}[height=1cm,width=1cm]
0&                                       &0&                                                                &0
\\
\uTo&                                    &\uTo&                                                            &\uTo
\\
{\cal G}_{_{|_U} }&\rTo^{\bar{\partial}} &{\cal G}_{_{|_U}}\otimes_{_{{\cal E}_{_U}}}{\cal E}^{0,1}_{_U}
                                                                               &\rTo^{\bar{\partial}} &    {\cal G}_{_{|_U}}
                                                                                                        \otimes_{_{{\cal E}_{_U}}}{\cal E}^{0,2}_{_U} 
\\
\uTo^{\psi}&                             &\uTo_{\psi \otimes \I_{(0,1)} }&                               &\uTo_{\psi \otimes \I_{(0,2)} }
\\
{\cal E}^{\oplus p_0} _{_U}&\rTo^{\bar{\partial}_{_J }+\omega ^{0,0} }&({\cal E}^{0,1}_{_U})^{\oplus p_0} &\rTo^{\bar{\partial}_{_J }+\omega ^{0,0} }&({\cal E}^{0,2}_{_U})^{\oplus p_0} 
\\
\uTo^{\varphi_1}&                             &\uTo_{\varphi_1\otimes \I_{(0,1)} }&                               &\uTo_{\varphi_1 \otimes \I_{(0,2)} }
\\
{\cal E}^{\oplus p_1} _{_U}&\rTo^{\bar{\partial}_{_J }+\omega ^{1,0} }&({\cal E}^{0,1}_{_U})^{\oplus p_1}&\rTo^{\bar{\partial}_{_J }+\omega ^{1,0} }&({\cal E}^{0,2}_{_U})^{\oplus p_1}
\\
\uTo&         &\uTo&          &\uTo
\end{diagram} 
La figure $\ref{fig2}$  montre les  matrices $\omega^{s,k}$  qui sont représentées par des flèches dans le diagramme suivant lequel représente le complexe déterminé par la 
${\cal E}$-résolution locale $(\varphi ,\psi)$ dans le cas de longueur $m=4$.
\begin{figure}[hbtp]
\begin{center} 
\input dgt1.pstex_t
    \caption{}
    \label{fig2}
\end{center} 
\end{figure}    
\\
Dans la relation $(\ref{3*} )$ on utilise les conventions formelles $\omega ^{0,-1} :=0,\;\omega ^{-1,j}:=0$ et $\omega ^{s,k}:=0$ si $s\geq m+1$ ou $k\geq m-s+1$.
Le diagramme suivant montre les matrices qui interviennent dans la relation $(\ref{3*} )$ pour $(s,k)=(1,2)$, dans le cas de longueur $m=4$. 
\begin{figure}[hbtp]
\begin{center} 
\input dgt2.pstex_t       
    \caption{}
    \label{fig3}
\end{center} 
\end{figure}    
\\
On considère le cas $m=4$ ceci étant  le cas de longueur minimale pour laquelle on voit le problème en toute sa généralité. Les relations $(\ref{3*})$ pour les indices $k\geq 0$ constituent les expressions locales de la condition d'intégrabilité $\bar{\partial}^2=0$ de la connexion $\bar{\partial}$ relativement à la ${\cal E}$-résolution locale choisie. Les relations $(\ref{3*})$ pour les indices $(s,k)=(\bullet,-1)$ représentent simplement des identités de commutation. La liberté homologique qui caractérise le choix des matrices $\omega ^{s,k}$ est exprimée par une action de semi-groupe qui aura une importance considérable dans la preuve du théorème $\ref{cdif}$  et qu'on expose dans la sous-section suivante.
\subsection{Deuxième étape: la notion de recalibration}\label{rec} 
On commence avec les définitions suivantes. On pose par définition
$$
\Gamma(U):={\displaystyle  \bigoplus_{\scriptstyle  s=0,...,m} }GL(p_s,{\cal E}_X(U)) .
$$
On considère $\Gamma(U)$ avec la loi de groupe naturelle induite par les groupes $GL(p_{\bullet},{\cal E}_X(U))$.   
\begin{defi}
La classe $[\varphi,\psi]$ de ${\cal E}$-isomorphisme de la ${\cal E}$-résolution locale $(\varphi ,\psi)$ est l'ensemble des ${\cal E}$-résolutions locales $(\tilde{\varphi} ,\tilde{\psi})$ de longueur $m$ au dessus de l'ouvert $U$ pour lesquelles il existe $g\in\Gamma(U)$ tel que le diagramme suivant soit commutatif.
\begin{diagram}[height=1cm,width=1cm]
0&\rTo &{\cal E}^{\oplus p_m}_{_U} &\rTo^{\tilde{\varphi}_m } &\cdot\cdot\cdot&\rTo^{\tilde{\varphi}_2}  &{\cal E}^{\oplus p_1}_{_U}&\rTo^{\tilde{\varphi}_1}   &{\cal E}_{_U}^{\oplus p_0}&\rTo^{\tilde{\psi}} & {\cal G}_{_{|_U} }&\rTo &0
\\
&       &\dTo^{g_m}_{\wr}&                                         &          &                                &\dTo^{g_1}_{\wr}&          &\dTo^{g_0}_{\wr}&   &\dTo_{\I} 
\\
0&\rTo &{\cal E}^{\oplus p_m}_{_U}  &\rTo^{\varphi_m}       &\cdot\cdot\cdot&\rTo^{\varphi _2}       &{\cal E}^{\oplus p_1}_{_U}  &\rTo^{\varphi _1}  &{\cal E}_{_U}^{\oplus p_0}&\rTo^{\psi}&{\cal G}_{_{|_U} }&\rTo & 0 
\end{diagram} 
\end{defi}
On a alors que $[\varphi,\psi]=\{(\varphi _g,\psi_g)\,|\,g\in\Gamma (U)\} $ où $\psi_g:=\psi\cdot g_0$ et $\varphi _{s,g}:=g_{s-1}^{-1}\cdot \varphi _s \cdot g_s$. Ensuite on désigne par
\begin{eqnarray*} 
\Omega (U,\varphi_g,\psi_g,\bar{\partial})\subset 
 \bigoplus_{(s,k)\in I_m} 
M_{p_{s+k} ,p_s}({\cal E}^{0,k+1}_X(U))
\end{eqnarray*} 
l'ensemble, non vide par l'hypothèse d'exactitude de la ${\cal E} $-résolution locale $(\varphi ,\psi)$, constitué par les éléments 
$\omega =(\omega ^{s,k})_{s,k},\;\omega ^{s,k}\in M_{p_{s+k} ,p_s}({\cal E}^{0,k+1}_X(U))$  tels que 
$ \omega ^{\bullet,-1}=\varphi _{\bullet,g},\;\bar{\partial}\psi_g =\psi_g\cdot \omega ^{0,0}$ et la relation $(\ref{3*})$ soit satisfaite. Ensuite on définit la ``fibration''
$$
\Omega (U,[\varphi,\psi],\bar{\partial}):= {\displaystyle  \coprod_{g\in\Gamma(U) }  }\Omega (U,\varphi_g,\psi_g,\bar{\partial})
$$
au dessus du groupe $\Gamma(U)$, dont les éléments seront appelés $calibrations$ et l'ensemble des paramètres au dessous de l'ouvert $U$
$$
{\cal P}(U)\subset {\displaystyle  \bigoplus_{\scriptstyle  s=0,...,m\atop\scriptstyle  k=0,...,m-s} }M_{p_{s+k} ,p_s}({\cal E}^{0,k}_X(U))
$$
constitué par les éléments $\eta=(\eta^{s,k})_{s,k}$, $ \eta^{s,k}\in M_{p_{s+k} ,p_s}({\cal E}^{0,k}_X(U))$, tels que 
$\I_{p_s}+\eta^{s,0}\in GL(p_s,{\cal E}_X(U))$, $s=0,...,m$. On munit ${\cal P}(U)$ de la loi de semi-groupe donnée par le produit extérieur des paramètres qu'on définit de la façon suivante; si $\eta_1,\;\eta_2\in{\cal P}(U) $ on désigne par $\eta_1\wedge\eta_2\in{\cal P}(U) $ le paramètre dont les composantes sont définies par la formule   
$$
(\eta_1\wedge \eta_2)^{s,k}:= \eta_1^{s,k}+\eta_2^{s,k}+ \displaystyle{\sum_{j=0}^k\eta^{s+j,k-j}_1\wedge \eta^{s,j}_2 } .
$$
L'élément neutre de cette loi est le zéro.  Considérons maintenant l'application
$$
\begin{array}{ccc}
R :{\cal P}(U)\times\Omega (U,[\varphi,\psi],\bar{\partial})& \longrightarrow &\Omega (U,[\varphi,\psi],\bar{\partial})
\\
\\
(\eta,\omega )&\longmapsto &\omega _{\eta}  
\end{array}
$$   
 où les composantes de $\omega _{\eta}$ sont définies par récurrence  sur $k$, pour tous les indices $(s,k)\in I_m$ par la formule (ici on utilise les définitions formelles $\eta^{s,m-s+1}:=0$, $\eta^{-1,j}:=0$ et $\eta^{s,-1}:=0$) 
\begin{eqnarray}\label{4*} 
\omega^{s,k}_{\eta}= \Big(\I_{p_{s+k}} +\eta^{s+k,0}\Big)^{-1}\cdot \Big( \bar{\partial}_{_J }\eta^{s,k}+\sum_{j=0}^{k+1}  \omega^{s+j,k-j}\wedge \eta^{s,j}-\sum_{j=-1}^{k-1} (-1)^{k-j}\eta^{s+j,k-j}\wedge \omega^{s,j}_{\eta} +\omega^{s,k}\Big) .  
\end{eqnarray} 
\begin{figure}[hbtp]
\begin{center} 
\input dgt3.pstex_t
    \caption{}
    \label{fig4}
\end{center} 
\end{figure}    
\\
Le diagramme précédent montre les matrices qui interviennent dans la définition de la matrice
$\omega^{1,2}_{\eta}$ dans le cas où la longueur  de la résolution est égale à 4. 
Le lecteur peut alors essayer avoir une perception visuelle de la formule de recalibration $(\ref{4*}) $. On appelle $R$ application de recalibration et on dit que $ \omega _{\eta} $ est la recalibration de $\omega $ avec paramètre de recalibration $\eta$. Dans la suite on utilisera aussi les notations $g_s\equiv g_s(\eta):=\I_{p_s}+\eta^{s,0}$, $\psi_{\eta}\equiv\psi_{g(\eta)}$. Avec la première notation on a évidemment 
$\omega ^{\bullet,-1}_{\eta} =\omega ^{\bullet,-1}_{g(\eta)}$.       
On remarque aussi que si $\omega \in\Omega (U,\varphi ,\psi,\bar{\partial})$ alors pour tout $\eta \in {\cal P}(U)$, la recalibration 
$R(\eta,\omega )\equiv\omega _{\eta}$ de $\omega $ détermine le diagramme commutatif suivant:
\\
\\
\begin{diagram}[height=1cm,width=1cm]
0&                                       &0&                                                                &0
\\
\uTo&                                    &\uTo&                                                            &\uTo
\\
{\cal G}_{_{|_U} }&\rTo^{\bar{\partial}} &{\cal G}_{_{|_U}}\otimes_{_{{\cal E}_{_V}}}{\cal E}^{0,1}_{_U}
                                                                               &\rTo^{\bar{\partial}} &    {\cal G}_{_{|_U}}
                                                                                                        \otimes_{_{{\cal E}_{_U}}}{\cal E}^{0,2}_{_U} 
\\
\uTo^{\psi_{\eta} }&                             &\uTo_{\psi_{\eta}  \otimes \I_{(0,1)} }&                               &\uTo_{\psi_{\eta}  \otimes \I_{(0,2)} }
\\
{\cal E}^{\oplus p_0} _{_U}&\rTo^{\bar{\partial}_{_J }+\omega^{0,0}_{\eta} }&({\cal E}^{0,1}_{_U})^{\oplus p_0} &\rTo^{\bar{\partial}_{_J }+\omega^{0,0}_{\eta} }&({\cal E}^{0,2}_{_U})^{\oplus p_0} 
\\
\uTo^{\omega^{1,-1}_{\eta} }&                             &\uTo_{\omega^{1,-1}_{\eta} \otimes \I_{(0,1)} }&                               &\uTo_{\omega^{1,-1}_{\eta}  \otimes \I_{(0,2)} }
\\
{\cal E}^{\oplus p_1} _{_U}&\rTo^{\bar{\partial}_{_J }+\omega^{1,0}_{\eta}}&({\cal E}^{0,1}_{_U})^{\oplus p_1}&\rTo^{\bar{\partial}_{_J }+\omega^{1,0}_{\eta} }&({\cal E}^{0,2}_{_U})^{\oplus p_1} 
\\
\uTo&         &\uTo&          &\uTo
\\
\vdots&                                       &\vdots&                                                                &\vdots
\end{diagram}
\\
\\
On a la proposition fondamentale suivante.
\begin{prop}
L'application de recalibration $R$ est bien définie et constitue une action de semi-groupe  transitive sur l'ensemble $\Omega (U,[\varphi,\psi],\bar{\partial})$.
\end{prop}  
\subsection{Troisième étape: introduction à la formulation du problème différentiel}\label{iforpb} 
La partie principale de la preuve consiste à prouver l'existence, pour tout $x\in U$, d'un voisinage ouvert $V \subset U$ de $x$ et $g\in \Gamma(V)$ solution du système différentiel
$$
(\Sigma )\;
\left  \{
\begin{array}{lr}
\bar{\partial}(\psi\cdot g_0)=0 
\\
\\
\bar{\partial}_{_J}(g_{s-1}^{-1} \cdot \varphi _s \cdot g_s)=0  
\\
\\
s=1,...,m . \qquad\qquad\qquad
\end{array}
\right.
$$
Bien évidemment résoudre ce système différentiel équivaut à trouver une autre ${\cal E}$-résolution de  ${\cal G}_{_{|_V}}$ dans la classe $[\varphi ,\psi]$, à partir de la  ${\cal E}$-résolution donnée $(\varphi ,\psi)$, de telle sorte qu'elle admet des matrices de connexion $\omega^{s,0}$  nulles. La fidélité plate de ${\cal E}_{X,x} $ sur ${\cal O}_{X,x} $ permet alors de conclure que le faisceau $Ker\;\bar{\partial}$ est analytique cohérent, (voir la section $\ref{s3.5}$ pour les détails).
\\
On a le lemme suivant.
\begin{lem}
Pour tout choix de calibration $\omega \in \Omega (U,\varphi ,\psi,\bar{\partial})$ on a que l'existence d'une solution $g\in \Gamma(U)$ du système différentiel $(\Sigma)$ est équivalente à l'existence d'une solution $\eta\in {\cal P}(U),\; g=g(\eta)$, du système différentiel quasi-linéaire
$$
(S_{\omega} )\;
\left  \{
\begin{array}{lr}
\displaystyle{\bar{\partial}_{_J }\eta^{s,k}+\sum_{j=0}^{k+1}  \omega^{s+j,k-j}\wedge \eta^{s,j}+(-1)^{k}\eta^{s-1,k+1}\wedge \omega^{s,-1}_{\eta} +\omega^{s,k}=0}
\\
\\
k=0,...,m 
\\
s=0,...,m-k 
\end{array}
\right.
$$ 
\end{lem} 
qui n'est rien d'autre que le système différentiel  
$$
\left  \{
\begin{array}{lr}
\omega^{s,t}_{\eta}=0
\\
s=0,...,m
\\
t=0,...,m-s .
\end{array}                                                                           
\right.
$$
On a besoin de l'hypothèse d'exactitude de la ${\cal E} $-résolution locale seulement pour prouver que les systèmes $\Sigma$ et $(S_{\omega} )$ sont équivalents. A partir du moment où on s'intéresse seulement au système $(S_{\omega} )$ l'hypothèse d'exactitude n'a plus aucun intérêt. On a seulement besoin d'avoir le complexe vertical $(\omega ^{\bullet,-1})$. La notion de recalibration existe encore et elle est une action de semi-groupe. La proposition suivante permet de traduire notre problème en termes purement différentiels.
\begin{prop}
Supposons données des matrices $\omega ^{s,k}\in M_{p_{s+k},p_s }({\cal E}_X^{0,k+1}(U))$, $(s,k)\in I_m$ telles que 
$$
\omega ^{s-1,-1}\cdot\omega ^{s,-1}=0,\;s=2,...,m
$$ et 
$$\bar{\partial}_{_J}\omega ^{s,-1}+\omega ^{s-1,0}\cdot \omega ^{s,-1}= \omega ^{s,-1}\cdot \omega ^{s,0},\;s=1,...,m .
$$
Alors, pour $k\geq 0$, les relations $(\eqref{3*}_{s,k})$
\begin{eqnarray*}
\bar{\partial}_{_J }\omega^{s,k} +\sum_{j=-1}^{k+1} (-1)^{k-j}  \omega^{s+j,k-j}\wedge \omega^{s,j}=0 
\end{eqnarray*}
constituent les conditions d'intégrabilité du système différentiel $(S_{\omega} )$.
\end{prop} 
A partir de maintenant on désignera par
\begin{eqnarray*} 
\Omega (U,m)\subset 
  \bigoplus_{(s,k)\in I_m} 
M_{p_{s+k} ,p_s}({\cal E}^{0,k+1}_X(U))
\end{eqnarray*}  l'ensemble constitué par les éléments $\omega$ dont les composantes vérifient la condition $\omega ^{s-1,-1}\cdot\omega ^{s,-1}=0,\;s=2,...,m$ 
et la relation $(\ref{3*})$.
\subsection{Quatrième étape: introduction au schéma de convergence rapide de type Nash-Moser}\label{intNM} 
\subsubsection{Idée de la preuve de l'estimation fondamentale du schéma de convergence rapide }\label{idest}  
Soit $\omega \in \Omega (B_1, m)$. Pour $r\in (0,1)$ on définit la quantité  
\begin{eqnarray*} 
a_h(\omega ,r):=\max\{\|\omega ^{s,k}\|_{r,h}\,|\,0\leq s \leq m\,,0\leq k \leq m-s\} .
\end{eqnarray*} 
On remarque que, par définition de la norme de Hölder, la quantité $a_h(\omega ,r)$ tend vers zéro lorsque le rayon $r$ tend vers zéro. Pour tout $\sigma \in (0,1)$ on définit les rayons
$$
r_l:=r(1-l\cdot \sigma _m)
$$
pour $l=0,...,m+1$ où on pose par définition $\sigma _m:=\sigma /(m+1)$. A partir de maintenant on désigne par $\varepsilon \in (0,1/2)$ une constante fixée telle que pour toutes les matrices $A\in M_{p_s,p_s}(\C)$ telle que $\|A\|<\varepsilon$ on a l'inversibilité de la matrice $\I_{p_s}+A $. 
\\
On désignera par $S(\omega^{\bullet,-1})$ une suite de poids qui vérifie l'inégalité $ \|\omega ^{\bullet,-1}\|_{1,S(\omega )}<+\infty$.
Venons-en maintenant à la partie essentielle de la preuve du théorème. Avec les notations introduites précédemment on a la proposition suivante.
\begin{prop}
Supposons donnés $\omega \in \Omega (B_1, m),\;r,\,\sigma \in (0,1),\;h\in \N$ et les poids $0<S_j\leq S_j(\omega ),\;j=0,...,h+1$ de la norme de Hölder 
$\|\cdot\|_{r,h+1}$. Supposons que le rayon $r$ soit suffisamment petit pour assurer l'estimation 
$$
L\cdot \sigma _m^{-s(m,h)}\cdot a_h(\omega ,r)<\varepsilon  
$$
où  $L=L(\omega ^{\bullet,-1} )>0$ constante positive et $s(m,h)\in \N$, $(m\geq 0,\,h\geq 0)$ une fonction affine strictement croissante par rapport à la variable $h$ et la quantité $a_h(\omega ,r)$ est calculée par rapport au poids $S_j,\;j=0,...,h$. Supposons de plus que  le poids  $S_{h+1}$ soit suffisamment petit pour pouvoir assurer l'estimation:
\begin{eqnarray*}
S_{h+1} \|\partial^{h+1} \omega^{s,k}_I \|_{r ,\mu } \leq \| \omega^{s,k}_I \|_{r ,\mu }
\end{eqnarray*}
pour tout $k=0,...,m,\;s=0,...,m-k$ et $|I|=k+1$. Il existe alors le paramètre de recalibration $\eta\in {\cal P}(\bar{B} _{r(1-\sigma )})$   dont   les composantes sont définies par la formule de récurrence décroissante sur 
$k=m,...,0$
$$
\eta^{s,k}:=-T_{r_{m-k}} \,\Big(\omega ^{s,k}+\omega ^{s+k+1,-1}\wedge\eta^{s,k+1}+(-1)^k\eta^{s-1,k+1} \wedge\omega ^{s,-1}\Big)
$$ 
tel que les estimations suivantes
\begin{eqnarray}
&\|\eta^{\bullet,k}  \|_{r(1-\sigma ) ,h+1 } \leq L\cdot\sigma_m^{-s(m,h)}\cdot a_h(\omega ,r) ,& \label{9*}   
\\\nonumber
\\
&\|\omega^{\bullet,k}_{\eta} \|_{r(1-\sigma ),h+1} \leq L \cdot\sigma_m^{-s(m,h)}\cdot a_h(\omega ,r) ^2 &\label{10*} 
\end{eqnarray}
soient satisfaites pour tout $k= 0,...,m$.
\end{prop} 
Sans l'hypothèse sur le poids $S_{h+1}$ l'estimation $(\ref{10*} )$ est valable pour $h$ à la place de $h+1$, pour toutes les matrices $\omega ^{\bullet, k}\;k\geq 0$ qui vérifient la relation $(\ref{3*} )$, une fois qu'on a fixé les matrices du complexe vertical $\omega ^{\bullet,-1}$ et les poids
$0<S_j\leq S_j(\omega ),\;j=0,...,h+1$. L'hypothèse sur les poids $S_{h+1}$, comme on verra mieux ensuite, joue un rôle fondamental pour la convergence vers une solution $\ci$ du problème $(S_{\omega } )$. On ''pose'' (voir la preuve de la proposition $\ref{3.4}$ pour l'interprétation correcte!)
$$
\omega^{s,k}_{\eta_{[k+1]}}=\omega ^{s,k}+\omega ^{s+k+1,-1}\wedge\eta^{s,k+1}+(-1)^k\eta^{s-1,k+1} \wedge\omega ^{s,-1}.
$$
On obtient l'estimation $(\ref{10*})$ à l'aide d'une récurrence croissante sur les indices $k\geq 0$. Le problème ici consiste dans le fait qu'on ne peut pas avoir un contrôle quadratique sur le terme $\omega ^{s,-1}_{\eta}$. Précisément les termes qui posent un problème pour obtenir l'estimation quadratique relativement aux matrices $\omega ^{s,k}_{\eta}$ sont les termes qui apparaissent dans la parenthèse de la définition des composantes $\eta^{s,k}$. 
En effet on peut écrire $\omega ^{s,-1}_{\eta}$ sous la forme 
$$
\omega ^{s,-1}_{\eta}=(\I+\theta ^{s-1,0})\cdot\omega ^{s,-1}\cdot (\I +\eta^{s,0})  
$$
avec un contrôle 
$$
\|\theta ^{s,0}\|_{r(1-\sigma ), h+1}\leq 2\|\eta^{s,0}\|_{r(1-\sigma ), h+1}
$$
 sur la norme de la matrice $\theta^{s,0}$. Ensuite en décomposant l'expression de la matrice $\omega ^{s,-1}_{\eta}$ à l'aide de l'expression précédente on arrive à séparer  le terme gênant  
$\eta^{\bullet,\bullet}\wedge \omega ^{\bullet, -1}$ pour les indices adéquats dans la deuxième somme. On a donc:
$$
\omega ^{s,k}_{\eta}= g^{-1}_{s+k}(  \bar{\partial}_{_J }\eta^{s,k}+\text{(termes quadratiques)}+ \omega^{s,k}_{\eta_{[k+1]} }) .
$$ 
On obtient alors en utilisant la formule d'homotopie pour l'opérateur $\bar{\partial}_{_J }$ l'expression:
$$
\omega ^{s,k}_{\eta} =g^{-1}_{s+k}( T_{r_{m-k}}\,\bar{\partial}_{_J }\omega^{s,k}_{\eta_{[k+1]} }+\text{(termes quadratiques)}) .
$$
Le lecteur comprend donc l'exigence d'avoir une estimation quadratique de la norme $h+1$ du terme $T_{r_{m-k}}\,\bar{\partial}_{_J }\omega^{s,k}_{\eta_{[k+1]} } $.
On observe que dans le cas $(s,k)=(0,m),\;(k=m \;\Rightarrow\; s=0)$ on a que le terme $\omega^{s,k}_{\eta_{[k+1]}} $ se réduit à $\omega ^{0,m}$. On remarque que la relation $(\eqref{3*} _{0,m})$ est la seule, parmi les autres relations $(\eqref{3*} _{\bullet,\bullet})$, qui ne présente pas de facteurs de type $\omega^{\bullet,-1} $ dans les termes quadratiques. 
Le diagramme suivant montre les matrices qui interviennent dans la relation $(\eqref{3*} _{0,m})$.
\begin{figure}[hbtp]
\begin{center} 
\input dgt6.pstex_t
    \caption{}
\end{center} 
\end{figure}    
\\
\\
On déduit alors à l'aide de l'inégalité $(\ref{Webineq})$, l'estimation quadratique de la norme $h+1$ du terme $T_r\, \bar{\partial}_{_J }\omega ^{0,m}$. 
Une récurrence décroissante sur les indices $k=m,...,0$ combinée avec le fait que les matrices $\omega$ vérifient la condition d'intégrabilité $(\ref{3*} )$ montre l'estimation quadratique 
\begin{eqnarray*}
\|T_{r_{m-k}}\,\bar{\partial}_{_J }\omega^{s,k}_{\eta_{[k+1]} }\|_{r(1-\sigma ),h+1}\leq Q\cdot\sigma_m^{-s(m,h)} \cdot a_h^2 
\end{eqnarray*}
pour tout $k=0,...,m$, (où $Q>0$ est une constante positive). 
\\
On explique  maintenant où intervient  l'hypothèse sur le poids $S_{h+1}$ dans la preuve de l'estimation $(\ref{10*})$. Pour obtenir celle-ci  on doit estimer les normes du type $\|\omega ^{\bullet,\bullet}\wedge \eta^{\bullet,\bullet}\|_{r(1-\sigma ),\, h+1}$. L'hypothèse faite sur le poids $S_{h+1}$ nous permet  d'obtenir l'estimation
$$
\|\omega ^{\bullet,\bullet}\wedge \eta^{\bullet,\bullet}\|_{r(1-\sigma ),\, h+1}\leq 2\|\omega ^{\bullet,\bullet}\|_{r,\, h}\cdot\| \eta^{\bullet,\bullet}\|_{r(1-\sigma ),\, h+1} .
$$
Voici certaines des idées principales de la preuve de la proposition.
\subsubsection{Esquisse du procédé itératif}\label{idprocit} 
Les calibrations $\omega _k\in \Omega (\bar{B}_{r_k},m)$ obtenues au $k$-ième pas du procédé itératif 
sont définies par la formule récursive $\omega _{k+1}:=\omega _{k,\,\eta_{k+1} }$
où $r_{k+1}:=r_k(1-\sigma _k)$ et où $\sigma _k\in (0,1)$ est un paramètre qui contrôle la décroissance des rayons des boules, (le rayon initial $r_0$ étant choisi suffisamment petit). On choisit les quantités $\sigma _k\in (0,1)$ de telle sorte que la série $\sum \sigma _k$ soit convergente. Le rayon limite $r_{\infty}:=\lim_{k\rightarrow +\infty}r_k$ est alors non nul. 
Le paramètre  $\eta_{k+1} \in {\cal P}(\bar{B}_{r_{k+1}})$
qui contrôle la recalibration des éléments $\omega _k$, $k\geq 0$ au $k$-ième pas du procédé itératif est défini de façon analogue à celle de la proposition précédente en fonction de la calibration $\omega _k$. 
Les poids sont choisis de telle sorte que les estimations suivantes soient satisfaites pour tout  entier $k\geq 0$ et $t=0,...,m$.
\begin{eqnarray}
&S_{k+1} \|\partial^{k+1} \omega^{\bullet,t}_{k,\,I}  \|_{r_k ,\,\mu }\leq\| \omega^{\bullet,t}_{k,\,I}  \|_{r_k ,\,\mu } , &\label{15*} 
\\\nonumber
\\
&S_{k+1} \|\partial^{k+1} g_{\bullet}(k)^{\pm 1}  \|_{r_k ,\,\mu }\leq 2^{-k-1}\|g_{\bullet}(k)^{\pm 1}   \|_{r_k ,\,\mu }& \label{16*} 
\end{eqnarray}
où 
$$
g_{s}(k) :={\displaystyle \prod_{0\leq j\leq k }^{\longrightarrow}}g_{s,j} 
$$
(ici on rappelle qu'on pose par définition $g_0:=\I_{p_0} $ et que le symbole de produit avec une flèche vers la droite désigne le produit non commutatif de termes qui sont écrits en ordre croissant de l'indice vers la droite).
Cette dernière inégalité sert à assurer un bon fonctionnement du procédé itératif, plus précisément elle permet d'appliquer la proposition précédente à toutes les étapes du procédé. On pose par définition
$$
a_{k} :=\max\{\|\omega ^{s,t}_k\|_{r_k,\,k}\,|\,0\leq s \leq m\,,0\leq t \leq m-s\} 
$$
et 
$$
b_k:=H\cdot\sigma_{m,\,k}  ^{- s(m,k) }\cdot a_k \qquad \sigma_{m,\,k}  :=\sigma _k/(m+1)
$$
Avec les notations introduites précédemment on a la proposition suivante.
\begin{prop} 
Pour tout entier $k\geq 0$ on a les estimations suivantes;
\begin{eqnarray}
&a_{k+1} \leq H\cdot\sigma_{m,\,k} ^{-s(m,k)}\cdot a^2_{k}\leq 1 ,&
\\\nonumber
\\
&\|\eta^{\bullet,t}_{k+1} \|_{r_{k+1}  ,\,k+1 } \leq  b_k   <\varepsilon <1/2 ,&
\\\nonumber
\\
&\|\omega ^{s,-1} _{k+1}\|_{r_{k+1},k+1 } \leq 4 \|\omega ^{\bullet,-1}\|_{1,S(\omega ^{\bullet,-1}) } & \label{19*}
\end{eqnarray}
et les quantités $a_k,\;b_k$ tendent $($avec la bonne vitesse$)$ vers zéro lorsque $k$ tend vers plus l'infini.
\end{prop}
L'estimation $(\ref{19*})$  montre que la norme du complexe vertical n'explose pas. La condition $(\ref{16*})$ sert aussi à assurer cette inégalité.
Si on pose par définition
$$
\eta(k):={\displaystyle \bigwedge_{1\leq j\leq k }^{\longrightarrow}}\eta_j 
$$
(ici aussi le symbole de produit extérieur avec une flèche vers la droite désigne le produit non commutatif de termes qui sont écrits en ordre croissant de l'indice vers la droite), le paramètre de recalibration sur la boule $B_{r_{\infty}}$ on aura la formule  $\omega_k=\omega_{\eta(k)}$, grâce au fait que la recalibration $R$ est une action.  On a alors la proposition suivante.
\begin{prop} La limite  
$$
\eta:=\lim_{k\rightarrow \infty} \eta (k)
$$
existe en topologie ${\cal C}^{h,\mu}$ pour tout $h\geq 0$ et constitue  un paramètre de recalibration $\eta\in {\cal P}(B_{r_{\infty}})$, solution du problème différentiel $(S_{\omega})$.
\end{prop} 
Le fait que $\eta$ constitue une solution 
(de classe ${\cal C}^{\infty}$) pour le problème différentiel $(S_{\omega} )$ est clair. En effet 
\begin{eqnarray*}
\omega ^{s,t}_{\eta}=\lim_{k\rightarrow \infty}\omega ^{s,t}_k
\end{eqnarray*} 
et
\begin{eqnarray*}
\|\omega ^{s,t}_{\eta} \|_{r_{\infty} ,0}=\lim_{k\rightarrow \infty}\|\omega ^{s,t}_k\|_{r_{\infty} ,0} \leq\lim_{k\rightarrow \infty}a_k=0  
\end{eqnarray*} 
ce qui montre que $\eta \in {\cal P}(B_{r_{\infty} })$ est une solution du  système différentiel $(S_{\omega} )$. 
\section{La preuve complète du théorème $\ref{cdif}$  dans le cas général de longueur arbitraire de la ${\cal E}$-résolution locale}
\subsection{Première étape: Preuve de l'expression locale de la condition d'intégrabilité $\bar{\partial}^2=0$} 
On commence par rappeler la définition de l'ensemble d'indices 
$$
I_m:=\{(s,k)\;|\;s=0,...,m,\; k=-1,...,m-s,\;(s,k)\not= (0,-1) \} .
$$
Avec cette notation on a le lemme élémentaire suivant.
\begin{lem} \label{3.1.1} 
Soit X une variété complexe et soit ${\cal G}$ un faisceau de  ${\cal E}_X$-modules.
\\
$(A)$ Supposons que le faisceau ${\cal G}$ admet des ${\cal E}$-présentations locales et soit
$$
 {\cal E}_{_U}^{\oplus p_1 }\stackrel{\varphi}{ \longrightarrow }{\cal E}_{_U}^{\oplus p_0 }\stackrel{\psi}{\longrightarrow}
{\cal G}_{_{|_U} }\rightarrow 0
$$
une ${\cal E}$-présentation au dessus d'un ouvert $U$. Alors l'existence d'une connexion 
$\bar{\partial}$ de type $(0,1)$ sur le faisceau ${\cal G}_{_{|_U} }$ telle que 
$\bar{\partial}^{2}=0$, implique l'existence de matrices $\omega ^{s,0}\in M_{p_s,p_s}({\cal E}^{0,1}_X(U)),\;s=0,1$ et 
$\omega ^{0,1}\in M_{p_1,p_0}({\cal E}^{0,2}_X(U))$ telles que $\bar{\partial}\psi=\psi\cdot \omega^{0,0}$ et
\begin{eqnarray}
&\bar{\partial}_{_J }\varphi+\omega ^{0,0}\cdot \varphi =  \varphi\cdot \omega ^{1,0} ,& \label{1}       
\\\nonumber
\\
&\bar{\partial}_{_J }\omega ^{0,0} +\omega ^{0,0}\wedge \omega ^{0,0}=\varphi \cdot\omega ^{0,1} .&  \label{2}  
\end{eqnarray}
Réciproquement l'existence des matrices $\omega ^{s,0},\;s=0,1 $ et $ \omega ^{0,1} $ qui vérifient les relations $(\ref{1} )$ et $(\ref{2})$, implique l'existence d'une connexion $\bar{\partial}$ de type $(0,1)$ sur le faisceau ${\cal G}_{_{|_U} }$ telle que 
$\bar{\partial}\psi=\psi\cdot \omega^{0,0}$ et
$\bar{\partial}^{2}=0$. 
\\
$(B)$ Supposons que le faisceau ${\cal G}$ admet des ${\cal E}$-résolutions locales de longueur finie et soit 
$$
0\rightarrow {\cal E}_{_U}^{\oplus p_m }\stackrel{\varphi_m }{\longrightarrow}  
{\cal E}_{_U}^{\oplus p_{m-1} }\stackrel{ \varphi _{m-1} }{\longrightarrow} \cdot \cdot\cdot
\stackrel{\varphi_2}{ \longrightarrow } {\cal E}_{_U}^{\oplus p_1 }\stackrel{\varphi_1}{ \longrightarrow }{\cal E}_{_U}^{\oplus p_0 }\stackrel{\psi}{\longrightarrow}
{\cal G}_{_{|_U} }\rightarrow 0
$$
une telle ${\cal E}$-résolution. Alors l'existence d'une connexion 
$\bar{\partial}$ de type $(0,1)$ sur le faisceau ${\cal G}_{_{|_U} }$ telle que 
$\bar{\partial}^{2}=0$, implique l'existence des matrices $\omega^{s,k}  \in M_{p_{s+k} ,p_s}({\cal E}^{0,k+1}_X(U))$ pour $(s,k)\in I_m$ telles que si on utilise l'identification $\varphi _s\equiv\omega ^{s,-1}$ et les conventions formelles $\omega ^{0,-1} :=0,\;\omega ^{-1,j}:=0$ et $\omega ^{s,k}:=0$ si $s\geq m+1$ ou $k\geq m-s+1$, on aura les relations $\bar{\partial}\psi=\psi\cdot \omega^{0,0}$ et
\begin{eqnarray}
\bar{\partial}_{_J }\omega^{s,k} +\displaystyle{\sum_{j=-1}^{k+1} }(-1)^{k-j}  \omega^{s+j,k-j}\wedge \omega^{s,j}=0 \label{3} 
\end{eqnarray}
pour tout $(s,k)\in I_m$.
\\
Réciproquement l'existence des matrices $\omega ^{s,k}$ qui vérifient la relation $(\ref{3})$ implique l'existence d'une connexion  $\bar{\partial}$ de type $(0,1)$ sur le faisceau ${\cal G}_{_{|_U} }$ telle que 
$\bar{\partial}\psi=\psi\cdot \omega^{0,0}$ et
$\bar{\partial}^{2}=0$. 
\end{lem} 
$Preuve\; de\;(A) $. La preuve de $(A)$ a été donnée dans la sous-section $(\ref{explocm1})$. On rappelle qu'on obtient le diagramme commutatif suivant, ayant des flèches verticales exactes:
\\
\\
\begin{diagram}[height=1cm,width=1cm]
0&                                       &0&                                                                &0
\\
\uTo&                                    &\uTo&                                                            &\uTo
\\
{\cal G}_{_{|_U} }&\rTo^{\bar{\partial}} &{\cal G}_{_{|_U}}\otimes_{_{{\cal E}_{_U}}}{\cal E}^{0,1}_{_U}
                                                                               &\rTo^{\bar{\partial}} &    {\cal G}_{_{|_U}}
                                                                                                        \otimes_{_{{\cal E}_{_U}}}{\cal E}^{0,2}_{_U} 
\\
\uTo^{\psi}&                             &\uTo_{\psi \otimes \I_{(0,1)} }&                               &\uTo_{\psi \otimes \I_{(0,2)} }
\\
{\cal E}^{\oplus p_0} _{_U}&\rTo^{\bar{\partial}_{_J }+\omega ^{0,0} }&({\cal E}^{0,1}_{_U})^{\oplus p_0} &\rTo^{\bar{\partial}_{_J }+\omega ^{0,0} }&({\cal E}^{0,2}_{_U})^{\oplus p_0} 
\\
\uTo^{\varphi_1}&                             &\uTo_{\varphi_1\otimes \I_{(0,1)} }&                               &\uTo_{\varphi_1 \otimes \I_{(0,2)} }
\\
{\cal E}^{\oplus p_1} _{_U}&\rTo^{\bar{\partial}_{_J }+\omega ^{1,0} }&({\cal E}^{0,1}_{_U})^{\oplus p_1}&\rTo^{\bar{\partial}_{_J }+\omega ^{1,0} }&({\cal E}^{0,2}_{_U})^{\oplus p_1} 
\end{diagram} 
\\
\\
$Preuve\;de\;(B)$. Pour $(s,k)=(1,-1)$ et $(s,k)=(0,0)$ la relation  $\eqref{3}$ exprime les relations $\eqref{1}$ et $\eqref{2}$ de la partie (A) du lemme. En  appliquant l'opérateur $\bar{\partial}_{_J }$ aux identités $\varphi_{t-1} \circ \varphi_t=0$ on obtient inductivement, de la même façon utilisée pour obtenir la relation $\eqref{1}$, l'existence d'une matrice $\omega ^{s,0} \in M_{p_s,p_s}({\cal E}^{0,1}_X(U))$, $s=1,...,m$ telle que 
$$
\bar{\partial}_{_J }\varphi _s +\omega ^{t-1,0} \cdot  \varphi_s =  \varphi_s\cdot \omega ^{t,0}  .
$$
Ces relations constituent les relations $\eqref{3}$ pour $(s,-1),\;s=1,...,m$. On va montrer maintenant l'existence des matrices $\omega ^{s,k},\;k\geq 1 $ qui vérifient la relation $\eqref{3}$ à l'aide du procédé récursif triangulaire suivant. Pour un couple $(s,k)$, $s=0,...,m-1$ et $k=1,...,m-s$ on suppose avoir déjà défini $\omega^{\sigma ,\kappa}$ pour $\sigma +\kappa \leq s+k$, $\kappa \leq k+1$, et on applique l'opérateur $\bar{\partial}_{_J }$ à l'expression $(\eqref{3} _{s,k-1})$ pour $k\geq 1$. On obtient alors la relation suivante:
\begin{eqnarray*}
 \sum_{j=-1}^{k} (-1)^{k-j-1}\bar{\partial}_{_J }  \omega^{s+j,k-j-1}\wedge \omega^{s,j}- \sum_{j=-1}^{k}  \omega^{s+j,k-j-1}\wedge \bar{\partial}_{_J } \omega^{s,j} =0 .
\end{eqnarray*}
En explicitant les termes $\bar{\partial}_{_J }  \omega^{\bullet,\bullet}$ dans la relation précédente (qui bien évidemment, grâce à l'hypothèse de récurrence précédente, vérifient la relation $(\eqref{3}_{\bullet,\bullet})$ pour les indices voulus) on obtient:
\begin{eqnarray*}
\sum_{j=-1}^{k} (-1)^{k-j-1}  \omega^{s+k,-1}\wedge \omega^{s+j,k-j}\wedge \omega^{s,j} + \qquad\qquad\qquad
\\
\\
+\sum_{j=-1}^{k} \sum_{r=-1}^{k-j-1}(-1)^{r+1} 
\omega^{s+j+r,k-j-1-r}\wedge \omega^{s+j,r} \wedge \omega^{s,j}+\qquad\qquad
\\
\\
+ \sum_{j=-1}^{k-1} \sum_{r=-1}^{j+1}(-1)^{j-r} 
\omega^{s+j,k-j-1}\wedge \omega^{s+r,j-r} \wedge \omega^{s,r}-\omega^{s+k,-1}\wedge\bar{\partial}_{_J }  \omega^{s,k}  =0 .
\end{eqnarray*}
En faisant le changement d'indice $j'=j+r$, $r'=j$ dans la deuxième somme et en rappelant que $\omega^{s-1,-1} \wedge \omega^{s,-1}=0$ on obtient: 
\begin{eqnarray*}
\omega^{s+k,-1}\wedge\big(\bar{\partial}_{_J }\omega^{s,k}+\sum_{j=-1}^{k} (-1)^{k-j} \omega^{s+j,k-j}\wedge \omega^{s,j}\big)=0  .
\end{eqnarray*}
L'hypothèse d'exactitude  nous permet de choisir $\omega^{s,k+1}$ telle que la relation
\begin{eqnarray*}
\bar{\partial}_{_J }\omega^{s,k}+\sum_{j=-1}^{k} (-1)^{k-j} \omega^{s+j,k-j}\wedge \omega^{s,j}=\omega^{s+k+1,-1}\wedge\omega^{s,k+1}  
\end{eqnarray*}
soit satisfaite.
 Ce type de récurrence peut se visualiser grâce au tableau de la figure $\ref{fig1}$.
\begin{figure}[hbtp]
\begin{center} 
\input dgt0.pstex_t       
   \caption{}
    \label{fig1}
\end{center} 
\end{figure}    
L'hypothèse de finitude de la longueur des ${\cal E}$-résolutions locales de ${\cal G}$ permet d'arrêter  ce procédé  après un nombre fini d'étapes. On a donc prouvé la première implication de la partie (B) du lemme. Le réciproque dans la partie (B) est évidemment une conséquence banale de la partie (A) du lemme.\hfill $\Box$
\\
Les relations $(\ref{3})$ pour les indices $k\geq 0$ constituent les expressions locales de la condition d'intégrabilité $\bar{\partial}^2=0$ de la connexion relativement à la ${\cal E}$-résolution locale choisie. Les relations $(\ref{3})$ pour les indices $(s,k)=(\bullet,-1)$ représentent simplement des identités de commutation.
\\
Si on désigne par ${\cal P}_0(U)\subset {\cal P}(U)$ le sous-semigroupe des paramètres tels que $\eta^{\bullet,0}=0$ on a que la restriction de la recalibration 
$R_0 :{\cal P}_0(U)\times\Omega (U,[\varphi,\psi],\bar{\partial})\longrightarrow \Omega (U,[\varphi,\psi],\bar{\partial})$ est une application fibrée qui contrôle la liberté homologique qui caractérise le choix des matrices $\omega^{\bullet,\bullet}$ relativement aux ${\cal E} $-résolutions locales $(\varphi _g,\psi_g)$, $g\in \Gamma(U)$. 
\\
{\bf{Remarque.} } A partir de maintenant le lecteur doit tenir compte du fait que certains des calculs et formules qui suivront n'existent pas dans le cas de longueur $m=0$ de la ${\cal E}$-résolutions locale, autrement dit dans le cas des faisceaux localement libres. Cependant les calculs qui survivent ont encore  sens et ils font partie de notre preuve (différente de la preuve donnée par Grothendieck-Koszul-Malgrange) dans ce cas. On rappelle aussi qu'on utilise la convention qui consiste à négliger les termes d'une somme ou d'un produit si l'ensemble des indices sur lesquels on effectue ces opérations est vide. 
\newpage
\subsection{Deuxième étape: le formalisme du procédé itératif}\label{ss3.2} 
On commence par prouver la proposition suivante.
\begin{prop}\label{3.2} 
L'application de recalibration $R$ est bien définie et constitue une action de semi-groupe  transitive sur l'ensemble $\Omega (U,[\varphi,\psi],\bar{\partial})$
\end{prop}  
$Preuve$. Nous commençons par prouver que l'application $R $ est bien définie. 
On prouve d'abord les relations $\bar{\partial}\psi_{\eta} =\psi_{\eta}\cdot\omega^{0,0}_{\eta}$. En effet on a :
\begin{eqnarray*}
\bar{\partial}\psi_{\eta} =\bar{\partial}\psi\cdot g_0+ \psi \cdot \bar{\partial}\eta^{0,0} =
\psi\cdot(\bar{\partial}\eta^{0,0}+\omega ^{0,0}\wedge \eta^{0,0}+ \omega ^{0,0})=
\\
\\
=\psi\cdot(\bar{\partial}\eta^{0,0}+\omega ^{0,0}\wedge \eta^{0,0}+\omega ^{1,-1}\wedge \eta^{0,1} + \omega ^{0,0})=\psi_{\eta}\cdot \omega ^{0,0}_{\eta} .\;\;
\end{eqnarray*}
On prouve maintenant que les matrices $\omega^{\bullet,\bullet}_{\eta} $ vérifient la relation $\eqref{3}$ pour l'indice $k=-1$, autrement dit on veut montrer la relation:
$$
\bar{\partial}_{_J } \omega ^{s,-1}_{\eta} + \omega ^{s-1,0}_{\eta}\cdot \omega ^{s,-1}_{\eta} =\omega ^{s,-1}_{\eta}\cdot\omega ^{s,0}_{\eta} .
$$
On commence par développer le terme $\bar{\partial}\omega ^{s,-1}_{\eta} $, en utilisant la relation $\eqref{3}$ pour l'indice $k=-1$, relativement aux matrices $\omega^{\bullet,\bullet}$. On obtient alors les égalités suivantes:
\begin{eqnarray*}
&\bar{\partial}_{_J } \omega ^{s,-1}_{\eta}=-g_{s-1}^{-1}\cdot\bar{\partial}_{_J } g_{s-1}\cdot g_{s-1}^{-1}\cdot\omega ^{s,-1} \cdot g_s+&
\\
\\
&+g_{s-1}^{-1}\cdot  \bar{\partial}_{_J } \omega ^{s,-1}\cdot g_s+  g_{s-1}^{-1}\cdot   \omega ^{s,-1}\cdot \bar{\partial}_{_J }g_s=&
\end{eqnarray*}
\begin{eqnarray*}
&=-g_{s-1}^{-1}(\bar{\partial}_{_J }\eta^{s-1,0}+\omega ^{s-1,0}\wedge \eta^{s-1,0}+\omega ^{s-1,0})\cdot\omega ^{s,-1}_{\eta} +&
\\
\\
&+\omega ^{s,-1}_{\eta} \cdot g_s^{-1}(\bar{\partial}_{_J }\eta^{s,0}+\omega ^{s,0}\wedge \eta^{s,0}+\omega ^{s,0}) .&
\end{eqnarray*}
En rappelant que $\omega ^{s-1,-1}\cdot \omega ^{s,-1}=0$ et en rajoutant et en soustrayant le terme 
$-g_{s-1}^{-1}\cdot \omega ^{s,-1}\cdot\eta^{s-1,1}\cdot \omega ^{s,-1}_{\eta}$ à la dernière expression de $\bar{\partial}\omega ^{s,-1}_{\eta} $, on obtient:
\begin{eqnarray*}
\bar{\partial}_{_J } \omega ^{s,-1}_{\eta}=-\omega ^{s-1,0}_{\eta}\cdot \omega ^{s,-1}_{\eta}+ \qquad \qquad \qquad 
\\
\\
+\omega ^{s,-1}_{\eta}\cdot g_s^{-1}(\bar{\partial}_{_J }\eta^{s,0}+\omega ^{s,0}\wedge \eta^{s,0}+\eta^{s-1,1}\wedge \omega ^{s,-1}_{\eta}+ \omega ^{s,0})= 
\\
\\
=-\omega ^{s-1,0}_{\eta}\cdot \omega ^{s,-1}_{\eta} + \omega ^{s,-1}_{\eta}\cdot\omega ^{s,0}_{\eta} . \qquad \qquad\qquad
\end{eqnarray*}
On va montrer maintenant la validité de la formule $\eqref{3}$ pour tous les indices, relativement aux matrices $\omega^{\bullet,\bullet}_{\eta}$, avec un procédé récursif analogue à celui qui nous a permis de définir les matrices $\omega^{\bullet,\bullet}$. Voici les détails de la récurrence. 
Pour un couple $(s,k)$, $s=0,...,m$, $k=0,...,m-s$ on suppose avoir déjà montré la relation
$$ 
\qquad \qquad \qquad \qquad \qquad \qquad \bar{\partial}_{_J }\omega^{\sigma ,\kappa}_{\eta} +\sum_{j=-1}^{\kappa+1} (-1)^{\kappa-j}  \omega^{\sigma +j,\kappa-j}_{\eta} \wedge \omega^{\sigma ,j}_{\eta}=0  \qquad \qquad \qquad \qquad   (\eqref{3}^{\sigma ,\kappa}_{\eta}) 
$$
pour $\sigma=s$, $\kappa=-1,...,k-1$. En développant le terme en $\bar{\partial}_{_J }$ de l'expression suivante on a l'égalité:
\begin{eqnarray*}
\bar{\partial}_{_J }\omega^{s,k}_{\eta} +\sum_{j=-1}^{k} (-1)^{k-j}  \omega^{s+j,k-j}_{\eta}\wedge \omega^{s,j}_{\eta}= \qquad\qquad\qquad\qquad\qquad
\\
\\
=g_{s+k}^{-1}\Big(-\bar{\partial}_{_J } \eta^{s+k,0}\wedge\omega^{s,k}_{\eta}+
\sum_{j=0}^{k+1}\bar{\partial}_{_J }\omega^{s+j,k-j}\wedge \eta^{s,j}- \sum_{j=0}^{k+1}(-1)^{k-j}\omega^{s+j,k-j}\wedge\bar{\partial}_{_J } \eta^{s,j}-
\\
\\
- \sum_{j=-1}^{k-1}(-1)^{k-j}\bar{\partial}_{_J }\eta^{s+j,k-j}\wedge\omega^{s,j}_{\eta}- \sum_{j=-1}^{k-1}\eta^{s+j,k-j}\wedge \bar{\partial}_{_J }\omega^{s,j}_{\eta}+\bar{\partial}_{_J }\omega^{s,k} \Big)+\qquad
\\
\\
+\sum_{j=-1}^{k} (-1)^{k-j}  \omega^{s+j,k-j}_{\eta}\wedge \omega^{s,j}_{\eta}=(A_1) .\qquad\qquad\qquad\qquad\qquad
\end{eqnarray*}
En développant les termes $\bar{\partial}_{_J }\omega^{s+j,k-j}$ et $\bar{\partial}_{_J }\omega^{s,j}_{\eta}$ à l'aide respectivement des expressions $\eqref{3}$ et $(\eqref{3}^{s ,j}_{\eta})$ on obtient:
\begin{eqnarray*}
(A_1)=g_{s+k}^{-1}  \omega^{s+k+1,-1}\Big(\bar{\partial}_{_J } \eta^{s,k+1} +\sum_{j=0}^{k+1}\omega^{s+j,k+1-j}\wedge\eta^{s,j}\Big)+\qquad\qquad\qquad\quad
\\
\\
+g_{s+k}^{-1} \Big(\sum_{j=0}^{k+1}\sum_{r=-1}^{k-j}(-1)^{k-j-r+1}\omega^{s+j+r,k-j-r}\wedge\omega^{s+j,r}\wedge\eta^{s,j}- \sum_{j=0}^{k}(-1)^{k-j}\omega^{s+j,k-j}\wedge\bar{\partial}_{_J } \eta^{s,j}-
\\
\\
-\sum_{j=-1}^{k}(-1)^{k-j}\bar{\partial}_{_J }\eta^{s+j,k-j}\wedge\omega^{s,j}_{\eta}    + \sum_{j=-1}^{k-1}\sum_{r=-1}^{j+1}(-1)^{j-r}\eta^{s+j,k-j}\wedge \omega^{s+r,j-r}_{\eta}\wedge \omega^{s,r}_{\eta} +\bar{\partial}_{_J }\omega^{s,k}\Big)+
\end{eqnarray*}
\begin{eqnarray*}
+\sum_{j=-1}^{k} (-1)^{k-j}  \omega^{s+j,k-j}_{\eta}\wedge \omega^{s,j}_{\eta}=(A_2) .
\end{eqnarray*}
En faisant le changement d'indice $j'=j+r$, $r'=j$ dans la première somme double et en développant les premiers facteurs $\omega^{s+j,k-j}_{\eta}$ de la dernière somme on a:
\begin{eqnarray*}
(A_2)=g_{s+k}^{-1}\omega^{s+k+1,-1}\Big(\bar{\partial}_{_J } \eta^{s,k+1} +\sum_{j=0}^{k+1}\omega^{s+j,k+1-j}\wedge\eta^{s,j}\Big)+
g_{s+k}^{-1}\bar{\partial}_{_J }\omega^{s,k} 
\\
\\
+\sum_{j=-1}^{k}(-1)^{k-j+1}g_{s+k}^{-1}\omega^{s+j,k-j}\wedge\Big(\bar{\partial}_{_J } \eta^{s,j} + \sum_{r=0}^{j+1}\omega^{s+r,j-r}\wedge\eta^{s,r}\Big)+
\end{eqnarray*}
\begin{eqnarray*}
+\sum_{j=-1}^{k} (-1)^{k-j}g_{s+k}^{-1} \Big(\sum_{r=0}^{k-j+1}\omega^{s+j+r,k-j-r}\wedge\eta^{s+j,r}+\omega^{s+j,k-j}\Big)\wedge \omega^{s,j}_{\eta}=(A_3) .
\end{eqnarray*}
En rappelant que $g_{s+j,0}:=\I_{s+k}+\eta^{s+j,0}$, en décomposant les termes extrêmes de la somme $\sum_{r=0}^{k-j+1}$ et en décomposant les facteurs $\omega^{s,j}_{\eta}$ qui apparaissent dans les produits $\omega^{s+j,k-j}\wedge\omega^{s,j}_{\eta}$ on obtient les égalités suivantes:
\begin{eqnarray*}
(A_3)=g_{s+k}^{-1}\omega^{s+k+1,-1}\Big(\bar{\partial}_{_J } \eta^{s,k+1} +\sum_{j=0}^{k+1}\omega^{s+j,k+1-j}\wedge\eta^{s,j} -\sum_{j=-1}^{k} (-1)^{k+1-j}\eta^{s+j,k+1-j}\wedge \omega^{s,j}_{\eta}\Big)+ 
\\
\\
+g_{s+k}^{-1}\bar{\partial}_{_J }\omega^{s,k}+\sum_{j=-1}^{k} (-1)^{k-j}g_{s+k}^{-1}\omega^{s+j,k-j}\wedge\Big(\sum_{r=-1}^{j-1}(-1)^{j-r+1}\eta^{s+r,j-r}\wedge\omega^{s,r}_{\eta}+\omega^{s,j}\Big)\qquad
\\
\\
+\sum_{j=-1}^{k}\sum_{r=1}^{k-j} (-1)^{k-j}g_{s+k}^{-1}\omega^{s+j+r,k-j-r}\wedge\eta^{s+j,r}\wedge \omega^{s,j}_{\eta}=(A_4) .\qquad\qquad\qquad\quad
\end{eqnarray*}
En faisant le changement d'indice $j'=j+r$, $r'=j$ dans la dernière somme double on a finalement:
\begin{eqnarray*}
(A_4)=g_{s+k}^{-1}\omega^{s+k+1,-1}\Big(\bar{\partial}_{_J } \eta^{s,k+1} +\sum_{j=0}^{k+1}\omega^{s+j,k+1-j}\wedge\eta^{s,j} -\sum_{j=-1}^{k} (-1)^{k+1-j}\eta^{s+j,k+1-j}\wedge \omega^{s,j}_{\eta}\Big)+ 
\\
\\
+g_{s+k}^{-1}\Big(\bar{\partial}_{_J }\omega^{s,k}+\sum_{j=-1}^{k} (-1)^{k-j}\omega^{s+j,k-j}\wedge\omega^{s,j}\Big)=
\omega^{s+k+1,-1}_{\eta}\wedge \omega^{s,k+1}_{\eta}\qquad\qquad\quad
\end{eqnarray*}
ce qui justifie la formule $(\eqref{3} ^{s,k}_{\eta})$. On a alors qu'à la fin de cette récurrence toutes les matrices $\omega^{s,k}_{\eta} $ vérifient la relation $(\eqref{3} _{\eta})$.
Montrons maintenant que l'application $R$ est une action de semi-groupe. On se propose donc de montrer la formule $R(\eta_2,\omega_{\eta_1})=: \omega_{\eta_1,\eta_2}=\omega_{\eta_1\wedge\eta_2}$ qui en termes de composantes s'exprime sous la forme 
$\omega^{s,k}_{\eta_1,\eta_2}=\omega^{s,k}_{\eta_1\wedge\eta_2} $ pour tout $k\geq -1$. On montre la formule précédente par récurrence sur $k$. On remarque que la formule est évidente pour $k=-1$. En explicitant l'expression de $\omega^{s,k}_{\eta_1,\eta_2}$ et en utilisant l'hypothèse de récurrence on a:
\begin{eqnarray*}
\omega^{s,k}_{\eta_1,\eta_2}= g_{s+k,2}^{-1}\Big(\bar{\partial}_{_J }\eta^{s,k}_2+\sum_{j=1}^{k+1}  \omega^{s+j,k-j}_{\eta_1}\wedge \eta^{s,j}_2-\sum_{j=-1}^{k-1} (-1)^{k-j}\eta^{s+j,k-j}_2\wedge \omega^{s,j}_{\eta_1\wedge \eta_2} +\omega^{s,k}_{\eta_1}\cdot g_{s,2}\Big) =
\\
\\
=(g_{s+k,1}\cdot g_{s+k,2})^{-1}\Big[g_{s+k,1}\cdot \bar{\partial}_{_J }\eta^{s,k}_2+\sum_{j=1}^{k+1}\bar{\partial}_{_J } \eta^{s+j,k-j}_1\wedge \eta^{s,j}_2 +\qquad\qquad\qquad\quad
\\
\\
\underbrace{+\sum_{j=1}^{k+1}\sum_{r=1}^{k-j+1}\omega^{s+j+r,k-j-r}\wedge\eta^{s+j,r}_1\wedge\eta^{s,j}_2+\sum_{j=1}^{k+1} \omega^{s+j,k-j}\wedge g_{s+j,1}\cdot\eta^{s,j}_2-}_{(1)} \qquad\qquad\quad
\\
\\
\underbrace{ -\sum_{j=1}^{k+1}\sum_{r=-1}^{k-j-1}(-1)^{k-j-r} \eta^{s+j+r,k-j-r}_1\wedge\omega^{s+j,r}_{\eta_1}\wedge\eta^{s,j}_2}_{(2)}
-\sum_{j=-1}^{k-1}(-1)^{k-j}g_{s+k,1}\cdot \eta^{s+j,k-j}_2\wedge\omega^{s,j}_{\eta_1\wedge\eta_2}+
\end{eqnarray*}
\begin{eqnarray*}
+\Big(\bar{\partial}_{_J }\eta^{s,k}_1+\underbrace{ \sum_{j=1}^{k+1}  \omega^{s+j,k-j}\wedge \eta^{s,j}_1}_{(1)} -\underbrace{ \sum_{j=-1}^{k-1} (-1)^{k-j}\eta^{s+j,k-j}_1\wedge \omega^{s,j}_{\eta_1}}_{(2)}  +\omega^{s,k}\cdot g_{s,1}\Big)g_{s,2}\Big]  =(B_1) .\quad 
\end{eqnarray*}
En rappelant l'expression du terme $\bar{\partial}_{_J}(\eta_1\wedge \eta_2)^{s,k} $, en faisant le changement d'indice $j'=j+r$, $r'=j$ dans la première  et deuxième somme double et en regroupant opportunément les termes on obtient:
\begin{eqnarray*}
(B_1)= (g_{s+k,1}\cdot g_{s+k,2})^{-1}\Big[\bar{\partial}_{_J}(\eta_1\wedge \eta_2)^{s,k}
\underbrace{ +\sum_{j=0}^{k+1}  \omega^{s+j,k-j}\wedge(\eta_1\wedge \eta_2)^{s,j}}_{(1)} +\omega^{s,k} -
\\
\\ 
-\sum_{j=-1}^{k-1}(-1)^{k-j}\eta^{s+j,k-j}_1\wedge\Big(\bar{\partial}_{_J } \eta^{s,j}_2    
\underbrace{+\sum_{r=0}^{j+1}\omega^{s+r,j-r}_{\eta_1} \wedge\eta^{s,r}_2+\omega ^{s,j}_{\eta_1}}_{(2)} \Big)-\qquad\quad
\\
\\
-\sum_{j=-1}^{k-1}(-1)^{k-j}g_{s+k,1}\cdot\eta^{s+j,k-j}_2\wedge\omega^{s,j}_{\eta_1\wedge\eta_2}\,\Big]=(B_2) .\qquad\qquad\qquad\quad
\end{eqnarray*}
En utilisant l'hypothèse de récurrence et la définition des matrices $\omega^{s,j}_{\eta_1,\eta_2}$, on peut écrire le terme entre parenthèses rondes sous la forme suivante:
\begin{eqnarray*} 
\bar{\partial}_{_J } \eta^{s,j}_2  + \sum_{r=0}^{j+1}\omega^{s+r,j-r}_{\eta_1} \wedge\eta^{s,r}_2+\omega ^{s,j}_{\eta_1} =
\sum_{r=-1}^{j}(-1)^{j-r}\eta^{s+r,j-r}_2\wedge\omega^{s,r}_{\eta_1\wedge\eta_2}+\omega^{s,j}_{\eta_1\wedge\eta_2} .
\end{eqnarray*}
On aura alors
\begin{eqnarray*} 
&(B_2)=(g_{s+k,1}\cdot g_{s+k,2})^{-1}\Big[ \bar{\partial}_{_J}(\eta_1\wedge \eta_2)^{s,k}+
\displaystyle{ \sum_{j=0}^{k+1}\omega^{s+j,k-j}\wedge(\eta_1\wedge \eta_2)^{s,j} }
+\omega^{s,k} -&
\\
\\
&\displaystyle{ -\sum_{j=-1}^{k-1}\sum_{r=-1}^{j}(-1)^{k-r}\eta^{s+j,k-j}_1 }
\wedge\eta^{s+r,j-r}_2\wedge\omega^{s,r}_{\eta_1\wedge\eta_2}-&
\end{eqnarray*}
\begin{eqnarray*} 
-\sum_{j=-1}^{k-1}(-1)^{k-j}\Big(\eta^{s+j,k-j}_1+ \eta^{s+j,k-j}_2\Big)\wedge\omega^{s,j}_{\eta_1\wedge\eta_2}-
\sum_{j=-1}^{k-1}(-1)^{k-j}\eta^{s+k,0}\wedge\eta^{s+j,k-j}_2\wedge\omega^{s,j}_{\eta_1\wedge\eta_2}\,\Big] .
\end{eqnarray*}
En faisant le changement d'indice $j'=r$, $r'=j-r$ dans la somme double et en regroupant ce terme avec les deux dernières sommes, on obtient le terme cherché $\omega^{s,k}_{\eta_1\wedge\eta_2}$.  La transitivité de l'action $R$  est complètement claire par les  calculs qui ont permis de prouver que l'action même est bien définie\hfill $\Box$
\\
On va considérer maintenant quelques formules utiles pour le procédé itératif de la convergence rapide qui sera exposé en détail dans la sous-section $\ref{ss3.4}$. On explique formellement les  étapes  du procédé itératif. On désigne par $\omega_0:=\omega \in\Omega (U,\varphi,\psi,\bar{\partial})$ le choix initial de la calibration $\omega$. Au 
$k$-ième pas du procédé itératif on suppose avoir obtenu la calibration $\omega_k \in\Omega (U,[\varphi,\psi],\bar{\partial})$  et avoir déterminé le paramètre $\eta_{k+1}\in {\cal P}(U)$ en fonction de $\omega _k$. On définit alors la calibration $\omega _{k+1}:=R(\eta_{k+1},\omega _k)\equiv \omega _{k,\,\eta_{k+1} }$. Si on pose par définition
$$
\eta(k):={\displaystyle \bigwedge_{1\leq j\leq k }^{\longrightarrow}}\eta_j 
$$
où le symbole de produit extérieur avec une flèche vers la droite désigne le produit non commutatif de termes qui sont écrit en ordre croissant de l'indice vers la droite, 
on aura la formule  $\omega_k=\omega_{\eta(k)}$, grâce au fait que la recalibration $R$ est une action de semi-groupe. Si on pose par définition 
$g(k):=g(\eta(k))\in \Gamma(U)$ on aura que les composantes de $g(k)$ sont définies par la formule
$$
g_s(k) :={\displaystyle \prod_{0\leq j\leq k }^{\longrightarrow}}g_{s,j} ,
$$
pour $s=0,...,m$, où $g_j:=g(\eta_j)$ et $g_0:=\I_{p_0} $, (ici aussi le symbole de produit avec une flèche vers la droite désigne le produit non commutatif de termes qui sont écrits en ordre croissant de l'indice vers la droite). On écrit maintenant, à l'aide de cette dernière définition, les composantes 
$\eta(k)^{s,t},\;t\geq 1$ du paramètre $\eta(k)$ défini précédemment, sous une forme utile pour la convergence vers une solution d'un problème différentiel qu'on exposera dans la troisième partie.
\begin{lem} Pour tout entier $k \geq 1$ on peut écrire les composantes $\eta(k)^{s,t},\;t\geq 1$ du paramètre  $\eta(k)$ sous la forme
\begin{eqnarray}\label{5} 
\eta(k)^{s,t} =\Big(\sum_{\tau \in\Delta _t}\;\sum_{J\in J_k(\rho (\tau) )}\;{\displaystyle \bigwedge_{1\leq r\leq \rho (\tau) }^{\longrightarrow}} 
\;g_{s+\sigma' (\tau,r)}  (j_r-1)\cdot \eta^{s+\sigma (\tau,r),\,\tau_{\rho (\tau)+1-r } }_{j_r} \cdot g_{s+\sigma (\tau,r)} (j_r)^{-1}\,\Big) g_s(k)  
\end{eqnarray}
où
\begin{eqnarray*}
\Delta _t:=\left\{\tau\in\N^t\,|\,\tau_j\not=0 \Rightarrow \tau_{j-1} \not=0\,,\,\sum_{j=1}^t\tau_j =t\right\} ,
\\
\\  
\rho(\tau):=\max\{j\,|\,\tau_j\not=0\} ,\qquad\qquad\quad
\\
\\
J_k(\rho (\tau)):=\{J\in\{1,...,k\}^{\rho(\tau)}\,|\,j_1<...<j_{\rho(\tau)}\} ,
\\
\\
\sigma'(\tau,r):=\sum_{j=1}^{\rho (\tau)+1-r}\tau_j \quad\text{et}\quad\sigma(\tau,r):=\sum_{j=1}^{\rho (\tau)-r}\tau_j .
\end{eqnarray*}
\end{lem} 
Remarquons que $J_k(\rho (\tau))=\emptyset$ si $k<\rho (\tau)$.
\\
$Preuve$. On remarque que l'expression $\eqref{5}$ est évidente dans le cas $k=1,2$. Il est immédiat de vérifier à l'aide d'une récurrence élémentaire la validité de l'expression $\eqref{5}$ pour $t=1$ et $k\geq 1$ entier quelconque. Dans ce cas la formule $\eqref{5}$ s'écrit sous la forme
\begin{eqnarray*}
\eta(k)^{s,1}=\Big(\sum_{j=1}^kg_{s+1}(j-1)\cdot \eta^{s,1}_j \cdot g_s(j)^{-1}\Big)\cdot g_s(k)      
\end{eqnarray*}
On montre maintenant la validité de l'expression $\eqref{5}$ en général à l'aide du procédé récursif suivant. On suppose vraie la formule $\eqref{5}$ pour les composantes 
$\eta(k)^{\bullet,j} $, $j=1,...,t+1$ , $k \geq 1$ et on prouve la formule pour la composante $\eta(k+1)^{\bullet ,t+1}$ En effet on a 
\begin{eqnarray*}
\eta(k+1)^{s,t+1}:=(\eta(k)\wedge \eta_{k+1})^{s,t+1}:=\qquad\qquad\qquad
\\
\\
=\eta(k)^{s,t+1}\cdot g_{s,k+1}+g_{s+t+1} (k)\cdot \eta_{k+1}^{s,t+1}+      
\sum_{j=1}^t\eta(k)^{s+j,t+1-j}\wedge \eta_{k+1}^{s,j}=
\\
\\
=\Big(\sum_{\tau \in\Delta _{t+1} }\;\sum_{J\in J_k(\rho (\tau))}...\Big)\cdot g_s(k+1)  +g_{s+t+1} (k) \cdot \eta_{k+1}^{s,t+1}+\qquad      
\\
\\
\sum_{j=1}^t \sum_{\tau \in\Delta _{t+1-j} }\;\sum_{J\in J_k(\rho (\tau) )}\;{\displaystyle \bigwedge_{1\leq r\leq \rho (\tau) }^{\longrightarrow}} \;(...)\wedge
 g_{s+j} (k)\cdot  \eta_{k+1}^{s,j}=\qquad\quad
\end{eqnarray*}
\begin{eqnarray*}
= \Big(\sum_{j=1}^{k+1} g_{s+t+1}(j-1)\cdot \eta^{s,t+1}_j \cdot g_s(j)^{-1}\Big)\cdot g_s(k+1) +\qquad\qquad\qquad
\\
\\
+\Big({\displaystyle \sum_{\scriptstyle \tau \in\Delta _{t+1} 
\atop
\scriptstyle \rho (\tau) \geq 2}\;\sum_{J\in J_k(\rho (\tau))}...}\Big)\cdot g_s(k+1) + 
\Big({\displaystyle \sum_{\scriptstyle \tau \in\Delta _{t+1} 
\atop
\scriptstyle \rho (\tau) \geq 2}\;\sum_{\scriptstyle J\in J_{k+1} (\rho (\tau))
\atop
\scriptstyle  j_{\rho (\tau)}=k+1 }...}\Big)\cdot g_s(k+1) =\qquad
\\
\\ 
\Big(\sum_{j=1}^{k+1} g_{s+t+1}(j-1)\cdot \eta^{s,t+1}_j \cdot g_s(j)^{-1}\Big)\cdot g_s(k+1) + \Big({\displaystyle \sum_{\scriptstyle \tau \in\Delta _{t+1} 
\atop
\scriptstyle \rho (\tau) \geq 2}\;\sum_{J\in J_{k+1} (\rho (\tau))}...}\Big)\cdot g_s(k+1) 
\end{eqnarray*}
ce qui prouve la formule $\eqref{5}$ pour la composante $\eta(k+1)^{s,t+1}$.\hfill $\Box$

\subsection{Troisième étape: la formulation du problème différentiel}

La partie principale de la preuve consiste à prouver l'existence, pour tout $x\in U$, d'un voisinage ouvert $V \subset U$ de $x$ et $g\in\Gamma(V)$ solution du système différentiel
$$
(\Sigma )\;
\left  \{
\begin{array}{lr}
\bar{\partial}(\psi\cdot g_0)=0 
\\
\\
\bar{\partial}_{_J}(g_{s-1}^{-1} \cdot \varphi _s \cdot g_s)=0  
\\
\\
s=1,...,m .\qquad\qquad\qquad
\end{array}
\right.
$$
Bien évidemment résoudre ce système différentiel équivaut à trouver une autre ${\cal E}$-résolution de  ${\cal G}_{_{|_V}}$ dans la classe $[\varphi ,\psi]$, à partir de la  ${\cal E}$-résolution donnée $(\varphi ,\psi)$ de telle sorte qu'elle admet des matrices de connexion $\omega^{s,0}$  nulles. 
Maintenant on va prouver les deux résultats suivants.
\begin{lem}\label{3.3.1} 
Pour tout choix de calibration $\omega \in \Omega (U,\varphi ,\psi,\bar{\partial})$ on a que l'existence d'une solution $g\in \Gamma(U)$ du système différentiel $(\Sigma)$ est équivalente à l'existence d'une solution $\eta\in {\cal P}(U),\; g=g(\eta)$, du système différentiel quasi-linéaire
$$
(S_{\omega} )\;
\left  \{
\begin{array}{lr}
\displaystyle{\bar{\partial}_{_J }\eta^{s,k}+\sum_{j=0}^{k+1}  \omega^{s+j,k-j}\wedge \eta^{s,j}+(-1)^{k}\eta^{s-1,k+1}\wedge \omega^{s,-1}_{\eta} +\omega^{s,k}=0}
\\
\\
k=0,...,m 
\\
s=0,...,m-k .
\end{array}
\right.
$$ 
\end{lem}  
La proposition suivante permet de traduire notre problème en termes purement différentiels.
\begin{prop}\label{3.3} 
Supposons données des matrices $\omega ^{s,k}\in M_{p_{s+k},p_s }({\cal E}_X^{0,k+1}(U))$, $(s,k)\in I_m$ telles que $\omega ^{s-1,-1}\cdot\omega ^{s,-1}=0$ pour $s=2,...,m$ et $\bar{\partial}_{_J}\omega ^{s,-1}+\omega ^{s-1,0}\cdot \omega ^{s,-1}= \omega ^{s,-1}\cdot \omega ^{s,0}$ pour $s=1,...,m$. Alors, pour $k\geq 0$, les relations $(\eqref{3}_{s,k})$
\begin{eqnarray*}
\bar{\partial}_{_J }\omega^{s,k} +\sum_{j=-1}^{k+1} (-1)^{k-j}  \omega^{s+j,k-j}\wedge \omega^{s,j}=0 
\end{eqnarray*}
constituent les conditions d'intégrabilité du système différentiel $(S_{\omega} )$.
\end{prop} 
Venons-en maintenant à la preuve du lemme $\ref{3.3.1}$.
\\
\\
$Preuve\; du\; lemme\;\ref{3.3.1} $.
Soit $g\in \Gamma (U)$ une solution du système $(\Sigma)$. On écrit les composantes de $g$ sous la forme $g_s=\I_{p_s}+\eta^{s,0}$. Avec ces notations le système 
$(\Sigma)$ s'écrit sous la forme 
$$
(\Sigma_1)\;
\left  \{
\begin{array}{lr}
\bar{\partial}\psi_{\eta}=0 
\\
\\
\bar{\partial}_{_J }\omega ^{s,-1}_{\eta}=0  \qquad
\\
\\
s=1,...,m .
\end{array}
\right.
$$
On rappelle que les calculs utilisés pour montrer que l'application de recalibration $R$ est bien définie, (proppsition $\ref{3.2}$) nous donnent les égalités:
$$
\left  \{
\begin{array}{lr}
\bar{\partial}\psi_{\eta}=\psi \cdot ( \bar{\partial}_{_J }\eta^{0,0}+\omega ^{0,0}\wedge \eta^{0,0}+\omega ^{0,0})    
\\
\\
\bar{\partial}_{_J }\omega ^{s,-1}_{\eta}=-\omega ^{s-1,0}_{\eta}\cdot \omega ^{s,-1}_{\eta}+
\omega ^{s,-1}_{\eta}\cdot g_s^{-1} ( \bar{\partial}_{_J }\eta^{s,0}+\omega ^{s,0}\wedge \eta^{s,0}+\eta^{s-1,0}\wedge\omega ^{s,-1}_{\eta} +\omega ^{s,0})  
\\
\\
s=1,...,m .
\end{array}
\right.
$$
L'hypothèse d'exactitude de la ${\cal E}$-résolution locale implique alors l'existence d'une matrice $\eta^{0,1} \in M_{p_1,p_0}({\cal E}^{0,1}_X(U))$ telle que 
$$
\bar{\partial}_{_J }\eta^{0,0}+\sum_{j=0}^1 \omega ^{j,-j}\wedge \eta^{0,j}+\omega ^{0,0}=0 .
$$
L'équation précédente est bien évidemment équivalente à l'équation $\omega ^{0,0}_{\eta}=0$ (remarquons que la dépendance effective des matrices $\omega ^{\bullet,0}_{\eta}$ du paramètre $\eta$ est limitée aux composantes $\eta^{\bullet ,k},\;k=0,1$). On a donc que l'équation du système $(S_{\omega} )$, relative aux indices $(s,k)=(0,0)$ est satisfaite. On obtient alors à l'aide d'une récurrence croissante sur les indices $s=1,...,m$ relatifs aux expressions précédentes des matrices $\bar{\partial}_{_J }\omega ^{s,-1}_{\eta}$ et de l'exactitude de la  ${\cal E}$-résolution locale, l'existence de matrices $\eta^{s,1} \in M_{p_{s+1} ,p_s}({\cal E}^{0,1}_X(U)),\;s=0,...,m$ telles que $\omega ^{s,0}_{\eta}=0$ pour les indices en question. Ces équations ne représentent rien d'autre que le système différentiel quasi-linéaire
$$
(S_1)\;
\left  \{
\begin{array}{lr}
\displaystyle{\bar{\partial}_{_J }\eta^{s,0}+\sum_{j=0}^{1}  \omega^{s+j,-j}\wedge \eta^{s,j}+\eta^{s-1,1}\wedge \omega^{s,-1}_{\eta} +\omega^{s,0}=0}
\\
\\
s=0,...,m .
\end{array}
\right.
$$ 
qui est évidemment équivalent au système $(\Sigma)$. On rappelle aussi que les calculs relatifs à la bonne définition de l'application de recalibration $R$, nous donnent les égalités
\begin{eqnarray*}
\bar{\partial}_{_J }\omega^{s,k}_{\eta} +\sum_{j=-1}^{k} (-1)^{k-j}  \omega^{s+j,k-j}_{\eta}\wedge \omega^{s,j}_{\eta}=\qquad\qquad\qquad\qquad\qquad
\\
\\
=g^{-1}_{s+k} \omega^{s+k+1,-1}\Big(\bar{\partial}_{_J } \eta^{s,k+1} +\sum_{j=0}^{k+1}\omega^{s+j,k+1-j}\wedge\eta^{s,j} -\sum_{j=-1}^{k} (-1)^{k+1-j}\eta^{s+j,k+1-j}\wedge \omega^{s,j}_{\eta}+\omega ^{s,k+1} \Big)
\end{eqnarray*}
On obtient alors à l'aide d'une récurrence triangulaire, analogue à celle qui nous a permis de choisir les matrices $\omega ^{\bullet, \bullet}$ dans la première étape de la preuve, l'existence des matrices $\eta^{s,k}$ telles que $\omega^{s,k}_{\eta}=0$ pour $s=0,...,m,\;k=0,....,m$. Le système formé par ces équations est bien évidemment équivalent au système $(S_{\omega } )$, ce qui prouve le lemme.\hfill $\Box$
\\
Venons maintenant à la proposition $\ref{3.3}$. On commence par prouver la nécessité des conditions d'intégrabilité pour le système $(S_{\omega } )$. La suffisance de ces conditions sera prouvée dans l'étape suivante de la preuve du théorème $\ref{cdif}$.
\\
$Preuve\;de \;la \;n\acute{e}cessit\acute{e}\; des \; conditions \;d'int\acute{e}grabilit\acute{e}\;pour \;le\; syst\grave{e} me\;(S_{\omega} ).$
On commence par prouver la validité des relations $(\eqref{3} _{\bullet,k})$, $k\geq 0$ à l'aide d'une récurrence  croissante sur $k=0,..,m$. On rappelle que, le fait que par hypothèse on dispose des relations $\omega ^{s-1,-1}\cdot\omega ^{s,-1}=0,\;s=2,...,m$ et $(\eqref{3} _{\bullet,-1})$, combiné avec le fait que les équations du système $(S_{\omega } )$, relatives aux indices $(s,0)$ ne représentent rien d'autre que les équations $\omega^{s,0}_{\eta}=0$, implique la validité des équations 
$\bar{\partial}_{_J }\omega^{s,-1}_{\eta}=0$. On suppose par hypothèse de récurrence la validité des relations $(\eqref{3} _{\bullet,j})$ pour $j=-1,...,k-1$. En utilisant la validité des équations du système $(S_{\omega } )$ et en appliquant l'opérateur $\bar{\partial}_{_J }$ à l'équation relative aux indices $(s,k)$ de ce système on obtient les égalités suivantes
\begin{eqnarray*}
0=\bar{\partial}_{_J }\omega ^{s,k}\cdot g_s-(-1)^k \omega ^{s,k}\wedge \bar{\partial}_{_J }\eta ^{s,0} +
\sum_{j=1}^{k+1} \bar{\partial}_{_J } \omega ^{s+j,k-j}\wedge\eta^{s,j}-\quad
\\
\\
-\sum_{j=1}^{k+1} (-1)^{k-j} \omega ^{s+j,k-j}\wedge \bar{\partial}_{_J } \eta^{s,j} +(-1)^k\bar{\partial}_{_J }\eta^{s-1,k+1}\wedge \omega^{s,-1}_{\eta}=\quad\;
\\
\\
=\bar{\partial}_{_J }\omega ^{s,k}\cdot g_s-\sum_{j=1}^{k+1}\;\sum_{r=-1}^{k-j+1}(-1)^{k-j-r}\omega ^{s+j+r,k-j-r}\wedge \omega ^{s+j,r}\wedge \eta^{s,j}-
\\
\\
-\sum_{j=1}^{k+1} (-1)^{k-j} \omega ^{s+j,k-j}\wedge \bar{\partial}_{_J } \eta^{s,j} +(-1)^k\bar{\partial}_{_J }\eta^{s-1,k+1}\wedge \omega^{s,-1}_{\eta} .\qquad 
\end{eqnarray*}
En faisant le changement d'indice $j'=j+r$, $r'=j$ dans la somme double on obtient
\begin{eqnarray*}
0=\bar{\partial}_{_J }\omega ^{s,k}\cdot g_s-\sum_{j=0}^{k+1}(-1)^{k-j}
 \omega ^{s+j,k-j}\wedge\Big(\bar{\partial}_{_J } \eta^{s,j}+\sum_{r=1}^{j+1}\omega ^{s+r,j-r}\wedge \eta^{s,r}\,\Big)+
\\
\\
+(-1)^k\bar{\partial}_{_J }\eta^{s-1,k+1}\wedge \omega^{s,-1}_{\eta} =\Big(\bar{\partial}_{_J }\omega^{s,k} +\sum_{j=0}^{k+1} (-1)^{k-j}  \omega^{s+j,k-j}\wedge \omega^{s,j}\Big)\cdot g_s+
\\
\\
+(-1)^k\Big(\bar{\partial}_{_J }\eta^{s-1,k+1}+\sum_{j=0}^{k+1}\omega ^{s+j,k-j}\wedge \eta^{s-1,j+1}\Big)\wedge \omega^{s,-1}_{\eta}=\qquad\quad
\end{eqnarray*}
\begin{eqnarray*}
=\Big(\bar{\partial}_{_J }\omega^{s,k} +\sum_{j=-1}^{k+1} (-1)^{k-j}  \omega^{s+j,k-j}\wedge \omega^{s,j}\Big)\cdot g_s  .
\end{eqnarray*}
L'inversibilité de $g_s$ permet alors de conclure la  preuve de la nécessité des conditions d'intégrabilité $(\eqref{3} _{\bullet,k})$, $k\geq 0$ pour le système $(S_{\omega } )$. \hfill $\Box$
\\
La proposition $\ref{3.3}$ nous suggère de considérer les définitions suivantes. On définit l'ensemble non vide
\begin{eqnarray*} 
\Omega (U,p)\subset \bigoplus_{(s,k)\in I_m} M_{p_{s+k} ,p_s}({\cal E}^{0,k+1}_X(U)) ,
\end{eqnarray*} 
$p=(p_1,...,p_m)$, constitué par les éléments $\omega =(\omega ^{s,t})_{s,t}$ tels que $\omega^{s-1,-1}\cdot \omega ^{s,-1}=0$ et la relation $\eqref{3}$ soit satisfaite. Si on dispose de matrices $\omega^{s,-1}_0\in M_{p_{s-1},p_s}({\cal E}_X(U)) $, $s=1,...,m$ qui vérifient la relation écrite précédemment on peut définir l'ensemble
$$
\Omega (U,\omega ^{\bullet,-1}_0):=\{\omega \in \Omega (U,p)\,|\,\exists g\in \Gamma(U):\omega ^{\bullet,-1}=\omega ^{\bullet,-1}_{0,\,g}\}  
$$
dont les éléments seront encore appelés calibrations. Les calculs relatifs à la proposition $\ref{3.2}$ permettent d'étendre la recalibration $R$ à l'application
$$
R:{\cal P}(U)\times \Omega (U,\omega ^{\bullet,-1}_0)     \longrightarrow \Omega (U,\omega ^{\bullet,-1}_0)
$$
laquelle est encore une action de semi-groupe. A partir de maintenant on va considérer plus généralement la recalibration en termes de l'application définie précédemment.
Venons-en maintenant  à un  préliminaire technique avant d'exposer la preuve de l'existence des solutions pour le système différentiel $(S_{\omega})$.

\subsection{Détails sur le choix des normes et sur l'opérateur de Leray-Koppelman}

A partir de maintenant on va supposer que $U=B_1(0)$ est la boule ouverte de $\C^n$ de centre l'origine et de rayon unité. 
Si $A\in M_{k,l}(\C)$  on définit la norme $\|A\|:=\sup_{v\in\C^l-\{0 \} }\|Av\|/\|v\|$ et si $u=\sum'_{|I|=q}u_I\, d\bar{z}_I$ est une $(0,q)$-forme à coefficients des $(k,l)$-matrices à coefficients dérivables jusque à l'ordre $h\geq 0$, on définit une norme de Hölder invariante par changement d'échelle
$$
\|u\|_{r,\,h,\,\mu ,\,q}:=\sum _{\displaystyle
\scriptstyle   |I|=q 
\atop
\scriptstyle |\alpha |\leq h }
S_{|\alpha |} \,r^{|\alpha |+q}\,\|\partial^\alpha u_I\|_{r,\mu } 
$$    
où 
$$
\|f\|_{r,\,\mu}:=\sup_{z\in B_r} \,\|f(z)\|+\sup_{\displaystyle
\scriptstyle  z,\zeta \in B_r
\atop
\scriptstyle  z\neq \zeta }
 r^{\mu}\,\frac{ \|f(z) -f(\zeta ) \|}{\|z-\zeta \|^{\mu } } ,
$$    
avec $\mu \in (0,1)$ une constante fixée une fois pour toutes dans notre problème et $(S_k)_{k\geq 0}\subset (0,\infty),\;S_0:=1$ est une suite de réels qu'on construira ensuite et qui vérifie l'inégalité
$$
S_k\leq\left[\max_{|\alpha +\beta |=k}{\alpha +\beta\choose \alpha  }\right]^{-1} S_j\,S_{k-j}   
$$
pour tout $k\geq 1$ et $j=1,...,k-1$ On remarque que si le degré $q\geq 1$ on a que la norme $\|u\|_{r,h,q}$ tend vers zéro lorsque le rayon $r$ tend vers zéro. On désignera par ${\cal C}^{h,\mu} _{0,q}(\bar{B}_r, M_{k,l}(\C))$ l'espace de Banach de $(0,q)$-formes sur la boule fermée $\bar{B}_r $ à coefficients des $(k,l)$-matrices à coefficients dérivables jusque à l'ordre $h\geq 0$, telles que la norme $\|\cdot\|_{r,h,q}$ soit finie (on ne notera pas les dimensions des matrices).  On remarque que si $u\in {\cal C}^{h,\mu } _{0,q}(\bar{B}_r , M_{k,l}(\C))$ et  $v\in {\cal C}^{h,\mu } _{0,p}(\bar{B}_r, M_{l,t}(\C))$ alors on a l'inégalité $\|u \wedge v\|_{r,h ,p+q}\leq \|u\|_{r,h,q}\cdot\|v\|_{r,h ,p}$. On rappelle  très rapidement la définition de l'opérateur de Leray-Koppelman (voir les ouvrages classiques de Henkin-Leiterer \cite{He-Le}, de Range \cite{Ra} et l'article de Harvey-Polkin \cite{Ha-Po}). L'opérateur de Leray-Koppelman de la boule unité 
$$
T_q: {\cal C}^{h,\mu } _{0,q+1}(\bar{B}_1, M_{k,l}(\C))\longrightarrow  {\cal C}^{h,\mu } _{0,q}(B_1, M_{k,l}(\C))
$$
pour $q\geq 0$ est défini par une formule du type
$$
T_q\,u(z):=\int\limits_{\zeta \in B_1} u(\zeta )\wedge K_q (\zeta ,z)+ 
\int\limits_{\zeta \in \partial B_1} u(\zeta )\wedge k_q (\zeta ,z)
$$
où le premier opérateur intégral s'exprime en termes des coefficients $u_I$ de la forme $u$, par des termes du type
$$
\int\limits_{\zeta \in B_1} u_{I} (\zeta )\cdot K (\zeta ,z)\,d\lambda (\zeta )\qquad\text{avec}\qquad K (\zeta ,z)=\frac{\bar{\zeta }_l-\bar{z}_l  }{|\zeta -z|^{2n} }  
$$
et le deuxième par des termes du type
$$
\int\limits_{\zeta \in \partial B_1} u_{I} (\zeta )\cdot k (\zeta ,z)\,d\sigma  (\zeta )\qquad\text{avec}\qquad k(\zeta ,z)=
\frac{(\bar{\zeta }_j-\bar{z}_j)\cdot\bar{\zeta }_k  }{|\zeta -z|^{2l+2}\cdot\left[\bar{\zeta }\cdot(\bar{\zeta }-\bar{z})\right]^{n-1-l}  }  
$$
$ l=0,...,n-2$. L'opérateur de Leray-Koppelman $T_{r,q} : {\cal C}^{h,\mu} _{0,q+1}(\bar{B} _r, M_{k,l}(\C))\longrightarrow  {\cal C}^{h,\mu } _{0,q}(B_r, M_{k,l}(\C))$, $q\geq 0$ de la boule de rayon $r$ et de centre l'origine est défini par la formule $T^{r,q}:=(\lambda _r^{-1})^*\circ T_q \circ \lambda ^*_r$ avec $\lambda _r:B_1\rightarrow B_r$ l'homothétie de rapport $r$. Les propriétés de l'opérateur de Leray-Koppelman qui nous intéressent sont les suivantes:
\\
1) Pour toute forme différentielle $u\in {\cal C}^{h,\mu } _{0,q+1}(\bar{B} _r , M_{k,l}(\C))$ on a la formule d'homotopie:
\begin{eqnarray}
u= \bar{\partial}_{_J }T_{r,q}\, u+T_{r,q+1}\,\bar{\partial}_{_J }u
\end{eqnarray}
valable sur la boule $B_r$. 
\\
2) Il existe une suite de poids $S=(S_k)_{k\geq 0}$ de la norme de Hölder introduite précédemment telle que pour toute forme différentielle $u\in {\cal C}^{h,\mu } _{0,q+1}(\bar{B} _r , M_{k,l}(\C))$ on a l'estimation intérieure:
\begin{eqnarray}\label{7} 
\|T_{r,q} \,u\|_{r(1-\sigma ),\,h+1,\,\mu ,\,q} \leq C\cdot \sigma ^{-s(h)}\cdot\|u\|_{r,\,h,\,\mu ,\,q+1 }  , 
\end{eqnarray}
avec $\sigma \in (0,1)$, $s(h)=2n+k+2$ et $C=C(n,\mu )>0$ une constante $ind\acute{e} pendante$ de la régularité $h$. La preuve de la formule d'homotopie est exposée dans les ouvrages classiques mentionnés précédemment. Une estimation analogue à la $\eqref{7}$ a été déjà  montrée par S.Webster (voir \cite{We-1}). Nous utiliserons essentiellement les mêmes arguments de Webster pour montrer celle-ci. La différence avec l'estimation obtenue par Webster consiste dans le fait que la constante $C>0$ est indépendante de la régularité $h$. A partir de maintenant on désignera par $C$ une constante strictement positive indépendante de la régularité des formes. Pour prouver l'estimation $\eqref{7}$ il suffit de se restreindre au cas $r=1$, la norme étant choisie invariante pour changement d'échelle. On considère à ce propos une fonction $\rho \in  {\cal C}^{\infty}(\R,\,[0,1])$ telle que $\rho (x)=1$ pour $x\leq 0$ et $\rho (x)=0$ pour $x\geq 1$. On définit alors la fonction de cutt-off $\chi_{\sigma }$, avec $\sigma \in(0,1)$, par la loi 
$$
\chi _{\sigma } (z):=\left  \{
\begin{array}{lr}
1  &  \text{si}\quad  |z|\leq 1-\sigma/2 
\\
\rho (2\sigma ^{-1}(|z|-1+\sigma/2 )) & \text{si} \quad 1-\sigma/2\leq |z| 
\end{array}
\right.
$$ 
On aura alors, comme conséquence de l'invariance par translation de $ K(\zeta ,z)$, les égalités suivantes: 
\begin{eqnarray*}
J_1^{\alpha }(z):= \partial^{\alpha }_z \int\limits_{\zeta \in B_1} u_I (\zeta )\cdot K (\zeta ,z)\,d\lambda (\zeta )=\qquad\qquad\qquad\qquad\qquad
\\
\\
=\partial^{1_{\alpha } } _z \int\limits_{\zeta \in B_1} \partial^{\alpha-1_{\alpha }  }_{\zeta }(\chi_{\sigma }\cdot u_I) (\zeta )\cdot K (\zeta ,z)\,d\lambda (\zeta )+
 \int\limits_{\zeta \in B_1-B_{1-\sigma/2 } }((1-\chi_{\sigma })\cdot u_I) (\zeta )\cdot\partial^{\alpha }_z K (\zeta ,z)\,d\lambda (\zeta )
\end{eqnarray*}
pour $|z|\leq 1-\sigma$ et pour tout multi-indice $|\alpha |=k+1, \;k\geq 0$. Ici on désigne par $1_{\alpha}$ un multi-indice tel que $|1_{\alpha}|=1$ et $1_{\alpha }\leq \alpha $. La théorie classique du potentiel (voir par exemple \cite{Gi-Tru}) nous fournit alors les estimations 
\begin{eqnarray*}
\left\| \partial^{1_{\alpha } } _z \int\limits_{\zeta \in B_1} \partial^{\alpha-1_{\alpha }}_{\zeta }(\chi_{\sigma }\cdot u_I) (\zeta )\cdot K (\zeta ,z)\,d\lambda (\zeta )\right\|_{1-\sigma ,\,\mu } \leq C(n,\mu )\cdot\sigma^{-2n-1}\cdot\|\partial^{\alpha-1_{\alpha }}(\chi_{\sigma }\cdot u_I)\|_{1 ,\,\mu }  
\end{eqnarray*}
et
$
|\partial^{\alpha }_z K (\zeta ,z)|\leq  C(n,\mu )\cdot|\zeta -z|^{1-2n-|\alpha |}  
$
. On aura alors l'estimation  suivante:
\begin{eqnarray*}
\|J_1^{\alpha } \|_{1-\sigma ,\,\mu } \leq C(n,\mu)\cdot \sigma ^{-2n-1}  \cdot \|\partial^{\alpha-1_{\alpha}}(\chi_{\sigma }\cdot u_I)\|_{1,\,\mu }   
+C(n,|\alpha |)\cdot(\sigma/2)^{1-2n-|\alpha |} \cdot\|u_I \|_{1,\,0} 
\end{eqnarray*}
et donc
\begin{eqnarray*}
\|J_1\|_{1-\sigma ,\,h+1,\,\mu } \leq \sum_{|\alpha |\leq h+1}S_{|\alpha |}  \|J_1^{\alpha } \|_{1-\sigma, \,\mu }\leq \qquad\qquad\qquad\qquad
\\
\\
\leq\|J_1\|_{1-\sigma ,\,\mu } +  C(n,\mu )\cdot \sigma ^{-2n-1} \cdot\sum_{|\alpha |\leq h} \sum_{\beta \leq\alpha }S_{|\alpha |}{\alpha \choose\beta }  
\,\sigma ^{-|\beta |}\|\partial^{\beta  }\rho  \|_{1,\,\mu } \cdot\|\partial^{\alpha-\beta  }u_I \|_{1,\,\mu }+
\\
\\
+\sigma ^{-2n-h} \cdot\|u_I \|_{1,\,0}\cdot\sum_{|\alpha |\leq h+1} S_{|\alpha |} \cdot C(n,|\alpha |)\cdot 2^{|\alpha |+2n-1} .\qquad\qquad\qquad
\end{eqnarray*}
Pour un choix convenable de la suite $ S:=(S_k)_{k\geq 0}$, qu'on présentera ensuite, on peut se ramener à supposer que $\|\rho  \|_{1,\,S,\,\mu }:=\sum_{\alpha \geq 0}S_{|\alpha |}\|\partial^{\alpha  }\rho  \|_{1,\,\mu }<+\infty $ et $\sum_{\alpha \geq 0}S_{|\alpha |} \cdot C(n,|\alpha |)\cdot 2^{|\alpha |+2n-1} <+\infty  $. On aura alors l'estimation 
\begin{eqnarray}
\|J_1^{\alpha } \|_{1-\sigma ,\,h+1,\,\mu }
\leq C(n,\mu ) \cdot\sigma ^{-2n-1}\cdot\|u_I \|_{1,\,\mu } +C(n,\mu ) \cdot\sigma ^{-2n-h-1}\cdot \|\rho \|_{1,\,h,\,\mu } \cdot\|u_I \|_{1,\,h,\,\mu }+\nonumber
\\\nonumber
\\
+C\cdot \sigma ^{-2n-h} \cdot\|u_I \|_{1,\,0} \leq C \cdot\sigma ^{-2n-h-1}\cdot\|u_I \|_{1,\,h,\,\mu} .\qquad\qquad\qquad\qquad\label{8} 
\end{eqnarray}
Enfin pour estimer les termes du type
\begin{eqnarray*}
J_2^{\alpha }(z):= \partial^{\alpha }_z \int\limits_{\zeta \in\partial  B_1} u_I (\zeta )\cdot k (\zeta ,z)\,d\sigma  (\zeta )=\int\limits_{\zeta \in\partial  B_1} u_I (\zeta )\cdot\partial^{\alpha }_z k (\zeta ,z)\,d\sigma  (\zeta )
\end{eqnarray*}
 avec $|z|\leq 1-\sigma $, il suffit de dériver $|\alpha |+1$ fois le noyau $k(\zeta ,z)$, de remarquer l'estimation élémentaire:
\begin{eqnarray*}
\frac{|J_2^{\alpha }(z)-J_2^{\alpha }(\bar{z})| }{|z-\bar{z}|^{\mu } } \leq |z-\bar{z}|^{1-\mu }\cdot
\int\limits_{\zeta \in\partial  B_1} |u_I (\zeta )|\cdot|\nabla_z \,\partial^{\alpha }_z k (\zeta ,\hat{z}_{\zeta}  )|  \,d\sigma  (\zeta )
\end{eqnarray*}
(où $\hat{z}_{\zeta}$ est un point entre $z$ et $\bar{z} $) et les inégalités $|\zeta -z|\geq 3\sigma $, $ |\bar{\zeta }\cdot(\zeta -z)|\geq 3\sigma  $ pour $|\zeta |=1$ et $|z|\leq 1-3\sigma $. 
On aura alors l'estimation  
$$
\|J_2^{\alpha } \|_{1-\sigma ,\,\mu } \leq C(n,|\alpha |) \cdot \sigma^{-2n-1-|\alpha |}\cdot\|u_I \|_{1,0} .
$$
Par l'hypothèse faite précédemment sur la suite de poids on aura l'estimation 
$$
\|J_2^{\alpha } \|_{1-\sigma ,\,h+1,\,\mu } \leq C\cdot \sigma ^{-s(h)}\cdot\|u_I \|_{1,0}
$$
laquelle combinée avec l'estimation $\eqref{8}$ nous donne l'estimation $\eqref{7}$ sur la boule de rayon unité. 
\\
Venons-en à la définition partielle de la suite $S =(S_k)_{k\geq 0}$ laquelle sera déterminée en partie par l'exigence de satisfaire les hypothèses faites dans les calculs précédents. 
\\
\\
{\bf{Définition partielle de la suite de poids}} 
\\
On pose par définition
$$
A_k:=\sum_{|\alpha |=k}\max\{\|\partial^{\alpha}\rho \|_{1,\,\mu },\,\|\partial^{\alpha}\omega ^{s,-1}  \|_{1,\,\mu},\,s=0,..,m,\} , 
$$
$B_k:=(\max\{A_k,C(n,k)\})^{-1}$ si $\max\{A_k,C(n,k)\}\not=0$  et 1 sinon, $D_k:=[ \max_{|\alpha +\beta |=k}{\alpha +\beta \choose\alpha  }]^{-1}$. On pose par définition $S_0=1,\;S_1=B_1>0$ et on définit $S_k$, $k\geq 2$ à l'aide de la formule récursive
$$
0<S_k:=\min\{2^{-k}B_k,\,R_k,\,L_k,\,D_k\cdot\min_{1\leq j\leq k-1}S_j\cdot S_{k-j}\} 
$$ 
où $R_k$, $L_k$ sont des constantes qui seront déterminées dans la deuxième étape. On désigne par $S(\omega )$ la suite de poids obtenue si on pose 
$R_k=L_k=+\infty$, dans la définition précédente des poids.
Avec ces définitions on aura $\|\omega ^{\bullet,-1}\|_{1,S}\leq \|\omega ^{\bullet,-1}\|_{1,S(\omega )}<+\infty$. Pour simplifier les notations on identifiera dans la suite
$\|\cdot\|_{r,h,\mu ,q}\equiv \|\cdot\|_{r,h},\;T_{r,q}\equiv T_r $  et $\|\partial^h f\|_{\bullet}\equiv\sum_{|\alpha |=h}\|\partial^{\alpha }  f\|_{\bullet}$. 
\subsection{Quatrième étape: présentation du schéma de convergence rapide de type Nash-Moser et existence d'une solution du problème différentiel $(S_{\omega} )$}\label{s3.4}
\subsubsection{Preuve de l'estimation fondamentale du schéma de convergence rapide }\label{ss3.4.1}  

Dans cette partie de la preuve on va montrer l'existence d'un paramètre de recalibration $\eta$ des calibrations $\omega$, lequel permettra un contrôle quadratique de la norme des matrices $\omega ^{\bullet ,t}_{\eta},\;t\geq 0$ en termes de la norme des matrices $\omega ^{\bullet ,t},\;t\geq 0$. Ce contrôle est essentiel pour montrer la convergence vers zéro de la norme des matrices $\omega ^{\bullet ,t}_k,\;t\geq 0$ obtenues au $k$-ième pas du procédé itératif de la convergence rapide. La convergence vers une solution du problème différentiel $(S_{\omega} )$ est appelée rapide en raison de de l'estimation quadratique mentionnée précédemment. Avant de prouver l'estimation en question on va introduire quelques notations utiles pour la suite. Soit $\omega \in \Omega (B_1, p)$. Pour $r\in (0,1)$ on définit les quantités  
\begin{eqnarray*} 
a_h(\omega ,r):=\max\{\|\omega ^{s,k}\|_{r,h}\,|\,0\leq s \leq m\,,0\leq k \leq m-s\} ,
\\
\\
c(\omega ):=\max\{\|\omega ^{s,-1}\|_{1,\,S(\omega )}\,|\,1\leq s \leq m\} .\qquad\quad
\end{eqnarray*} 
On remarque que par définition de la norme de Hölder, la quantité $a_h(\omega ,r)$ tend vers zéro lorsque le rayon $r$ tend vers zéro. Pour tout $\sigma \in (0,1)$ on définit les rayons
$
r_l:=r(1-l\cdot \sigma _m)
$
pour $l=0,...,m+1$ où on pose par définition $\sigma _m:=\sigma /(m+1)$. Ensuite on définit par récurrence décroissante sur $k=m,...,0$, les constantes $L_k=L_k(C,c(\omega ))>0$ par les formules $L_m:=C$ et $L_{k-1}:=\max\{C, 2c(\omega )\cdot C\cdot L_k\}$. A partir de maintenant on désigne par $\varepsilon \in (0,1/2)$ une constante fixée telle que pour toutes les matrices $A\in M_{p_s,p_s}(\C)$ telles que $\|A\|<\varepsilon$ on a l'inversibilité de la matrice $\I_{p_s}+A $. Avec ces notations on a la proposition suivante:
\begin{prop}\label{3.4} 
Supposons donnés $\omega \in \Omega (B_1, p),\;r,\,\sigma \in (0,1),\;h\in \N$ et les poids $0<S_j\leq S_j(\omega ),\;j=0,...,h+1$ de la norme de Hölder 
$\|\cdot\|_{r,h+1}$. Supposons que le rayon $r$ soit suffisamment petit pour assurer les estimations 
$$
L_0(C,c(\omega))\cdot \sigma _m^{-(m+1)\cdot s(h)}\cdot a_h(\omega ,r)<\varepsilon  
$$
et $a_h(\omega ,r)\leq 1$, où la quantité $a_h(\omega ,r)$ est calculée par rapport au poids $S_j,\;j=0,...,h$. Supposons de plus que  le poids  $S_{h+1}$ soit suffisamment petit pour pouvoir assurer l'estimation:
\begin{eqnarray*}
S_{h+1} \|\partial^{h+1} \omega^{s,k}_I \|_{r ,\mu } \leq \| \omega^{s,k}_I \|_{r ,\mu }
\end{eqnarray*}
pour tout $k=0,...,m,\;s=0,...,m-k$ et $|I|=k+1$. Si on définit les composantes  $\eta^{s,k}$ du paramètre de recalibration 
 $\eta\in {\cal P}(\bar{B} _{r(1-\sigma )})$ par la formule de récurrence décroissante sur 
$k=m,...,0$
$$
\eta^{s,k}:=-T_{r_{m-k}} \,\Big(\omega ^{s,k}+\omega ^{s+k+1,-1}\wedge\eta^{s,k+1}+(-1)^k\eta^{s-1,k+1} \wedge\omega ^{s,-1}\Big)\in M_{p_{s+k},p_s}({\cal E}^{0,k}_X(B_{r_{m-k}} )  )     
$$ 
pour $s=0,...,m-k$, alors on aura la validité des  estimations
\begin{eqnarray}
&\|\eta^{\bullet,k}  \|_{r(1-\sigma ) ,h+1 } \leq L\cdot\sigma_m^{-(m+1)\cdot s(h)}\cdot a_h(\omega ,r) ,&  \label{9} 
\\\nonumber
\\
&\|\omega^{\bullet,k}_{\eta} \|_{r(1-\sigma ),h+1} \leq R \cdot\sigma_m ^{\nu (m,h)}\cdot a_h(\omega ,r) ^2 &\label{10}
\end{eqnarray}
 pour tout $k= 0,...,m$, avec $L=L(C,c(\omega )):=\max_k L_k\,>0,\;R=R(C,c(\omega ))>0$ constantes positives et $\nu (m,h):=[(m+2)\cdot m+1]\cdot s(h)$, $(m\geq 0,\,h\geq 0)$. 
\end{prop} 
$Preuve$.
Dans les calculs qui suivront on utilisera les identifications $a_h\equiv a_h(\omega ,r)$ et $c=c(\omega) $. On commence par prouver l'estimation $\eqref{9}$. Si on définit $\sigma_{m,l}>0$ par la formule $r_{l+1}=r_l(1+\sigma_{m,l})$ on a que $\sigma_{m,l}\geq \sigma _m$. On obtient alors à l'aide de l'estimation $\eqref{7}$ et  d'une récurrence élémentaire décroissante sur $k=m,...,0$, l'estimation suivante
\begin{eqnarray}
\|\eta^{\bullet,k}  \|_{r_{m-k+1} ,h+1 } \leq L_k\cdot\sigma_m^{-(m-k+1)\cdot s(h)}\cdot a_h   ,\label{11}
\end{eqnarray}
laquelle prouve l'estimation $\eqref{9}$. Venons maintenant à la preuve de l'estimation $\eqref{10}$. Pour cela on définit les troncatures 
$\eta_{[t]}  :=({\eta_{[t]}}^{s,k} )_{s,k}\in{\cal P}(B_{r_{m-t}})$ du paramètre $\eta$ défini dans l'hypothèse de la proposition $\ref{3.4}$, de la façon suivante;  
${\eta_{[t]}}^{s,k} :=0$ si $k<t$ et  
${\eta_{[t]}}^{s,k}  :=\eta^{s,k}$ sur la boule $B_{r_{m-t}}$, si $k\geq t$. Par définition de la recalibration avec paramètre $\eta_{[k+1]}$ on aura alors: 
\begin{eqnarray*}
\omega ^{s,k}_{\eta_{[k+1]}}=  \omega^{s,k}+\omega^{s+k+1,-1} \wedge \eta^{s,k+1}+(-1)^k\eta^{s-1,k+1} \wedge\omega ^{s,-1}
\end{eqnarray*}
sur la boule $B_{r_{m-k-1}} $ pour $k=0,...,m$. 
On montre maintenant à l'aide d'une récurrence en ordre décroissant sur $k=m,...,0$, l'estimation quadratique 
\begin{eqnarray}
\|T_{r_{m-k}}\,\bar{\partial}_{_J }\omega^{s,k}_{\eta_{[k+1]} }\|_{r_{m-k+1} ,h+1}\leq Q_k\cdot\sigma^{-b(m,k,h)}_m \cdot a_h^2 \label{12} 
\end{eqnarray}
 où $Q_k=Q_k(C,c)>0$ est une constante positive, $b(m,k,h):= [2(m-k)+1]\cdot s(h)$. Pour $k=m$ on a banalement $\omega^{0,m}_{\eta_{[m+1]} } =\omega^{0,m}$. L'estimation précédente découle alors immédiatement de la relation $(\eqref{3} _{0,m})$, laquelle peut être écrite explicitement sous la forme
\begin{eqnarray*}
\bar{\partial}_{_J }\omega^{0,m}
= -\sum_{j=0}^{m} (-1)^{m-j}\,\omega^{j,m-j}\wedge \omega^{0,j}
\end{eqnarray*}
et de l'estimation $\eqref{7}$. Le diagramme de figure $\ref{fig6}$ montre les matrices qui interviennent dans la relation $\eqref{3}$ dans le cas où $k=m$.
\begin{figure}[hbtp]
\begin{center} 
\input dgt6.pstex_t
    \caption{}
    \label{fig6}
\end{center} 
\end{figure}    
On remarque que la relation $(\eqref{3} _{0,m})$ est la seule, parmi les autres relations $(\eqref{3} _{\bullet,\bullet})$, qui ne présente pas de facteurs de type $\omega^{\bullet,-1} $ dans les termes quadratiques. 
\\ 
Montrons maintenant l'estimation quadratique $\eqref{12}$ pour $0\leq k-1<m$ (si $m\geq 1$, autrement il n'y a plus rien à prouver en ce qui concerne l'estimation $\eqref{12}$) en admettant qu'elle soit vraie pour $1\leq k\leq m$. En effet en  utilisant la relation $\eqref{3}$ relativement aux matrices $\bar{\partial}_{_J } \omega ^{s,k-1} _{\eta_{[k]} } $ on obtient l'expression suivante:
\begin{eqnarray}
\bar{\partial}_{_J }\omega^{s,k-1}_{\eta_{[k]} } =(-1)^{k+1}  \omega^{s-1,k}_{\eta_{[k]} }\wedge \omega^{s,-1}+ 
 \omega^{s+k,-1}\wedge \omega^{s,k}_{\eta_{[k]} }+\nonumber
\\\nonumber
\\
+\sum_{j=0}^{k-1} (-1)^{k-j} \omega^{s+j,k-1-j}_{\eta_{[k]} }\wedge \omega^{s,j}_{\eta_{[k]} }\qquad\qquad\quad\label{13} 
\end{eqnarray}
(remarquons que $\omega^{\bullet,-1}_{\eta_{[k]} }=\omega^{\bullet,-1}$ étant $k\geq 1$). En explicitant à l'aide de la formule de recalibration relative au paramètre $\eta_{[k]}$ les termes $\omega^{\bullet,k}_{\eta_{[k]}}$ qui apparaissent dans l'expression précédente et en utilisant la formule d'homotopie pour l'opérateur $\bar{\partial}_{_J }$ on obtient:
\begin{eqnarray*}
&\omega^{s,k}_{\eta_{[k]} }= \bar{\partial}_{_J }\eta^{s,k}+  \omega^{s+k,0}\wedge \eta^{s,k}-(-1)^{k}\eta^{s,k}\wedge \omega^{s,0}_{\eta_{[k]} } +\omega^{s,k}_{\eta_{[k+1]} }=&
\\
\\
&=T_{r_{m-k}}\,\bar{\partial}_{_J }\omega^{s,k}_{\eta_{[k+1]} }+\omega^{s+k,0}\wedge \eta^{s,k}-(-1)^{k}\eta^{s,k}\wedge \omega^{s,0}_{\eta_{[k]} } .&        
\end{eqnarray*}
Le diagramme suivant montre les matrices qui interviennent dans la définition de la matrice $\omega^{1,2}_{\eta_{[2]}}$ dans le cas  $m=4$.
\begin{figure}[hbtp]
\begin{center} 
\input dgt5.pstex_t
    \caption{}
    \label{fig5}
\end{center} 
\end{figure}    
\\
On estime donc la norme des matrices $\omega^{s,k}_{\eta_{[k]} }$ à l'aide de l'expression obtenue précédemment, de l'hypothèse récursive et de l'estimation $\eqref{11}$.
\begin{eqnarray*}
& \|\omega^{s,k}_{\eta_{[k]} }\|_{r_{m-k+1},h }  \leq Q_k \cdot \sigma_m^{-[2(m-k)+1]\cdot s(h)}\cdot a^2_h
+\|\omega^{s+k,0}\wedge \eta^{s,k} \|_{r_{m-k+1},h } +&
\\
\\
&+\|\eta^{s,k}\wedge \omega^{s,0} \|_{r_{m-k+1},h}
+2c\cdot\|\eta^{s,k}  \|_{r_{m-k+1},h}\cdot \|\eta^{\bullet,k } \|_{r_{m-k+1},h} \leq &
\end{eqnarray*}
\begin{eqnarray*}
&\leq Q_k \cdot \sigma_m^{-[2(m-k)+1]\cdot s(h)}\cdot a^2_h+&
\\
\\
&+2 L_k\cdot \sigma^{-(m-k+1)\cdot s(h)}_m \cdot a^2_h +2c\cdot L_k^2 \cdot \sigma^{-2(m-k+1)\cdot s(h)}_m \cdot a^2_h . &
\end{eqnarray*}
On a alors prouvé l'estimation quadratique
\begin{eqnarray*}
\|\omega^{s,k}_{\eta_{[k]} }\|_{r_{m-k+1},h }  \leq H_k\cdot \sigma^{-2(m-k+1)\cdot s(h)}_m \cdot a^2_h  
\end{eqnarray*}
$(k\geq 1)$, où $H_k =H_k (C,c)>0$ est une constante positive (on remarque  que le terme quadratique 
$2c\cdot\|\eta^{s,k}  \|_{r_{m-k+1},h}\cdot \|\eta^{\bullet,k } \|_{r_{m-k+1},h}  $, et donc le dernier terme 
$$
2c\cdot L_k^2 \cdot \sigma^{-2(m-k+1)\cdot s(h)}_m \cdot a^2_h  
$$
de la dernière inégalité précédente, sont présents seulement si $k=1$ du fait que $\omega^{\bullet,0}_{\eta_{[t]}}=\omega^{\bullet,0}$ si $t\geq 2$). On estime maintenant la norme 
$\|T_{r_{m-k+1}}\,\bar{\partial}_{_J }\omega^{\bullet,k-1}_{\eta_{[k]} } \|_{r_{m-k+2},h+1 }$ à l'aide de l'expression $\eqref{13}$ et de l'estimation obtenue précédemment. On a:
\begin{eqnarray*}
&\|T_{r_{m-k+1}}\,\bar{\partial}_{_J }\omega^{\bullet,k-1}_{\eta_{[k]} } \|_{r_{m-k+2},h+1 }  
\leq 2c\cdot C\cdot \sigma_m^{-s(h)}\cdot\|\omega^{\bullet,k}_{\eta_{[k]} }\|_{r_{m-k+1},h } + &
\\
\\
&+k\cdot C\cdot \sigma_m^{-s(h)}\Big(a_h +2c\cdot \|\eta^{\bullet,k}\|_{r_{m-k+1},h } \Big)^2\leq&
\\
\\
&\leq 
2c\cdot C\cdot H_k\cdot  \sigma_m^{-[2(m-k+1)+1]\cdot s(h)}\cdot a^2_h+&
\\
\\
&+k\cdot C\cdot \sigma_m^{-s(h)}\Big(a_h +2c\cdot L_k\cdot  \sigma^{-(m-k+1)\cdot s(h)}_m\cdot a_h\Big)^2  &
\end{eqnarray*}
ce qui prouve l'estimation $\eqref{12}$ pour l'indice $k-1$ et donc pour tous les indices $k=0,...,m$. [ici aussi on remarque que les termes quadratiques 
$(2c)^2\cdot \|\eta^{\bullet,k}\|^2_{r_{m-k+1},h }$, et donc les termes qui en dérivent
$$
\Big(2c\cdot L_k\cdot  \sigma^{-(m-k+1)\cdot s(h)}_m\cdot a_h\Big)^2
$$ 
 et qui apparaissent dans l'estimation précédente, sont présents seulement dans le cas $k=1$. En effet, pour arriver à cette conclusion il suffit de tenir compte de la relation $\omega^{\bullet,j}_{\eta_{[t]}}=\omega^{\bullet,j}$ si $j<t-1$ dans l'expression $\eqref{13}$. Le cas $k=1$ est celui pour lequel le ``poids'' $\sigma _m^{-b(m,k-1,h)}$ figurant dans l'estimation de la norme 
$\|T_{r_{m-k+1}}\,\bar{\partial}_{_J }\omega^{\bullet,k-1}_{\eta_{[k]} } \|_{r_{m-k+2},h+1 }$ est le plus grand]. On est maintenant en position de prouver l'estimation $\eqref{10}$.
On pose par définition
$\theta^{s,0}:= (\I_{p_s}+\eta^{s,0})^{-1}-\I_{p_s}\in M_{p_s,p_s}({\cal E}_X(B_{r_m})) $. 
Le fait que $\sup_{z\in B_{r(1-\sigma )}}\|\eta^{s,0}(z)\| <\varepsilon <1$ implique l'égalité
$$
\theta ^{s,0} ={\displaystyle \sum_{j=1}^{\infty}\,(-1)^j(\eta^{s,0})^j}  
$$
(a priori la série précédente est convergente en topologie ${\cal C}^{0}$ vers l'élément $\theta ^{s,0}$ de classe ${\cal C}^{\infty}$, mais une étude élémentaire plus précise, dont on n'aura pas besoin ici, montre que l'estimation précédente $\sup_{z\in B_{r(1-\sigma )}}\|\eta^{s,0}(z)\| <\varepsilon <1$ est suffisante pour assurer la convergence de la série en topologie ${\cal C}^h$ pour tout $h$). L'inégalité
$$
\|\eta^{s,0}\|_{r(1-\sigma ),h+1}\leq L_0\cdot \sigma _m^{-(m+1)\cdot s(h)}\cdot a_h<\varepsilon <1/2
$$ 
implique la convergence de la série précédente en topologie ${\cal C}^{h+1,\mu }(\bar{B}_{r(1-\sigma )})$ et permet d'effectuer les estimations 
$$
\|\theta^{s,0}\|_{r(1-\sigma ),h+1}\leq \sum_{j=1}^{\infty}\|\eta^{s,0}\|_{r(1-\sigma ),h+1}^j\leq 2\|\eta^{s,0}\|_{r(1-\sigma ),h+1}
$$
et $\|g^{\pm 1}_{\bullet}  \|_{r(1-\sigma ),h+1}<2 $. On utilisera ces deux inégalités et l'estimation $\eqref{12}$ obtenue précédemment pour prouver  l'estimation $\eqref{10}$ sous la forme plus précise suivante:
\begin{eqnarray}
\|\omega^{\bullet,k}_{\eta} \|_{r(1-\sigma ),h+1} \leq R'_k \cdot\sigma_m ^{-[(k+2)\cdot m+1]\cdot s(h)}\cdot a_h^2 \label{14} 
\end{eqnarray}
où $R'_k=R'_k(C,c)>0$ est une constante positive. Nous considérons pourtant l'expression suivante de $\omega^{s,k}_{\eta}$, pour tous les indices $k=0,...,m$
\begin{eqnarray*}
&\omega^{s,k}_{\eta}= g_{s+k}^{-1}\cdot \Big(  \bar{\partial}_{_J }\eta^{s,k}+\displaystyle{\sum_{j=0}^{k}  \omega^{s+j,k-j}\wedge \eta^{s,j}-\sum_{j=0}^{k-1} (-1)^{k-j}\eta^{s+j,k-j}\wedge \omega^{s,j}_{\eta} +}&
\\
\\
&+(-1)^k\eta^{s-1,k+1}\wedge \theta ^{s-1,0}\wedge \omega ^{s,-1}+(-1)^k\eta^{s-1,k+1}\wedge g_{s-1}^{-1}\cdot\omega ^{s,-1}\wedge \eta^{s,0}+\omega^{s,k}_{\eta_{[k+1]} }\,\Big)=&
\end{eqnarray*}
\begin{eqnarray*}
&= g_{s+k}^{-1}\cdot \Big( T_{r_{m-k}}\,  \bar{\partial}_{_J }\omega^{s,k}_{\eta_{[k+1]} }\displaystyle{ +\sum_{j=0}^{k}  \omega^{s+j,k-j}\wedge \eta^{s,j}-\sum_{j=0}^{k-1} (-1)^{k-j}\eta^{s+j,k-j}\wedge \omega^{s,j}_{\eta} +}&
\\
\\
&+(-1)^k\eta^{s-1,k+1}\wedge \theta ^{s-1,0}\wedge \omega ^{s,-1}+(-1)^k\eta^{s-1,k+1}\wedge g_{s-1}^{-1}\cdot\omega ^{s,-1}\wedge \eta^{s,0}\Big) .&
\end{eqnarray*}
Dans le cas $k=0$, l'expression précédente s'écrit sous la forme:
\begin{eqnarray*}
\omega ^{s,0}_{\eta} =g^{-1}_s\cdot\Big(T_{r_m} \,\bar{\partial}_{_J }\omega ^{s,0}_{\eta_{[1]} }   +\omega ^{s,0}\wedge \eta^{s,0}+
\eta^{s-1,1}\wedge \theta ^{s-1,0}\wedge \omega ^{s,-1}+ \eta^{s-1,1}\wedge g^{-1}_{s-1}\cdot \omega ^{s,-1}\wedge \eta^{s,0} \Big) .
\end{eqnarray*}
On estime maintenant la norme $\|\omega^{s,0}_{\eta }\|_{r(1-\sigma ),h+1}$. L'hypothèse faite sur le poids $S_{h+1}$ implique l'inégalité
 \begin{eqnarray*}
\|\omega^{s,0}\wedge \eta^{s,0} \|_{r(1-\sigma ),h+1} \leq 
2\|\omega ^{s,0}\|_{r,h}\cdot\|\eta^{s,0}\|_{r(1-\sigma ),h+1}
\end{eqnarray*}
laquelle combinée avec les inégalités 
$\|\theta^{s,0}\|_{r(1-\sigma ),h+1}\leq 2\|\eta^{s,0}\|_{r(1-\sigma ),h+1}$, $\|g^{\pm 1}_{\bullet}  \|_{r(1-\sigma ),h+1}<2$, l'estimation $\eqref{11}$ et l'estimation $\eqref{12}$ pour $k=0$, prouvée précédemment, donne les estimations suivantes:
\begin{eqnarray*} 
&\|\omega^{s,0}_{\eta}\|_{r(1-\sigma ),h+1} \leq 2Q_0\cdot \sigma_m ^{-(2m+1)\cdot s(h)}\cdot  a^2_h
+4\|\omega ^{s,0}\|_{r,h}\cdot\|\eta^{s,0}\|_{r(1-\sigma ),h+1}+&
\\
\\
&+2c\cdot\|\eta^{s-1,1}\|_{r(1-\sigma ),h+1}\cdot \|\theta ^{s-1,0}\|_{r(1-\sigma ),h+1}+&
\\
\\
&+4c\cdot \|\eta^{s-1,1}\|_{r(1-\sigma ),h+1}\cdot\|\eta^{s,0}\|_{r(1-\sigma ),h+1}\leq&
\\
\\
&\leq 2Q_0\cdot \sigma_m ^{-(2m+1)\cdot s(h)}\cdot  a^2_h+4L_0\cdot \sigma_m ^{-(m+1)\cdot s(h)}\cdot  a^2_h+&
\\
\\
&+4c\cdot L_1\cdot L_0\cdot \sigma_m^{-(2m+1)\cdot s(h)}\cdot  a^2_h .&
\end{eqnarray*} 
ce qui prouve l'estimation $\eqref{14}$ dans le cas $k=0$. On montre  maintenant l'estimation quadratique $\eqref{14}$ pour tous les indices à l'aide d'une récurrence croissante sur $k=0,...,m-s$, appliquée à l'expression précédente de la matrice $\omega^{s,k}_{\eta}$  et à l'aide de l'estimation quadratique $\eqref{12}$. On suppose vraie l'estimation $\eqref{14}$ pour tout $j=0,...,k-1$ et on considère l'estimation suivante dans laquelle on utilise comme précédemment l'hypothèse faite sur le poids $S_{h+1}$ de la norme de Hölder et l'hypothèse $a_h\leq 1$:
\begin{eqnarray*}
&\|\omega^{s,k}_{\eta}\|_{r(1-\sigma) ,h+1} \leq 2Q_k\cdot \sigma_m ^{-[2(m-k)+1]\cdot s(h)}\cdot a^2_h +&
\\
\\
&+\displaystyle{4\sum_{j=0}^k  \|\omega^{s+j,k-j}\|_{r,h}\cdot \| \eta^{s,j}\|_{r(1-\sigma ),h+1}+ }&
\\
\\
&+\displaystyle{\sum_{j=0}^{k-1}\| \eta^{s+j,k-j} \|_{r(1-\sigma ),h+1}\cdot\|  \omega^{s,j}_{\eta}\|_{r(1-\sigma ),h+1}+}&
\\
\\
&+8c\cdot\| \eta^{s-1,k+1}\|_{r(1-\sigma ),h+1} \cdot\| \eta^{\bullet,0}\|_{r(1-\sigma ) ,h+1}\leq &
\\
\\
&\leq (2 Q_k+4(k+1)\cdot L)\cdot \sigma_m ^{-(2m+1)\cdot s(h)}\cdot a^2_h +&
\end{eqnarray*}
\begin{eqnarray*}
&+\displaystyle{\sum_{j=0}^{k-1}(L_{k-j}\cdot R'_j) \cdot \sigma_m ^{-[m+(j+2)\cdot m+1]\cdot s(h)}\cdot a^2_h+}&
\\
\\
&+8c\cdot L_{k+1}\cdot L_0 \cdot \sigma _m^{-(2m+1)\cdot s(h)}\cdot a^2_h .&
\end{eqnarray*}
La dernière inégalité implique évidement l'estimation $\eqref{14}$ et donc l'estimation $\eqref{10}$, ce qui prouve la proposition $\ref{3.4}$. \hfill $\Box$

\subsubsection{Bon fonctionnement du procédé itératif}\label{ss3.4} 

Dans cette partie on va établir les hypothèses qui permettent d'appliquer l'estimation fondamentale $\eqref{10}$ à une étape quelconque du procédé itératif. On commence par préciser les paramètres qui contrôleront les étapes de la convergence rapide. Pour cela on définit d'abord les paramètres $\sigma _k:=e^{-k-2}$ qui contrôleront les restrictions des rayons des boules, lesquels sont définis de façon récursive par la formule $r_{k+1}:=r_k(1-\sigma _k)$ pour tout entier $k \geq 0$. Bien évidemment le rayon limite 
$$
r_{\infty}:=\lim_{k\rightarrow +\infty}r_k =r_0{\displaystyle \prod_{k=0}^{\infty}(1-\sigma _k)}
$$
est non nul étant  $-\sum_{k=0}^{\infty}\log(1-\sigma _k)\leq Cst \sum_{k=0}^{\infty}\sigma _k <\infty$. 
Ensuite on désigne par $r(k,l):=r_k(1-l \cdot \sigma _{m,k}),\;l=0,...,m+1$ où $\sigma _{m,k}:=\sigma _k/(m+1)$. Pour le choix du rayon initial on considère les suites numériques
\begin{eqnarray*} 
\alpha _k(r):= a_0(r)^{2^k}\cdot \prod_{j=0}^{k-1} H^{2^{k-1-j}}\cdot\sigma _j^{-\nu (m,j)\cdot2^{k-1-j}} ,
\\
\\   
\beta  _k(r):= b_0(r)^{2^k}\cdot \prod_{j=0}^{k-1} P^{2^{k-1-j}}\cdot e^{-\gamma (m,j)\cdot 2^{k-1-j}} 
\end{eqnarray*}
où $b_0(r):=H\cdot \sigma _{m,0}^{-(m+1)\cdot s(0) }\cdot a_0(r)$, $H>0,\;P>0$ sont deux constantes (dépendantes seulement de $m$) et $\gamma (m,j) $ est une fonction affine strictement croissante en $j$.  On rappelle que par définition de la norme de Hölder on a que la quantité $a_0=a_0(r)$ tend vers zéro lorsque $r>0$ tend vers zéro. Avant de faire le choix du rayon initial on a besoin du lemme élémentaire suivant.
\addtocounter{subsubsection}{-1}  
\begin{lem} \label{3.4.2} 
Il existe $\rho \in (0,1)$ tel que;
\\
$(A)$, pour tout $r\in (0,\rho ]$  les séries numériques $\sum_{k\geq 0}\alpha _k(r) $ et $\sum_{k\geq 0}\beta_k(r)$  sont convergentes.
\\
$(B)$
\begin{eqnarray*}
\lim_{r\rightarrow 0^+ } \sum_{k=0}^{\infty}\,\alpha _k(r)=0, \quad\lim_{r\rightarrow 0^+ } \sum_{k=0}^{\infty}\,\beta _k(r)=0 
\end{eqnarray*}
\\
$(C)$, pour tout $r\in (0,\rho ]$ et entier $k\geq 0$ on a les inégalités $\alpha _k(r)\leq 1$ et $\beta _k(r)<\varepsilon$
\end{lem}
$Preuve$. Par le critère du rapport il suffit de vérifier que les suites $\ln(\alpha _{k+1}/ \alpha _k)$ et $\ln(\beta  _{k+1}/\beta_k)$ tendent vers $-\infty$ lorsque $k$ tend vers $+\infty$. En explicitant les logarithmes on a:
\begin{eqnarray*} 
\ln(\alpha _{k+1}/ \alpha _k )=2^k\cdot\Big( \ln a_0(\rho)+2^{-1}(\ln H)\sum_{j=0}^{k-1}2^{-j}-2^{-1}\sum_{j=0}^{k-1}\nu (m,j)2^{-j}\ln \sigma _{m,j}  \Big)-\qquad
\\
\\
-\nu (m,k)\ln\sigma _{m,k} +\ln H\qquad
\end{eqnarray*}
\begin{eqnarray*}
\ln(\beta  _{k+1}/\beta  _k) =2^k\cdot\Big( \ln b_0(\rho )+2^{-1}(\ln P)\sum_{j=0}^{k-1}2^{-j} + 2^{-1}\sum_{j=0}^{k-1}\gamma (m,j)2^{-j}\Big)+\gamma (m,k)+\ln P .
\end{eqnarray*}
Il suffit donc choisir $\rho >0$ suffisamment petit pour assurer les inégalités
\begin{eqnarray*} 
 \ln a_0(\rho )+2^{-1}(\ln H)\sum_{j=0}^{\infty}2^{-j}+2^{-1}\sum_{j=0}^{\infty}\nu (m,j)[j+2+\ln (m+1)]2^{-j}<0 ,
\\
\\
\ln b_0(\rho )+2^{-1}(\ln P)\sum_{j=0}^{\infty}2^{-j} +  2^{-1}\sum_{j=0}^{\infty}\gamma (m,j)2^{-j}<0\qquad\qquad
\end{eqnarray*}
(se rappeler la définition de $\sigma _j$). On aura alors la convergence voulue pour les suites $\ln(\alpha _{k+1}/ \alpha _k)$ et $\ln(\beta  _{k+1}/\beta_k)$. Le fait que pour tout entier $k\geq 0$ les quantités $\alpha_k(r)$ et $\beta_k(r)$ tendent monotonement vers zéro lorsque $r$ tend vers zéro, implique par le théorème de la convergence dominée les conclusions $(B)$ et $(C)$ du lemme.\hfill $\Box$
\\
D'après le lemme précédent  on peut  choisir le rayon initial $r_0\in (0,\rho ]$ de telle sorte que l'inégalité $\sum_{k=0}^{\infty}\,\beta _k(r_0)<1/2 + 1/(4\ln 2)$ soit satisfaite. Dans la suite on notera $\alpha _k:=\alpha _k(r_0)$ et $\beta _k:=\beta _k(r_0)$.
On définit maintenant le paramètre $\eta_{k+1}$, qui contrôle la recalibration des calibrations $\omega _k$, $k\geq 0$ au $k$-ième pas du procédé itératif (on désigne avec $\omega _0=\omega$ le choix initial de $\omega$ associée au système différentiel $(S_{\omega})$), de la façon suivante; on définit récursivement en ordre décroissant sur $t=m,...,0$ les matrices 
$$
\eta^{s,t}_{k+1} :=-T_{r(k,m-t)} \,\Big(\omega ^{s,t}_k+\omega ^{s+t+1,-1}_k\wedge\eta^{s,t+1}_{k+1} +(-1)^t\eta^{s-1,t+1}_{k+1}
\wedge\omega ^{s,-1}_k\Big)\in M_{p_{s+t},p_s}({\cal E}_{X}^{0,t}(B_{r(k,m-t)}))     
$$
et on pose $\eta_{k+1} :=(\eta^{s,t}_{k+1})_{s,t}\in {\cal P}(\bar{B}_{r(k,m)}) $. On va justifier ensuite l'estimation $\sup_{z\in B_r} \|\eta^{s,0} _{k+1} (z)\| <\varepsilon $ qui assure l'invertibilité de la matrice $g_{s,k+1}$. On définit alors les matrices
$\omega^{s,t}_{k+1}:=\omega^{s,t}_{k,\;\eta_{k+1}}\in M_{p_{s+k},p_s}({\cal E}_{X}^{0,t+1}(\bar{B}_{r_{k+1} }))$ 
pour tout $t=-1,...,m$. Les constantes $R_k$, $L_k$ qui apparaissent dans la définition des poids $S=(S_k)_{k\geq 0}$  de la norme de Hölder sont définies par les formules:
\begin{eqnarray*}
R_{k+1} :=\max\{\|\omega ^{s,t}_{k,\,I} \|_{r_k,\,\mu }\,/\|\partial^{k+1} \omega^{s,t}_{k,\,I}  \|_{r_k ,\,\mu }  
 \,|\,0\leq s \leq m,\,0\leq t \leq m-s,\,|I|=t+1\}       ,
\\
\\
L_{k+1} :=\max\{ 2^{-k-1} \|g_s(k)^{\pm 1}   \|_{r_k ,\,\mu }\,/\|\partial^{k+1} g_s(k)^{\pm 1}  \|_{r_k ,\,\mu }  \,|\,0\leq s \leq m\}  \qquad\qquad
\end{eqnarray*}
pour tout entier $k\geq 0$. Ici on suppose que 
$\max\{\|\omega ^{s,t}_{k,\,I} \|_{r_k,\,\mu }   
 \,|\,0\leq s \leq m,\,0\leq t \leq m-s,\,|I|=t+1\}>0$, autrement il n'y a rien à prouver.
 Avec ce choix des poids $S_k>0$ on aura que les inégalités
\begin{eqnarray}
&S_{k+1} \|\partial^{k+1} \omega^{\bullet,t}_{k,\,I}  \|_{r_k ,\,\mu } \leq \| \omega^{\bullet,t}_{k,\,I}  \|_{r_k ,\,\mu } ,&\label{15}
\\\nonumber
\\
&S_{k+1} \|\partial^{k+1} g_{\bullet}(k)^{\pm 1}  \|_{r_k ,\,\mu }\leq 2^{-k-1}\|g_{\bullet}(k)^{\pm 1}   \|_{r_k ,\,\mu }&  \label{16}
\end{eqnarray}
seront satisfaites pour tout  entier $k\geq 0$ et $t=0,...,m$. On pose par définition
$$
a_{k} :=\max\{\|\omega ^{s,t}_k\|_{r_k,\,k}\,|\,0\leq s \leq m\,,0\leq t \leq m-s\} ,
$$
$b_k:=H\cdot\sigma_{m,\,k}  ^{-(m+1)\cdot s(k) }\cdot a_k $ et $c\equiv c(\omega _0)$. Avec les notations introduites précédemment on a la proposition suivante.
\begin{prop}\label{3.5}  
Pour tout entier $k\geq 0$ on a les estimations suivantes;
\begin{eqnarray}
&a_{k+1} \leq H\cdot\sigma_{m,\,k} ^{-\nu (m,k)}\cdot a^2_{k}\leq 1 ,& \label{17} 
\\\nonumber
\\
&\|\eta^{s,t}_{k+1} \|_{r_{k+1}  ,\,k+1 } \leq  b_k   <\varepsilon <1/2 ,&\label{18}
\\\nonumber
\\
&\|\omega ^{s,-1} _{k+1}\|_{r_{k+1},k+1 } \leq 4c&\label{19}
\end{eqnarray}
où $H:=\max\{R(C,4c), L(C,4c)\}>0$ et les inégalités $a_k\leq \alpha _k,\;b_k\leq \beta _k$.
\end{prop} 
$Preuve$. On montre les trois estimations précédentes à l'aide d'une récurrence sur $k\geq 0$. Pour $k=0$ on a d'après le lemme $\ref{3.4.2}$, la validité des inégalités
$\alpha _0(r)\leq 1$ et $\beta _0(r)<\varepsilon <1/2$. On est donc en position d'appliquer la proposition $\ref{3.4}$, laquelle assure les inégalités $\eqref{17}$ et $\eqref{18}$ pour $k=0$. On obtient alors l'estimation $\|\eta^{\bullet,0}_1 \|_{r_1,1}<1/2$ laquelle, pour les calculs faits dans la preuve de la proposition $\ref{3.4}$, assure les inégalités 
$\|g_{\bullet,1}^{\pm 1} \|_{r_1,1}<2$ qui assurent donc l'estimation $\eqref{19}$ pour $k=0$, car 
$$
\|\omega ^{s,-1} _1\|_{r_1,1} \leq \|g_{s-1,1}^{-1} \|_{r_1,1}\cdot \|\omega ^{s,-1}\|_{r_1,1}\cdot \|g_{s,1}\|_{r_1,1} .
$$ 
Supposons maintenant par hypothèse récursive que les estimations $\eqref{17}$, $\eqref{18}$ et $\eqref{19}$, sont vraies pour tout $l=0,...,k-1$. L'inégalité  $\eqref{17}$ implique alors l'inégalité 
$$
b_{l+1}\leq P\cdot e^{\gamma (m,l)}\cdot b_l^2
$$
pour les indices en question. On obtient donc les inégalités $a_l\leq \alpha _l\leq 1,\;b_l\leq \beta _l<\varepsilon $. Le fait que par hypothèse inductive on dispose de l'estimation $\|\omega ^{\bullet,-1}_k\|_{r_k,\,k}<4c$ 
permet alors d'appliquer la proposition $\ref{3.4}$ laquelle assure la validité des estimations $\eqref{17}$ et $\eqref{18}$ pour $l=k$. En particulier l'estimation $\eqref{18}$ assure la validité de l'estimation
$$
\|\eta^{\bullet,0}_{k+1} \|_{r_{k+1} ,\,k+1 } \leq \beta _k<\varepsilon <1/2 .
$$
On passe  maintenant à l'estimation des  normes $\|\omega^{s,-1}_{k+1} \|_{r_{k+1} ,\,k+1}$. On a par définition des matrices $\omega ^{s,-1}_{k+1}$ l'estimation  
$$
\|\omega ^{s,-1}_{k+1} \|_{r_{k+1} ,\,k+1} \leq \|g_{s-1}(k+1)^{-1} \|_{r_{k+1} ,\,k+1}\cdot \|\omega ^{s,-1}\|_{r_{k+1} ,\,k+1}\cdot \|g_s(k+1)\|_{r_{k+1} ,\,k+1} .
$$
La condition $\eqref{16}$ sur les poids implique que pour tout $k\geq 0$ on dispose de l'inégalité 
$$
\|g_s(k+1)^{\pm 1}\|_{r_{k+1} ,\,k+1}\leq (1+2^{-k-1})\cdot \|g_{s,k+1}^{\pm 1}\|_{r_{k+1} ,\,k+1}\cdot\|g_s(k)^{\pm 1}\|_{r_k ,\,k} .
$$
L' hypothèse inductive nous permet alors d'effectuer les estimations suivantes pour $k\geq 1$
\begin{eqnarray*}
&\|g_s(k+1)^{\pm 1}  \|_{r_{k+1} ,\,k+1}\leq\displaystyle{ \prod_{j=2}^{k+1}(1+2^{-j})\cdot  \prod_{j=1}^{k+1}\|g_{s,j}^{\pm 1} \|_{r_j,\,j}} \leq &
\\
\\
&\leq\sqrt{e}\cdot\displaystyle{ \prod_{j=1}^{k+1}(1\pm  \|\eta^{s,0}_j\|_{r_j,\,j} )^{\pm 1}  \leq }
 \sqrt{e}\cdot\exp\Big(2(\ln 2)\displaystyle{\sum_{j=1}^{k+1}\|\eta^{s,0}_j\|_{r_j,\,j}}\Big) \leq&
\\
\\
&\leq  \sqrt{e}\cdot \exp\Big(2(\ln 2)\displaystyle{\sum_{j=0}^{\infty}}\,\beta _j\Big)<2&
\end{eqnarray*}
(rappeler le choix initial du rayon $r_0$). On obtient donc l'estimation voulue $\|\omega ^{s,-1}_{k+1} \|_{r_{k+1} ,\,k+1} \leq 4c$, ce qui conclu la preuve des trois estimations $\eqref{17}$, $\eqref{18}$ et $\eqref{19}$ pour $j=k$ et donc pour tout entier positif $k$.\hfill $\Box$

\subsubsection{Preuve de la convergence vers une solution du problème différentiel $(S_{\omega} )$}
Avec les notations introduites précédemment on a la proposition suivante.
\begin{prop}. Les limites  
$$
\eta^{s,t}:=\lim_{k\rightarrow \infty} \eta (k)^{s,t}
$$
$s=0,...,m,\;t=0,...,m-s$ existent en topologie ${\cal C}^{h,\mu}(B_{r_{\infty}})$ pour tout $h\geq 0$ et elles constituent les composantes du paramètre de recalibration $\eta=(\eta^{s,t})_{s,t}\in {\cal P}(B_{r_{\infty}})$, solution du problème différentiel $(S_{\omega} )$.
\end{prop} 
$Preuve$. Nous commençons par prouver l'existence des limites
$$
g_s:={\displaystyle \prod_{j\geq 1 }^{\longrightarrow}}g_{s,j}=\lim _{k\rightarrow \infty}g_s(k)=\I_{p_s}+\lim_{k\rightarrow \infty} \eta (k)^{s,0} 
$$
 en topologie ${\cal C}^{h,\mu}(B_{r_{\infty}})$ pour  $h\geq 0$ quelconque et le fait que les matrices $g_s$ sont inversibles. On aura alors $\eta^{s,0}=g_s-\I_{p_s}$. On déduit immédiatement de la proposition $\ref{3.5}$ les estimations 
$\|\eta^{\bullet,\bullet}_k\|_{r_{\infty},h }<\beta _k<1/2$
et 
$\|g_{\bullet} (k)^{\pm 1}  \|_{r_{\infty},h }<2$ pour tout $k\geq h$. Le fait que 
$g_s(k+1)-g_s(k)=\eta^{s,0}_{k+1}\cdot g_s(k)$ implique les estimations suivantes:
$$
\|g_s(k+1)-g_s(k)\|_{r_{\infty} ,\,h}\leq \|\eta^{s,0}_{k+1}\|_{r_{\infty} ,\,h}\cdot \|g_s(k)\|_{r_{\infty} ,\,h}\leq 2\|\eta^{s,0} _{k+1}\|_{r_{\infty} ,\,h}\leq 2\beta _k
$$
pour tout $k\geq h$. On a alors
\begin{eqnarray*}
\sum_{k=0}^{\infty}\|g_s(k+1)-g_s(k)\|_{r_{\infty} ,\,h}\leq\sum_{k=0}^h\|g_s(k+1)-g_s(k)\|_{r_{\infty} ,\,h}+2\cdot \sum_{k=h+1}^{\infty}\,\beta _k <\infty 
\end{eqnarray*}
et donc l'existence des matrices $g_s\in{\cal C}^{h,\mu}(B_{r_{\infty}} ,M_{p_s,p_s}(\C))$ pour tout $h \geq 0$  telles que 
\begin{eqnarray*}
g_s=\lim_{k\rightarrow \infty}g_s(k)=\I_{p_s}+  \sum_{k=0}^{\infty}(g_s(k+1)-g_s(k))
\end{eqnarray*}
en topologie ${\cal C}^{h,\mu}(B_{r_{\infty}})$. D'autre part l'égalité
\begin{eqnarray*}
g_s(k+1)^{-1} -g_s(k)^{-1} =g_s(k)^{-1}\cdot\sum_{j=1}^{\infty} \,(-1)^j (\eta^{s,0}_{k+1})^j  
\end{eqnarray*}
permet d'effectuer les estimations suivantes:
\begin{eqnarray*}
\|g_s(k+1)^{-1}-g_s(k)^{-1}\|_{r_{\infty} ,\,h}\leq \|g_s(k)^{-1}\|_{r_{\infty} ,\,h}\cdot \sum_{j=1}^{\infty}\|\eta^{s,0}_{k+1} \|_{r_{\infty} ,\,h}^j
\leq 4\|\eta^{s,0} _{k+1} \|_{r_{\infty} ,\,h}\leq 4\beta _k
\end{eqnarray*}
pour tout $k\geq h$, lesquelles  assurent la convergence de la série
\begin{eqnarray*}
\sum_{k=0}^{\infty}\|g_s(k+1)^{-1} -g_s(k)^{-1} \|_{r_{\infty} ,\,h}\leq
\sum_{k=0}^h\|g_s(k+1)^{-1} -g_s(k)^{-1} \|_{r_{\infty} ,\,h} + 4\cdot\sum_{k=h+1}^{\infty} \beta_k<\infty ,
\end{eqnarray*}
ce qui prouve l'existence des matrices $\rho_s\in{\cal C}^{h,\mu}(B_{r_{\infty} } ,M_{p_s,p_s}(\C))$ telles que $\rho _s=\lim_{k\rightarrow \infty}g_s(k)^{-1}$ en topologie ${\cal C}^{h,\mu}(B_{r_{\infty}})$, pour tout $h\geq 0$. On a alors l'égalité  $\I_{p_s}=\lim_{k\rightarrow \infty}g_s(k)\cdot g_s(k)^{-1}=g_s\cdot \rho _s$, qui montre que 
$g_s\in GL(p_s, {\cal E}(B_{r_{\infty}}))$ et $\rho _s=g^{-1}_s$. Montrons maintenant l'existence des limites $\eta^{s,t}$ pour $t\geq 1$. En rappelant l'expression des composantes $\eta(k)^{s,t},\;t\geq 1$, du paramètre $\eta(k)$ introduite dans la sous-section $\ref{3.2}$ et en tenant compte de l'existence de la limite 
$g\in \Gamma (U)$ prouvée précédemment, on déduit qu'il suffit de prouver l'existence des limites
\begin{eqnarray*}
\lim_{k\rightarrow \infty } \sum_{J\in J_k(\rho (\tau) )}\;{\displaystyle \bigwedge_{1\leq r\leq \rho (\tau) }^{\longrightarrow}} 
\;g_{s+\sigma' (\tau,r)}  (j_r-1)\cdot \eta^{s+\sigma (\tau,r),\,\tau_{\rho (\tau)+1-r } }_{j_r} \cdot g_{s+\sigma (\tau,r)} (j_r)^{-1}
\end{eqnarray*}
 $\tau\in \Delta _t,\;t\geq 1$ en topologie ${\cal C}^{h,\mu}(B_{r_{\infty} })$, avec $h\geq 0$ quelconque, pour obtenir l'existence de la limite 
$\eta^{s,t}\in M_{p_{s+t},p_s}({\cal E}_{_X}^{0,t}(B_{r_{\infty} }))$. Il suffit donc de prouver que pour tout $h\geq \rho (\tau)$ la limite 
\begin{eqnarray*}
&\displaystyle{\lim_{k\rightarrow \infty } \sum_{J\in J_k(\rho (\tau) )}\;\prod_{r=1}^{\rho (\tau)}   
\;\|g_{\bullet}  (j_r-1)\cdot \eta^{\bullet,\bullet}_{j_r} \cdot g_{\bullet} (j_r)^{-1}\|_{r_{\infty},h} =}&
\\
\\
&=\displaystyle{ \lim_{k\rightarrow \infty }{\displaystyle\sum_{\scriptstyle l+p=\rho (\tau) 
\atop
\scriptstyle  l,p\geq 0} }\;\sum_{I\in J_h(l)}\;\prod_{r=1}^l 
\;\|g_{\bullet} (i_r-1)\cdot \eta^{\bullet,\bullet}_{i_r} \cdot g_{\bullet} (i_r)^{-1}\|_{r_{\infty},h}\times}&
\\
\\
&\displaystyle{ \times\sum_{I\in J_{h,k} (p)}\;\prod_{r=1}^p 
\;\|g_{\bullet} (j_r-1)\cdot \eta^{\bullet,\bullet}_{j_r} \cdot g_{\bullet} (j_r)^{-1}\|_{r_{\infty},h} }  &
\end{eqnarray*}
est finie. Ici on pose par définition $J_{h,k} (p):=\{J\in \{h+1,...,k\}^p\;|\;j_1<...<j_p \},\;J_h(0)=J_{h,k}(0):=\emptyset$, (rappeler qu'on utilise la convention qui consiste à négliger  les symboles de somme et de produit si l'ensemble d'indices sur lequel on effectue ces opérations est vide). On considère pourtant les estimations suivantes pour $1\leq p\leq \rho (\tau)$
\begin{eqnarray*}
\lim_{k\rightarrow \infty }\sum_{I\in J_{h,k} (p)}\;\prod_{r=1}^p 
\;\|g_{\bullet} (j_r-1)\|_{r_{\infty},h}\cdot \|\eta^{\bullet,\bullet}_{j_r}\|_{r_{\infty},h} \cdot \|g_{\bullet} (j_r)^{-1}\|_{r_{\infty},h}<\qquad\qquad\quad
\\
\\
<4^p\lim_{k\rightarrow \infty }\sum_{J\in J_{h,k} (p)}\;\prod_{r=1}^p \beta_{j_r-1}
<4^p\lim_{k\rightarrow \infty }\sum_{J\in J_{h,k} (p)}\beta_{j_p-1}
=4^p\lim_{k\rightarrow \infty }\sum_{j=h+p}^k \,{j-h-1\choose p-1 }\cdot \beta _{j-1}  <\infty .
\end{eqnarray*}
La dernière limite est finie par le critère du rapport, rappeler en fait que $\lim_{k\rightarrow \infty }\beta _{k+1}/\beta _k=0$, d'après la preuve du lemme 
$\ref{3.4.2}$. On a donc prouvé l'existence du paramètre limite $\eta \in {\cal P}(B_{r_{\infty} })$. Prouvons maintenant qu'il constitue une solution 
(de classe ${\cal C}^{\infty}$) pour le problème différentiel $(S_{\omega } )$. En effet l'existence de la limite $\eta$ en topologie 
${\cal C}^{h,\mu}(B_{r_{\infty} })$, pour $h\geq 1$ implique les égalités
\begin{eqnarray*}
\omega ^{s,t}_{\eta}=\lim_{k\rightarrow \infty}\omega ^{s,t}_{\eta(k)} =\lim_{k\rightarrow \infty}\omega ^{s,t}_k
\end{eqnarray*} 
en topologie ${\cal C}^{h-1,\mu}(B_{r_{\infty} })$. En rappelant l'inégalité $a_k\leq \alpha _k$ obtenue dans la preuve de la proposition $\ref{3.5}$ on obtient
\begin{eqnarray*}
\|\omega ^{s,t}_{\eta} \|_{r_{\infty} ,0}=\lim_{k\rightarrow \infty}\|\omega ^{s,t}_k\|_{r_{\infty} ,0} \leq\lim_{k\rightarrow \infty}a_k=0  ,
\end{eqnarray*} 
ce qui prouve que $\eta \in {\cal P}(B_{r_{\infty} })$ est une solution du  système différentiel
$$
\left  \{
\begin{array}{lr}
\omega^{s,t}_{\eta}=0
\\
s=0,...,m
\\
t=0,...,m-s 
\end{array}                                                                           
\right.
$$
qui n'est rien d'autre que le système différentiel $(S_{\omega } )$.  \hfill $\Box$

\subsection{Cinquième étape: fin de la preuve du théorème $\ref{cdif}$}\label{s3.5} 

L'étape précédente montre qu'on peut se ramener à considérer le diagramme commutatif suivant. 
\begin{diagram}[height=1cm,width=1cm]
&         &                              &                 &0&                                         &0                                                       \\
&         &                              &              &\uTo&                                      &\uTo                                                       \\
0&  \rTo  &(Ker\bar{\partial})_{_{|V} } &  \rTo   &{\cal G}_{|_V}& \rTo^{\bar{\partial}} 
 & {\cal G}_{|_V}\otimes_{_{{\cal E}_{_V}}}{\cal E}^{0,1}_{_V} 
\\
&          &\uTo^{\tilde{\psi}_{|..}}&                     &\uTo_{\tilde{\psi}}&                      &\uTo_{\tilde{\psi}\otimes\I_{(0,1)}}                   \\
0&\rTo  & {\cal O} _{_V}^{\oplus p_0} &  \rTo   &{\cal E}^{\oplus p_0}_{_V}&    \rTo^{\bar{\partial}_{_J }}
 & ({\cal E}^{0,1}_{_V} )^{\oplus p_0}&  
\\
&             &\uTo^{\tilde{\varphi}_{|..}} &                             &\uTo_{\tilde{\varphi}}&               &\uTo_{\tilde{\varphi}\otimes \I_{(0,1)} }
\\
0& \rTo  & {\cal O}
          _{_V}^{\oplus p_1} &  \rTo   &{\cal E}^{\oplus p_1}_{_V}&    \rTo^{\bar{\partial}_{_J }}
 & ({\cal E}^{0,1}_{_V} )^{\oplus p_1} 
\end{diagram} 
\\
avec $V:=B_{r_{\infty}}, \;\tilde{\psi}:=\psi_g$ et $\tilde{\varphi}:=\varphi_{1,\,g} $. Ce diagramme est exact, sauf pour l'instant au niveau des premières flèches verticales à gauche. Pour conclure il nous reste donc à montrer l'exactitude de la suite
$$
{\cal O}_{_V}^{\oplus p_1}\stackrel{\tilde{\varphi}_{|..}}{\longrightarrow}  
{\cal O}_{_V}^{\oplus p_0}\stackrel{\tilde{\psi}_{|..}}{\longrightarrow}
(Ker\bar{\partial})_{_{|V}}\rightarrow 0
$$
On identifie $\tilde{\varphi}=(\tilde{\varphi}_1,...,\tilde{\varphi}_{p_1}  )$, on désigne par
$\tilde{\varphi}_l^k \in{\cal O}_X (V)$ la $k$-ème composante de $\tilde{\varphi}_l$, et on pose par définition
$$
a_{k,x} :={\cal O}_x\tilde{\varphi}^k_{1,x}+...+{\cal O}_x\tilde{\varphi}^k_{p_1,x}\lhd {\cal O}_x
$$
Le fait que $(a_{k,x} \cdot {\cal E}_x)\cap {\cal O}_x=a_{k,x} $ (par définition de fidélité plate de l'anneau ${\cal E}_x$ sur l'anneau ${\cal O}_x$. On peut aussi  déduire l'égalité précédente en utilisant un résultat beaucoup plus simple, i.e la
 fidélité plate de l'anneau des séries formelles $\quad{\cal E}_x /m^{\infty}({\cal E}_x)=\hat{{\cal O}}_x$ en $x$, sur l'anneau ${\cal O}_x$ (voir\cite{mal})) implique la surjectivité du morphisme
$$
\tilde{\varphi}_{|..}:{\cal O}_{_V}^{\oplus p_1}\longrightarrow {\cal R}^{\cal O}(\tilde{\psi})
$$
où ${\cal R}^{\cal O}(\tilde{\psi})={\cal R}^{\cal E}(\tilde{\psi})\cap{\cal O}_{_V}^{\oplus p_0}$ désigne le faisceau des relations holomorphes de $\tilde{\psi}$. On pose par définition
$${\cal F}:=Im(\tilde{\psi}_{|..}: {\cal O}_{_V}^{\oplus p_1}\longrightarrow (Ker\bar{\partial})_{_{|V}})$$
L'exactitude de la suite
$
{\cal O}_{_V}^{\oplus p_1}\stackrel{\tilde{\varphi}_{|..}}{\longrightarrow}  
{\cal O}_{_V}^{\oplus p_0}\stackrel{\tilde{\psi}_{|..}}{\longrightarrow}
{\cal F} \rightarrow 0
$
et la platitude de l'anneau ${\cal E}_x$ sur l'anneau ${\cal O}_x$ impliquent l'existence du diagramme commutatif exact
\begin{diagram}[height=1cm,width=1cm]
{\cal E}^{\oplus p_1}_{_V}& \rTo^{\tilde{\varphi}_1}&{\cal E}^{\oplus p_0}_{_V}&\rTo^{{\tilde{\psi}}_{|..} \otimes \I}&{\cal F}\otimes_{_{{\cal O}_{_V}}}{\cal E}_{_V}&\rTo 0                                                                             \\
\dTo^{\I} &   &\dTo^{\I} &   &\dTo^{\alpha }&                                                                              \\
{\cal E}^{\oplus p_1}_{_V}& \rTo^{\tilde{\varphi}_1}&{\cal E}^{\oplus p_0}_{_V}&\rTo^{\tilde{\psi}}&{\cal G}_{|_V} &\rTo 0
\end{diagram} 
avec $\alpha :\sum_k\tilde{\psi}_{k,y}\otimes_{_{{\cal O}_y}} f_{k,y}\mapsto\sum_k\tilde{\psi}_{k,y} \cdot f_{k,y}\  $. On a alors 
${\cal G}_{|_V}\cong{\cal F}\otimes_{_{{\cal O}_{_V}}}{\cal E}_{_V}.$ Les égalités
\begin{eqnarray*}
\bar{\partial}_{_{\cal F}}\Big(\sum_{k=1}^{p_0}  \tilde{\psi}_{k,y}\otimes_{_{{\cal O}_y}} f_{k,y}\Big) &=&
 \sum_{k=1}^{p_0} \tilde{\psi}_{k,y}\otimes_{_{{\cal O}_y}}\bar{\partial}_{_J } f_{k,y}\cong
\sum_{k=1}^{p_0}\tilde{\psi}_{k,y}\otimes_{_{{\cal E}_y}}\bar{\partial}_{_J } f_{k,y} 
\\
\\
\bar{\partial}\Big(\sum_{k=1}^{p_0} \tilde{\psi}_{k,y}\cdot f_{k,y}\Big) & =&\sum_{k=1}^{p_0} \tilde{\psi}_{k,y}\otimes_{_{{\cal E}_y}}\bar{\partial}_{_J } f_{k,y}  
\end{eqnarray*} 
impliquent la commutativité du diagramme suivant:
\begin{diagram}[height=1cm,width=1cm]
{\cal G}_{|_V}&\rTo^{\bar{\partial}} &{\cal G}_{|_V}\otimes_{_{{\cal E}_{_V}}}{\cal E}^{0,1} _{_V}&             &
\\
\uTo^{\alpha}_{\wr}&                                      &\uTo_{\wr}&          \luTo^{\alpha\otimes\I_{(0,1)}}                          \\    
{\cal F}\otimes_{_{{\cal O}_{_V}}}{\cal E}_{_V}&\rTo^{\bar{\partial}_{_{\cal F}}}&{\cal F}\otimes_{_{{\cal O}_{_V}}}{\cal E}^{0,1} _{_V}&\simeq\ 
&{\cal F}^{\infty} \otimes_{_{{\cal E}_{_V}}}{\cal E}^{0,1} _{_V}
\end{diagram}  
Le fait que le faisceau ${\cal F}$ soit analytique cohérent implique, par les remarques faites dans la section $\ref{S1}$, que ${\cal F}=Ker\bar{\partial}_{_{\cal F}} $. La commutativité du diagramme précédent montre alors que 
${\cal F}=(Ker\bar{\partial})_{_{|V}}$, ce qui démontre le théorème $\ref{cdif}$.\hfill $\Box$

\section{Un résultat d'intégrabilité des connexions sur les faisceaux de ${\cal E}$-modules au dessus d'une variété différentiable}

Le présent travail s'est situé principalement dans le cadre complexe car c'était notre principal intérêt. Toutefois, il est possible de  déduire aussi un résultat d'intégrabilité dans le cas des variétés ${\cal C}^{\infty}$. Considérons en effet $(X,{\cal E}_X)$ une variété ${\cal C}^{\infty}$ (${\cal E}_X\equiv{\cal E}_X(\R)$ représente ici le faisceau des fonctions ${\cal C}^{\infty}$ à valeurs réelles) et 
$D:{\cal G}\longrightarrow  {\cal G} \otimes_{_{{\cal E}_X }}{\cal E}(T^*_X)$ une connexion sur le faisceau de ${\cal E}_X(\K)$-modules ${\cal G}$ où $\K=\R,\;\C$. Si le faisceau des sections parallèles $KerD$ engendre ${\cal G}$ sur ${\cal E}_X(\K)$ alors évidemment $D^2=0$. Le théorème suivant donne une réciproque de ce fait dans un cas particulier.
\begin{theoreme}\label{dcdif} 
Soit $(X,{\cal E}_X)$ une variété différentiable et soit ${\cal G}$ un faisceau de  ${\cal E}_X(\K)$-modules,  muni d'une connexion  
$D:{\cal G}\longrightarrow  {\cal G} \otimes_{_{{\cal E}_{_X}  }}{\cal E}(T^*_X)$ telle que $D^{2}=0$. Si de plus le faisceau ${\cal G}$ admet localement une ${\cal E}(\K)$-résolution de longueur finie, alors le faisceau des sections parallèles $KerD$ engendre sur ${\cal E}_X(\K)$ le faisceau   $\cal G$ qui est le faisceau  des sections d'un système local de coefficients $($le faisceau ${\cal G}$ est donc localement libre$)$ et le complexe $({\cal G} \otimes_{_{{\cal E}_X }}{\cal E}(\Lambda ^{\bullet} T^*_X)\; ; D)$ est une résolution acyclique du faisceau des sections parallèles. On a alors l'isomorphisme fonctoriel  de De Rham-Weil
$H^k(X,KerD)\cong H^k(\Gamma (X,{\cal G} \otimes_{_{ {\cal E}_X }}{\cal E}(\Lambda ^{\bullet} T^*_X))\; ; D)$.
\end{theoreme}
$Preuve$. On commence par substituer formellement dans les étapes de la preuve du théorème $\ref{cdif} $ la connexion $\bar{\partial}$ avec $D$, les opérateurs 
$\bar{\partial}_{_J }$  avec $d$ et l'opérateur de Leray-Koppelman avec l'opérateur d'homotopie de Poincaré 
$$
P_q: {\cal C}^{h} _{q+1}(B_r, M_{k,l}(\K))\longrightarrow  {\cal C}^h_q(B_r, M_{k,l}(\K))
$$
pour $q\geq 0$. Il est élémentaire de vérifier qu'on peut choisir une suite de poids $S=(S_k)_{k\geq 0} \subset (0,+\infty)$ pour les normes ${\cal C}^h$ de telle sorte à obtenir une estimation du type $\|P_q u\|_{r,h}\leq C\cdot \|u\|_{r,h}$, avec $C>0$ indépendante de la régularité $h\geq 0$ de la $q+1$-forme $u$. A condition de restreindre opportunément le rayon $r>0$ on obtient un schéma de convergence rapide considérablement plus simple que  celui expliqué dans la preuve du théorème $\ref{cdif} $. En effet dans le cas en considération  on a pas besoin de restreindre les rayons de la boule pendant les étapes du procédé itératif car on dispose de l'estimation précédente. Les détails de simplification et d'adaptation du procédé itératif  relatif à la preuve du théorème $\ref{cdif} $, au cas en examen sont laissés au lecteur. On obtient en conclusion une 
${\cal E}(\K) $-résolution locale, à partir d'une ${\cal E}(\K) $-résolution initiale, telle que les nouvelles matrices $\hat{\varphi}_j$ (ici on utilise les mêmes notations que dans le théorème $\ref{cdif} $) soient toutes constantes. En particulier le fait que $\hat{\varphi}_1=Cst$ implique que ${\cal G}$ est un faisceau de ${\cal E}(\K)$-modules localement libre. La suite du théorème $\ref{dcdif} $ dérive alors de résultats classiques bien connus du lecteur.\hfill $\Box$
\subsection{Effet gênant de la $\bar{\partial}_{_J }$-stabilité des faisceaux d'ideaux dans le cas des variétés presque-complexes}
Un corollaire du théorème $\ref{dcdif}$ est le suivant.
\begin{corol}\label{effetgenant} 
Soit $(X,J)$ une variété presque-complexe connexe telle que 
$$
T_{X,x} \otimes_{_{\R}}\C=\C\langle\, [\xi ,\eta](x)\;|\; \xi ,\eta \in{\cal E}(T^{0,1}_{X,J})(X)\, \rangle  
$$
pour tout $x\in X$ et soit ${\cal J}\subset {\cal E}(\C) $ un faisceau d'ideaux de fonctions $\ci$ à valeurs complexes sur $X$ admetant des ${\cal E}(\C)$-résolutions locales de longueur finie. Si ${\cal J}$ est $\bar{\partial}_{_J }$-stable alors soit ${\cal J}=0$ soit ${\cal J}={\cal E}(\C) $.
\end{corol} 
L'hypothèse sur le complexifié de l'espace tangent signifie que la structure presque-complexe est ``fortement non-intégrable''. Les variétés presque-complexes pour lesquelles le complexifié de l'espace tangent est engendré ponctuellement par les crochets des champs de vecteurs de type $(0,1)$ seront appellé ``fortement non-intégrables''. Les variétés presque-complexes fortement non-intégrables ont nécéssairement une dimension complexe supérieure à deux. En effet tout les structures presque complexes sont intégrables en dimension complexe un. En dimension complexe deux on a 
$\dim_{_{\C} } \Lambda ^{0,2}_{_J}T_X=1$. Si 
$$
\bar{\tau}_{_J }\in {\cal E}(\Lambda ^{0,2}_{_J}T_X^*\otimes_{_{\C}} T^{1,0}_{X,J})(X)
$$
désigne le tenseur conjugué de la torsion de la structure presque-complexe ($\bar{\tau}_{_J }(\xi ,\eta):= [\xi^{0,1},\eta^{0,1}]^{1,0}$ pour tout $\xi,\eta\in{\cal E}(T_X\otimes_{_{\R}}\C)(U)$, où $U\subset X$ désigne un ouvert quelconque), on a que $\bar{\tau}_{_J }(\Lambda ^{0,2}_{_J}T_X )\subset T^{1,0}_{X,J}$ est contenue strictement dans $T^{1,0}_{X,J}$ si $\dim_{_{\C} } X=2$. Le diagramme suivant 
\begin{diagram}[height=1cm,width=1cm]
&\bar{\tau}_{_J }(x) :\; &{\cal E}(T^{0,1}_{X,J})_x^{\oplus 2} &\rTo        &T_{X,x} \otimes_{_{\R}}\C &\rTo^{\;\;\;\pi^{1,0}_{_J}} &T^{1,0}_{X,J,x}&
\\
&& (\xi ,\eta)                         &\longmapsto &[\xi ,\eta](x)            &\longmapsto           & [\xi ,\eta]^{1,0}(x)&
\end{diagram}
(où $\pi^{1,0} _{_J}:T_X\otimes_{_{\R}}\C\longrightarrow T^{1,0}_{X,J}$ désigne la projection sur le fibré des $(1,0)$-vecteurs tangents) montre alors que pour tout point $x\in X$ le complexifié de l'espace tangent $T_{X,x} \otimes_{_{\R}}\C$ ne peut pas être engendré par les crochets des champs de vecteurs de type $(0,1)$. On donne maintenant un example de variété presque complexe fortement non-intégrable.
\\
{\bf{Example}. } Sur un voisinage ouvert $U\subseteq \C^n,\;n\geq 3$ de l'origine on considére la structure presque-complece dont le fibré  $T^{0,1}_{U,J}$ est engendré par les champs de vecteurs complexes
$$
\displaystyle{\xi _k:=\frac{\partial }{\partial \bar{z}_k}-\sum_{
\scriptstyle   1\leq l \leq n
\atop
\scriptstyle k<r\leq n }\,C_{l,k}^r\, \bar{z}_r\frac{\partial }{\partial z_l } } ,
$$
$k=1,...,n$, où $,C_{l,k}^r\in \C$ sera définie en suite, (on rappelle qu'on utilise la convention qui consiste à négliger les termes d'une somme si l'ensemble des indices sur lesquels on effectue cette opération est vide). On remarque que donner une structure presque-complexe sur  $U\subseteq \C^n$ est équivalent à donner un sous-fibré complexe $F\subset T_{U} \otimes_{_{\R}}\C,\; rg _{_{\C} }F=n$ tel que 
$F\oplus \overline{F} =T_{U} \otimes_{_{\R}}\C$. Bien évidement dans notre cas le voisinage ouvert $U$ de l'origine est choisie suffisement petit pour pouvoir assurer cette dernière condition. Il est facile de vérifier que pour tout $k<t$ on a:
$$
[\xi _k,\xi _t]=\sum_{r=1}^n\,C_{r,k}^t\, \frac{\partial }{\partial z_r }  .
$$
Le fait que $n\geq 3$ permet de choisire $n$-multi-indices $(L_k)_{k=1,...,n},\;L_k=(l_{1,k},l_{2,k}), 1\leq l_{1,k}<l_{2,k}\leq n$ différents. On définit alors les constantes $C_{r,k}^t$ par la règle suivante; 
$C_{r,k}^t=1$ si $(k,t)=(l_{1,r},l_{2,r})$ et $C_{r,k}^t=0$ autrement. On aurà alors 
$$
[\xi _{l_{1,k} },\xi _{l_{2,k} } ]=\frac{\partial }{\partial z_k } . 
$$
Le fait que $[\xi _k,z_k\xi _k]=\xi _k$ montre alors que 
$$
T_{U,0} \otimes_{_{\R}}\C=\C\langle\, [\xi _{l_{1,k} },\xi _{l_{2,k} } ](0),\;[\xi _k,z_k\xi _k](0),\;k=1,...,n\, \rangle  .
$$
Si on désigne par $X\subset U$ le voisinage ouvert de l'origine sur lequel la propriété précédente est vérifié on a que la variété presque-complexe $(X,J)$ est fortement non-intégrable. 
\\
La conclusion du corollaire précedent montre que sur une variété presque-complexe fortement non-intégrable la $\bar{\partial}_{_J }$-stabilité des faisceaux d'ideaux est une propriété très forte qui fait perdre d'intérêt à ces objets! Venons-en maintenant à la preuve du corollaire $\ref{effetgenant}$.
\\
$Preuve$. Le fait que le faisceau  ${\cal J}$ soit $\bar{\partial}_{_J }$-stable combiné avec le fait que le complexifié de l'espace tangent est engendré ponctuellement par les crochets des champs de vecteurs de type $(0,1)$, implique que le faisceau ${\cal J}$ est $d$-stable, en d'autre termes stable par rapport à tous les champs de vecteurs. On peut voir alors l'opérateur $d$ comme une connexion  
$$
d:{\cal J}\longrightarrow  {\cal J} \otimes_{_{{\cal E}_X }}{\cal E}(T^*_X)
$$ 
intégrable sur le faisceau  d'ideaux ${\cal J}$. Une consequence de la preuve du  théorème $\ref{dcdif}$ est que localement il existe un générateur $\psi$ du faisceaux d'ideaux  ${\cal J}$ tel que $d\psi=0$, ce qui permet de conclure. \hfill $\Box$

\vspace{1cm}
\noindent
Nefton Pali
\\
Institut Fourier, UMR 5582, Université Joseph Fourier
\\
BP 74, 38402 St-Martin-d'Hères cedex, France
\\
E-mail: \textit{nefton.pali@ujf-grenoble.fr}
\end{document}